\documentclass{article}
\usepackage{amsmath}
\usepackage{amsfonts}
\usepackage{amsthm}

\newtheorem{thm}{Theorem}[section]
\newtheorem{prop}[thm]{Proposition}
\newtheorem{cor}[thm]{Corollary}

\newtheorem{lem}[thm]{Lemma}

\begin{document}

\title{An Algebraic Description of the Exceptional Isogenies to Orthogonal
Groups}

\author{Shaul Zemel}

\maketitle

\section*{Introduction}

The orthogonal groups of finite dimensional quadratic spaces over fields appear
in many branches of mathematics. They form an infinite family of algebraic
groups, indexed by the dimension of the underlying space. Moreover, if the
field is algebraically closed (or more generally quadratically closed) then
there is only one quadratic space of any given dimension, up to isomorphism.
However, over a general field $\mathbb{F}$ (we consider only characteristic
different from 2 in this manuscript) there may be many different quadratic space
of the same dimension, hence many orthogonal groups of the same dimension. One
way to view these subtleties is as different $\mathbb{F}$-rational structures on
the orthogonal group over the algebraic closure of $\mathbb{F}$. 

In lower dimensions the orthogonal groups become isomorphic (up to finite
index) with members of other families of algebraic groups. These isomorphisms,
especially over $\mathbb{R}$, are very useful in several mathematical fields:
E.g., the relation between $SO^{+}(2,1)$ and $SL_{2}(\mathbb{R})$, between
$SO^{+}(2,2)$ and the product of two copies of $SL_{2}(\mathbb{R})$, and between
$SO^{+}(2,3)$ and $Sp_{4}(\mathbb{R})$ has far-reaching applications in the
theory of modular and automorphic forms---see \cite{[B]}, for example. Now,
these isomorphisms (or isogenies) are relatively simple to describe over a
quadratically closed in dimensions up to 6. Indeed, in the case of algebraically
closed fields and dimensions 4 and 6 such a description appears already in
\cite{[vdW]}, but for more general fields that reference simply refers to the
resulting groups as the $\mathbb{F}$-rational structures mentioned above.
Section III of \cite{[D]} also presents some results for the general orthogonal
groups in these dimensions, but again using this descent method. Moreover, some
aspects of the theory becomes simpler for the general orthogonal group. Hence
in some cases our method is a true refinement of that of \cite{[D]}. In
addition, there is an isogeny in dimension 8 over $\mathbb{R}$, namely signature
$(2,6)$, which is related to a symplectic quaternionic group. This relation
appears in detail in \cite{[SH]}. The phenomenon of triality (see \cite{[KMRT]}
for more details) in the isotropic case may also be viewed as a type of an
exceptional isogeny.

The special orthogonal group of a quadratic space $V$ over $\mathbb{F}$ comes
with a natural central extension, with kernel $\mathbb{F}^{\times}$, called the
Gspin group or the even Clifford group. It has a general definition in terms of
Clifford algebras, namely the set of elements of the even Clifford algebra
of $V$ conjugation by which preserves the embedding of $V$ into the full
Clifford algebra. The spin group is a subgroup of the Gspin group, which maps
onto the spinor norm kernel of $SO(V)$ with kernel $\pm1$. It may also be
described in similar terms using the Clifford algebra. It is these groups, the
spin and Gspin groups, which mainly appear in the exceptional isogenies.
However, the condition of preservation of $V$ under conjugation is not so easy
to verify without delving deep into the multiplicative structure of the Clifford
algebra. This makes the actual structure of the groups thus obtained less
visible.

We use a different, more elementary method in order to determine the spin and
Gspin groups of spaces of dimensions up to 8, thus providing the groups which
are isogenous to the special orthogonal groups in these dimensions. The idea is
simple: We first observe that these groups are invariant under rescaling of the
space, allowing we to choose our space with some extra, useful properties. We
then show how a group (which ends up being the Gspin group) acts naturally on
this (rescaled) space, with kernel $\mathbb{F}^{\times}$. The surjectivity is a
consequence of the Cartan--Dieudonn\'{e} Theorem, since we show how all the
reflections can be realized in an appropriate semi-direct product which maps to
the full orthogonal group, showing that we do get the full Gspin group.
Moreover, in all these cases the spinor norm, which takes values in
$\mathbb{F}^{\times}/(\mathbb{F}^{\times})^{2}$, factors through a map to
$\mathbb{F}^{\times}$, which takes $r$ in the kernel $\mathbb{F}^{\times}$ to
$r^{2}$. This allows us to determine the spin group as the subgroup of the Gspin
group which maps to 1 under this map. In some isotropic cases there are
particular choices or conjugations which one may apply in order for the
realization of the spin and Gspin groups becomes well-known classical groups.

In many cases, especially the isotropic ones, these groups have equivalent
representations, which appear in dimensions 3 and 4 in many places in the
literature. We give a general, simple construction for these equivalent
representations. Dimension 6 is strongly related to the second exterior power of
a 4-dimensional space. We provide the resulting equivalent representations in
this case too. 

The possible complexity of quadratic spaces increases with the size of the
group $\mathbb{F}^{\times}/(\mathbb{F}^{\times})^{2}$. It is thus expedient to
see what are the results for the cases where this group is small, in particular
when it is 1 or 2. In the latter case one must distinguish between two cases,
the Euclidean and the quadratically finite ones, according to whether the field
admits a non-split quaternion algebra or not. We give full details for these
cases.

This manuscript is divided into 12 sections. In Section \ref{Alg} we present the
required notation for central simple algebras over fields, and Section
\ref{Grp} contains the definitions for the groups which we shall encounter.
Section \ref{Dim123} is concerned with dimensions up to 3, and Section
\ref{Dim4} deals with the 4-dimensional case. In Section \ref{Dim6d1} we
examine the case of 6 dimensions and trivial discriminant, and Section
\ref{Dim5} considers 5-dimensional quadratic spaces. Section \ref{Dim6gen} we go
on to general 6-dimensional quadratic spaces, and Section \ref{Wedge} presents
the equivalent representations which arise from the exterior square of
4-dimensional spaces. Section \ref{Dim8id1} takes care of isotropic
8-dimensional quadratic spaces with trivial discriminant, while in Section
\ref{Dim7rd} spaces of dimension 7 which represent their discriminant are
examined. Section \ref{Dim8igen} considers general isotropic spaces of
dimension 8, and Section \ref{ManySq} presents the results for fields in which
the squares have index 1 or 2 in the group of non-zero elements in detail.

\section{Finite Dimensional Algebras \label{Alg}}

Let $\mathbb{F}$ be a field of characteristic different from 2. An
\emph{algebra} over $\mathbb{F}$ is a ring $R$ with identity together with an
embedding of $\mathbb{F}$ into the center of $R$ (taking $1\in\mathbb{F}$ to $1
\in R$). We shall only consider finite-dimensional algebras here, so that we
shall write \emph{algebra} for \emph{finite-dimensional algebra} throught.
Wedderburn's Theorem states that any $\mathbb{F}$-algebra the product of simple
rings $R$, and each such ring admits two maps, the \emph{reduced norm} and the
\emph{reduced trace}, to its center, which is a finite field extension
$\mathbb{K}$ of $\mathbb{F}$. The latter map, denoted $Tr^{R}_{\mathbb{K}}$, is
a homomorphism of additive groups. The former, which we denote
$N^{R}_{\mathbb{K}}$, is multiplicative, yielding a group homomorphism from the
group $R^{\times}$ of invertible elements of $R$ into $\mathbb{K}^{\times}$. The
field extension $\mathbb{K}/\mathbb{F}$ comes with its norm and trace maps
$N^{\mathbb{K}}_{\mathbb{F}}$ and $Tr^{\mathbb{K}}_{\mathbb{F}}$, and the total
norm and trace maps $N^{R}_{\mathbb{F}}$ and $Tr^{R}_{\mathbb{F}}$ are the
compositions $N^{\mathbb{K}}_{\mathbb{F}} \circ N^{R}_{\mathbb{K}}$ and
$Tr^{\mathbb{K}}_{\mathbb{F}} \circ Tr^{R}_{\mathbb{K}}$ respectively. The norm
and trace from the total ring $\prod_{i}R_{i}$ into $\mathbb{F}$ and the
product of the $N^{R_{i}}_{\mathbb{F}}$ and the sum of the
$Tr^{R_{i}}_{\mathbb{F}}$ respectively.

All the tensor products (of vector spaces or algebras) will be over
$\mathbb{F}$, hence the index will be omitted. In case one of the multipliers in
a tensor product is a commutative $\mathbb{F}$-algebra, we shall shorten using
subscripts, such as $R_{\mathbb{E}}$ for $R\otimes\mathbb{E}$.

Many of our algebras will come with an involution, sending $x \in R$ to
$\overline{x}$, which is $\mathbb{F}$-linear and inverts the order of products.
Then $R$ decomposes, as a vector space, as the direct sum of the space $R^{+}$
of \emph{symmetric} elements which are invariant under the involution, and
the space $R^{-}$ of the \emph{anti-symmetric} ones, which are inverted by it
(see Section 2 of \cite{[KMRT]}, and recall that $ch\mathbb{F}\neq2$). In fact,
we shall encounter two cases. In the first case $R$ will be simple with center
$\mathbb{F}$ (hence of dimension $n^{2}$ over $\mathbb{F}$ for some number $n$
called the \emph{degree} of $R$), in which case the involution is said to be
\emph{of the first kind}. Then the involution can be either \emph{orthogonal},
where $R^{+}$ has dimension $\frac{n(n+1)}{2}$ and $R^{-}$ is of dimension
$\frac{n(n-1)}{2}$, or \emph{symplectic}, where the dimensions are interchanged.

The second case is where the center of $R$ is a (\'{e}tale) quadratic
$\mathbb{F}$-algebra $\mathbb{E}$, which is either a field extension (which is
separable by the assumption on $ch\mathbb{F}$) with Galois automorphism $\rho$
(denoted $z \mapsto z^{\rho}$ for $z\in\mathbb{E}$), or
$\mathbb{E}=\mathbb{F}\times\mathbb{F}$, with $\mathbb{F}$ embedded diagonally
and $\rho$ interchanging the two coordinates. Note that for $z$ in a quadratic
$\mathbb{F}$-algebra $\mathbb{E}$ we have
$z^{\rho}=Tr^{\mathbb{E}}_{\mathbb{F}}(z)-z$. We shall encounter only such
algebras which come from central simple algebras over $\mathbb{F}$, i.e., those
$R$ over $\mathbb{E}$ such that there exists a central simple algebra $S$ over
$\mathbb{F}$, with involution $x\mapsto\overline{x}$, such that $R \cong
S_{\mathbb{E}}$. Then we write $y^{\rho}$ for the image of $y \in R \cong
S_{\mathbb{E}}$ under $Id_{S}\otimes\rho$, and the involution in question, which
is \emph{of the second type} or \emph{unitary}, is
$y\mapsto\overline{y}^{\rho}$. In this case $R$ has dimension $2n^{2}$ over
$\mathbb{F}$, and the spaces $R^{+}$ and $R^{-}$ are both of dimension $n^{2}$
over $\mathbb{F}$ (and are \emph{not} vector spaces over $\mathbb{E}$). In case
$\mathbb{E}=\mathbb{F}\times\mathbb{F}$ we have $R=S \times S$, and
$\overline{(s,t)}^{\rho}=(\overline{t},\overline{s})$.

For a finite-dimensional $\mathbb{F}$-algebra $R$, let $M_{n}(R)$ be the ring of
$n \times n$ matrices over $R$. In case $R$ is commutative, the reduced norm and
trace from $M_{n}(R)$ into $R$ are just the matrix determinant and matrix trace
respectively. If furthermore $R=\mathbb{E}$ is a quadratic $\mathbb{F}$-algebra,
we shorten $(M^{t})^{\rho}=(M^{\rho})^{t}$ to just $M^{t\rho}$, where $M^{t}$
denotes the matrix which is the transpose of $M$. Similarly, for
$z\in\mathbb{E}^{\times}$ we write $z^{-\rho}$ for
$(z^{-1})^{\rho}=(z^{\rho})^{-1}$.

\smallskip

A \emph{quaternion algebra} $B$ over $\mathbb{F}$ is a central simple
$\mathbb{F}$-algebra of degree 2. It comes with a natural (symplectic)
involution, called the \emph{main involution}, which we denote by
$\iota:x\mapsto\overline{x}=Tr^{B}_{\mathbb{F}}(x)-x$ (or sometimes $\iota_{B}$
where the quaternion algebra will not be clear from the context). This is the
only symplectic involution on $B$---see Proposition 2.21 of \cite{[KMRT]}.
Since $ch\mathbb{F}\neq2$, every quaternion algebra is generated by two
anti-commuting (traceless) elements with squares in $\mathbb{F}$. The algebra in
which the squares of these elements are $\alpha$ and $\beta$ respectively will
be denoted $\big(\frac{\alpha,\beta}{\mathbb{F}}\big)$. We may multiply each
generator by an element of $\mathbb{F}^{\times}$, so that muliplying $\alpha$
or $\beta$ by squares yields an isomorphic quaternion algebra. If $\mathbb{E}$
is a quadratic $\mathbb{F}$-algebra and $B$ is a quaternion algebra then the
norms $N^{\mathbb{E}}_{\mathbb{F}}$ and $N^{B}_{\mathbb{F}}$ are quadratic
functions, and we have
\begin{lem}
The equalities
$N^{\mathbb{E}}_{\mathbb{F}}(z+w)=N^{\mathbb{E}}_{\mathbb{F}}(z)+N^{\mathbb{E}}
_{\mathbb{F}}(w)+Tr^{\mathbb{E}}_{\mathbb{F}}(zw^{\rho})$ and
$N^{B}_{\mathbb{F}}(x+y)=N^{B}_{\mathbb{F}}(x)+N^{B}_{\mathbb{F}}(y)+Tr^{B}_{
\mathbb{F}}(x\overline{y})$ hold for every $z$ and $w$ in $\mathbb{E}$ and $x$
and $y$ in $B$. \label{BEpol}
\end{lem}

\begin{proof}
The equalities $N^{\mathbb{E}}_{\mathbb{F}}(t)=tt^{\rho}$ and
$Tr^{\mathbb{E}}_{\mathbb{F}}(t)=t+t^{\rho}$ hold for every $t\in\mathbb{E}$,
and we also have $N^{B}_{\mathbb{F}}(s)=s\overline{s}$ and
$Tr^{B}_{\mathbb{F}}(s)=s+\overline{s}$ for every $s \in B$. The lemma follows
directly from these equalities.
\end{proof}

For such $\mathbb{E}$ and $B$ we denote $\mathbb{E}_{0}$ and $B_{0}$ the spaces
of traceless elements in $\mathbb{E}$ and $B$. These are the spaces
$\mathbb{E}^{-}$ and $B^{-}$ with respect to the involutions $\rho$ and $\iota$
respectively, and they have respective dimensions 1 and 3 (as $\rho$ is unitary
and $\iota$ is symplectic).

A quaternion algebra $B$ over $\mathbb{F}$ either \emph{splits}, i.e., it is
isomorphic to $M_{2}(\mathbb{F})$, or a division algebra. In the former case
$Tr^{B}_{\mathbb{F}}$ is the matrix trace, $N^{B}_{\mathbb{F}}$ is the
determinant, and $\iota_{B}$ is the adjoint involution $\binom{a\ \ b}{c\ \
d}\mapsto\binom{\ \ d\ \ -b}{-c\ \ \ \ a}$. A \emph{splitting field} of $B$ is
an extension $\mathbb{K}$ of $\mathbb{F}$ such that the quaternion algebra
$B_{\mathbb{K}}$ over $\mathbb{K}$ splits. A quadratic extension $\mathbb{K}$ of
$\mathbb{F}$, with Galois automorphism $\sigma$, is a splitting field of $B$ if
and only if it can be embedded into $B$. By choosing an embedding, the subset of
$B$ which anti-commutes with $\mathbb{K}$ (namely $x \in B$ such that
$xz=z^{\sigma}x$ for all $z\in\mathbb{K}$) form a 1-dimensional subspace over
$\mathbb{K}$, which is contained in $B_{0}$. The choice of the square $\delta$
of an invertible element there (which is a representative of a class in
$\mathbb{F}^{\times}/N^{\mathbb{K}}_{\mathbb{F}}(\mathbb{K}^{\times})$) yields
an embedding of $B$ into $M_{2}(\mathbb{K})$ in which $z+wj \in B$ (with
$z$ and $w$ from $\mathbb{K}$) is taken to $\binom{\ z\ \ \ \delta
w}{w^{\sigma}\ \ z^{\sigma}}$. We denote the image of this algebra as
$(\mathbb{K},\sigma,\delta)$. For this algebra we shall use
\begin{lem}
The action of $\sigma=Id_{B}\otimes\sigma$ on
$B_{\mathbb{K}}=M_{2}(\mathbb{K})$ for $B=(\mathbb{K},\sigma,\delta)$ is
defined by $\binom{a\ \ b}{c\ \ d}\mapsto\binom{0\ \ \delta}{1\ \
0}\binom{a^{\sigma}\ \ b^{\sigma}}{c^{\sigma}\ \ d^{\sigma}}\binom{0\ \
\delta}{1\ \ 0}/\delta$. \label{KsplitB}
\end{lem}

\begin{proof}
As $B_{\mathbb{K}}=M_{2}(\mathbb{K})$, it suffices to show that
$(\mathbb{K},\sigma,\delta)$ is the set of matrices which are stable under the
operation which is asserted to be $Id_{B}\otimes\sigma$. As this operation
takes $M=\binom{a\ \ b}{c\ \ d}$ to $\binom{\ d^{\sigma}\ \ \ \ \delta
c^{\sigma}}{b^{\sigma}/\delta\ \ \ a^{\sigma}\ }$, we find that $M$ is
invariant under this operation if and only if $d=a^{\sigma}$ and $b=\delta
c^{\sigma}$. As these conditions indeed characterize
$(\mathbb{K},\sigma,\delta)$, this proves the lemma.
\end{proof}

For the split algebra $B=M_{2}(\mathbb{F})$ we have the following
\begin{lem}
For any $g \in M_{2}(\mathbb{F})$, conjugation by $\binom{0\ \ -1}{1\ \ \ \ 0}$
takes $\overline{g}$ and $g^{t}$ to one another. \label{Sadjt}
\end{lem}

\begin{proof}
When $g=\binom{a\ \ b}{c\ \ d}$ we have $g^{t}=\binom{a\ \ c}{b\ \ d}$ and
$\overline{g}=\binom{\ \ d\ \ -b}{-c\ \ \ \ a}$. The result now follows from a
simple calculation.
\end{proof}

Lemma \ref{Sadjt} will allow us to obtain equivalent models for our spin and
Gspin groups, which are also used in the literature. Note that it relates a
symplectic involution on $M_{2}(\mathbb{F})$ with an orthogonal one.

Another type of algebras which we shall encounter are \emph{bi-quaternion
algebras}, which are simple algebras $A$ of degree 4 which may be presented as
$A \cong B \otimes C$ where $B$ and $C$ are quaternion algebras. Given such a
presentation, there is an involution
$\iota_{B}\otimes\iota_{C}:x\mapsto\overline{x}$ on $A$, which is orthogonal.
However, $A$ may be presented as the tensor product of two quaternion algebras
in many ways, each giving a different orthogonal involution. Moreover, not all
the orthogonal involutions on $A$ may be obtained in this way, and $A$ admits
also symplectic involutions---see Proposition 2.7 of \cite{[KMRT]}. However, we
shall use the notation $\overline{x}$ for a bi-quaternion algebra only when the
presentation as $B \otimes C$ is clear from the context.

\smallskip

We shall be needing subgroups of the groups of the form $R^{\times}$ arising
from simple $\mathbb{F}$-algebras $R$, which are defined by the norm. If
$\mathbb{K}$ is the center of $R$ and $H\subseteq\mathbb{E}^{\times}$, we shall
denote $R^{H}$ the subgroup of $\mathbb{E}^{\times}$ consisting of those
elements $x \in R^{\times}$ such that $N^{R}_{\mathbb{E}}(x) \in H$. We shall
extend this notation to algebras of the form $R=S \times S$ for a central
simple $\mathbb{F}$-algebra $S$, with subgroups of
$\mathbb{F}^{\times}\times\mathbb{F}^{\times}$. If $\mathbb{E}$ is any
commutative $\mathbb{F}$-algebra then $\mathbb{E}^{H}$ is defined similarly.
Note that for an algebra of the form $S_{\mathbb{E}}$ with $S$ central simple
over $\mathbb{F}$, we shall use the norm to the center $\mathbb{E}$ and not to
$\mathbb{F}$ for defining $S_{\mathbb{E}}^{H}$.

The group of invertible matrices in $M_{n}(R)$ will be denoted $GL_{n}(R)$ also
when $R$ is not commutative. If $R$ comes with an involution
$x\mapsto\overline{x}$, then $M\mapsto\overline{M}^{t}$ is an involution on
$M_{n}(R)$ of the same type as the involution on $R$. Note that if $R$ is not
commutative then $M\mapsto\overline{M}$ and $M \mapsto M^{t}$ do not behave
well with respect to products. If $R$ is simple then $M_{n}(R)$ is also simple,
with the same center $\mathbb{K}$, and the group $M_{n}(R)^{H}$ for a subgroup
$H\subseteq\mathbb{K}^{\times}$ will be denoted $GL_{n}^{H}(R)$. In case
$H=\{1\}$ we shall use just the superscript 1 (with no brackets), and where
$R=\mathbb{E}$ is a field extension of $\mathbb{F}$, we shall write
$SL_{n}(\mathbb{E})$ for $GL_{n}^{1}(\mathbb{E})$.

\section{Orthogonal and Other Groups \label{Grp}}

We shall be considering quadratic spaces over $\mathbb{F}$. As
$ch\mathbb{F}\neq2$, this is equivalent to spaces endowed with a symmetric
bilinear form, so that we use the bilinear and quadratic forms interchangeably.
All the spaces we consider will of (positive) finite dimension and
non-degenerate, and these assumptions will be made even when not stated
explicitly. Many of our vector spaces will be subsets of $\mathbb{F}$-algebras,
so that we write the pairing (or product) of two vectors $v$ and $w$ of a
quadratic space $v$ as $\langle v,w \rangle$. Moreover, the number $\langle v,v
\rangle$ will be written $|v|^{2}$ (in order to distinguish it from $v^{2}$ in
the algebra involved), and will be called the called the \emph{vector norm} (and
not just norm) of $v$. In case confusion may arise as to which vector space is
considered, we may write also $\langle v,w \rangle_{V}$ and $|v|_{V}^{2}$. We
have
\begin{lem}
the equality $|v+w|^{2}=|v|^{2}+|w|^{2}+2\langle v,w \rangle$ holds for any $v$
and $w$ in $V$. \label{Vpol}
\end{lem}

\begin{proof}
The lemma follows directly from the definition of the vector norm and the
symmetry of the bilinear form.
\end{proof}

While no ``absolute value'' $|v|$ exists in general, we shall write the $m$th
power of $|v|^{2}$ as $|v|^{2m}$. Two vectors $v$ and $w$ are said to be
\emph{orthogonal} or \emph{perpendicular} if $\langle v,w \rangle=0$, and
$v^{\perp}$ denotes the (1-codimensional) subspace of elements of $V$ which are
perpendicular to $v$. A vector $0 \neq v \in V$ is called \emph{isotropic} is
$|v|^{2}=0$, and \emph{anisotropic} otherwise. A quadratic space $V$ is called
\emph{isotropic} if it contains (non-zero) isotropic vectors, and
\emph{anisotropic} otherwise. An \emph{orthogonal basis} is a basis consisting
of vectors which are all orthogonal to one another (hence they must all be
anisotropic), and every quadratic space admits such a basis since
$ch\mathbb{F}\neq2$. A \emph{rescaling} of a quadratic space $V$ is the same
quadratic space but with all the pairings and vector norms multiplied by a
global scalar from $\mathbb{F}^{\times}$. The \emph{determinant} of a quadratic
space $V$ is defined as the determinant of a Gram matrix representing the
bilinear form of $V$ in some basis (which reduces to the product of the vector
norms of an orthogonal basis) in
$\mathbb{F}^{\times}/(\mathbb{F}^{\times})^{2}$. However, it turns out more
useful to consider the \emph{discriminant} of $V$, which is the determinant
multiplied by $(-1)^{n(n-1)/2}$ where $n$ is the dimension of $V$. By some abuse
of notation, we shall sometime treat the discriminant as an actual
representative from $\mathbb{F}^{\times}$ mapping to the appropriate class in
$\mathbb{F}^{\times}/(\mathbb{F}^{\times})^{2}$, but this will always be
independent of the representative chosen. Note that the discriminant is
invariant under rescaling if the dimension is even, but it is multiplied by the
rescaling factor when the dimension is odd.

\smallskip

Given a quadratic space $V$, we define the \emph{orthogonal group} $O(V)$ to be
the group of linear transformations of $V$ which preserves the bilinear form.
Given an anisotropic vector $v \in V$, the map which inverts $v$ and leaves the
space $v^{\perp}$ invariant belongs to $O(V)$. We call this map the
\emph{reflection inverting $v$}. The only property of $O(V)$ which we shall use
here is the \emph{Cartan--Dieudonn\'{e} Theorem}, namely
\begin{prop}
The group $O(V)$ is generated by reflections. \label{CDT}
\end{prop}

\begin{proof}
For a simple proof, see e.g., Corollary 4.3 of \cite{[MH]}.
\end{proof}

The determinant is a surjective homomorphism $O(V)\to\{\pm1\}$ (as reflections
have determinant $-1$) and the kernel, the \emph{special orthogonal group}, is
denoted $SO(V)$. It consists of those transformations which can be written as
the product of an \emph{even} number of reflections. There is a homomorphism
$O(V)\to\mathbb{F}^{\times}/(\mathbb{F}^{\times})^{2}$, called the \emph{spinor
norm}, which takes the reflection inverting an anisotropic vector $v$ to the
image of $|v|^{2}$ in $\mathbb{F}^{\times}/(\mathbb{F}^{\times})^{2}$. As with
the discriminant, we may sometimes say that an element of $O(V)$ has spinor norm
$t\in\mathbb{F}^{\times}$, meaning that its spinor norm is
$t(\mathbb{F}^{\times})^{2}\in\mathbb{F}^{\times}/(\mathbb{F}^{\times})^{2}$.
Note the rescaling the bilinear form leaves the spinor norms of elements of
$SO(V)$ invariant, but multiplies those of elements in
$O(V) \setminus SO(V)$ by the rescaling factor. Hence for $SO(V)$ the spinor
norm is well-defined also when we consider $V$ up to rescaling. The subgroup of
$SO(V)$ consisting of elements having spinor norm 1 (i.e., a square) is denoted
$SO^{1}(V)$. Note that the global inversion $-Id_{V}$ always has spinor norm
which equals the discriminant of the space, in correspondence with $-Id_{V}$
being in $SO(V)$ (hence having a spinor norm which is invariant under
rescalings) if and only if $V$ has even dimension.

We shall define the \emph{spin group} of $V$ to be a double cover of
$SO^{1}(V)$. The \emph{Gspin group}, or the \emph{even Clifford group}, of $V$
is defined as a group mapping onto $SO(V)$ with with kernel
$\mathbb{F}^{\times}$. We wish to construct these groups, in low dimensions,
without needing to investigate the Clifford algebra of $V$. This becomes much
simpler after some normalization by rescaling. Therefore we do not consider
groups like the pin group, $O^{1}(V)$, and the full Clifford group, which map
onto subgroups of $O(V)$ which are not contained in $SO(V)$, as they are not
invariant under rescaling.

\smallskip

Let $\mathbb{E}$ be a quadratic extension of $\mathbb{F}$, with Galois
automorphism $\rho$. A \emph{unitary space} over $\mathbb{E}$ (with respect to
$\rho$ is a vector space (again finite-dimensional and non-trivial) vector space
over $\mathbb{E}$ with a (non-degenerate) Hermitian sesqui-linear form, where
the conjugation is defined using $\rho$. Unitary spaces may be defined in terms
Hermitian Gram matrices, which may always be reduced (by the choice of an
appropriate basis) to regular diagonal matrices over the fixed field
$\mathbb{F}$ of $\rho$. A unitary space also has a determinant (and a
discriminant), which are similarly defined and lie in
$\mathbb{F}^{\times}/N^{\mathbb{E}}_{\mathbb{F}}(\mathbb{E}^{\times})$. We
shall allow ourselves the same abuse of notation for these unitary determinants
and discriminants. We shall present unitary spaces only through representing
Gram matrices (which are Hermitian and regular), and the \emph{unitary group} of
such a matrix $M$, denoted $U_{\mathbb{E},\rho}(M)$, is the group of linear
transformations of a unitary spaces which preserve the sesqui-linear form
defined by $M$, namely those $g \in GL_{n}(\mathbb{E})$ such that
$gMg^{t\rho}=M$. Elements of $U_{\mathbb{E},\rho}(M)$ have determinants in
$\mathbb{E}^{1}$, and we define $SU_{\mathbb{E},\rho}(M)$ to be the subgroup of
unitary matrices whose determinant is 1. Matrices which multiply the
sesqui-linear form (hence $M$) by a scalar from $\mathbb{E}^{\times}$, which
must then be in $\mathbb{F}^{\times}$, form the \emph{general unitary group}
$GU_{\mathbb{E},\rho}(M)$. Now, if $g \in GU_{\mathbb{E},\rho}(M)$ multiplies
the sesqui-linear form (hence $M$) by a scalar $s=s(g)\in\mathbb{F}^{\times}$,
and the space has dimension $n$, then we have $N^{\mathbb{E}}_{\mathbb{F}}(\det
g)=s(g)^{n}$. In case $n$ is even, we define the group
$GSU_{\mathbb{E},\rho}(M)$ consisting of elements $g$ of the latter group such
that $\det g=s(g)^{n/2}$ (note that this is not equivalent to the condition that
$\det g\in\mathbb{F}^{\times}$---see Lemma \ref{ind2int} and Corollary
\ref{iso6gen} below for an example with $n=4$). It follows that
$SU_{\mathbb{E},\rho}(M)=U_{\mathbb{E},\rho}(M) \cap GSU_{\mathbb{E},\rho}(M)$.
All these groups are invariant under rescaling of $M$ by an element of
$\mathbb{F}^{\times}$. In small dimension we may relate the (general) unitary
groups to other groups, as is seen in the following
\begin{lem}
If $M$ is 1-dimensional then $GU_{\mathbb{E},\rho}(M)$,
$U_{\mathbb{E},\rho}(M)$, and $SU_{\mathbb{E},\rho}(M)$ are
$\mathbb{E}^{\times}$, $\mathbb{E}^{1}$, and $\{1\}$ respectively. In the case
of 2-dimensional, $GSU_{\mathbb{E},\rho}(M)$ is conjugate to $B^{\times}$ for
some quaternion algebra $B$ over $\mathbb{F}$ which is split over $\mathbb{K}$,
and $SU_{\mathbb{E},\rho}(M)$ is conjugate to $B^{1}$. \label{uniEB}
\end{lem}

\begin{proof}
In case $M$ is just a scalar (which may be taken to be 1), the unitary
relations on $z \in GL_{1}(\mathbb{E})=\mathbb{E}^{\times}$ are just
$zz^{\rho}\in\mathbb{F}^{\times}$ (which poses no further restriction on $z$),
$zz^{\rho}=1$, and $z=1$ respectively. In the 2-dimensional case we may take
$M$ to be diagonal (this change might impose some conjugacy relation), and
after rescaling we may assume that $M=\binom{-\varepsilon\ \ 0}{\ \ 0\ \ 1}$
where $\varepsilon$ represents the discriminant of the unitary space. Now,
multiplying the defining relation of $g \in GSU_{\mathbb{E},\rho}(M)$ by
$\binom{0\ \ -1}{1\ \ \ \ 0}$ from the right and using Lemma \ref{Sadjt}
transforms this relation to $g\binom{0\ \ \varepsilon}{1\ \
0}\overline{g}^{\rho}=\det g\binom{0\ \ \varepsilon}{1\ \ 0}$. As $\det
g=g\overline{g}$, the latter relation shows that $\overline{g}$ is invariant
under the relation from Lemma \ref{KsplitB}, implying that $\overline{g}$ lies
in $B=(\mathbb{E},\rho,\varepsilon)$. As the latter algebra is closed under the
adjoint involution (which restricts to its main involution), we find that $g
\in B^{\times}$ as well. As $SU_{\mathbb{E},\rho}(M)$ is the subgroup of
determinant 1 elements in $GSU_{\mathbb{E},\rho}(M)$, it is taken to $B^{1}$ in
this map. This proves the lemma.
\end{proof}

For a subgroup $H$ of $\mathbb{F}^{\times}$, we write
$GU_{\mathbb{E},\rho}^{H}(M)$, as well as $GSU_{\mathbb{E},\rho}^{H}(M)$ of $n$
is even, for the subgroup of the appropriate groups consisting of matrices $g$
whose multiplier $t(g)$ lies in $H$. Thus
$GSU_{\mathbb{E},\rho}^{H}(M)=GSU_{\mathbb{E},\rho}(M) \cap
GU_{\mathbb{E},\rho}^{H}(M)$. The same argument as in Lemma \ref{uniEB} shows
that if $M$ is 1-dimensional then $GU_{\mathbb{E},\rho}^{H}(M)$ is
$\mathbb{E}^{H}$, while for 2-dimensional $M$ the group
$GSU_{\mathbb{E},\rho}^{H}(M)$ is conjugate to $B^{H}$ for some quaternion
algebra $B$ over $\mathbb{F}$.

\smallskip

The (classical) \emph{symplectic group} $Sp_{2n}(\mathbb{F})$ is the group
consisting of those matrices $g \in GL_{2n}(\mathbb{F})$ such that
$g\binom{0\ \ -I}{I\ \ \ \ 0}g^{t}=\binom{0\ \ -I}{I\ \ \ \ 0}$. More generally,
a symplectic group is the group of linear transformations of a vector space
(which must be of even dimension) which preserve a non-degenerate anti-symmetric
bilinear form on it, but every such group is isomorphic (or conjugate) to the
classical one. The \emph{general symplectic group} $GSp_{2n}(\mathbb{F})$
consists of those matrices whose action multiplies $\binom{0\ \ -I}{I\ \ \ \ 0}$
by a scalar. If we consider a quadratic extension $\mathbb{E}$ of $\mathbb{F}$
with Galois automorphism $\rho$, then preserving an anti-Hermitian matrix via
$g:M \mapsto gMg^{t\rho}$ is the same as preserving the Hermitian matrix which
is obtained from $M$ through multiplication by a scalar from $\mathbb{E}_{0}$,
so that no new groups are obtained in this way. On the other hand, using a
central simple algebra $R$ with an involution $x\mapsto\overline{x}$ of the
first kind, one defines further types of symplectic groups. We have the
operation of $GL_{n}(R)$ on $M_{n}(R)$ via $g:M \mapsto gM\overline{g}^{t}$, and
any matrix (Hermitian, anti-Hermitian, or neither) may be used to define such a
group. Note that $g:M \mapsto gMg^{t}$ and $g:M \mapsto gM\overline{g}$ may not
be used here, as they do not define actions of $M_{n}(R)$. Now, given any $M
\in GL_{n}(R)$, we let $Sp_{R}(M)$ is the group of elements $g \in GL_{n}(R)$
which preserve $M$, and $GSp_{R}(M)$ consists of those matrices whose action
multiplies $M$ by a scalar from $\mathbb{F}^{\times}$. Note that rescaling $M$
by a factor from $\mathbb{F}^{\times}$ still does not change these groups. If
$B$ is a quaternion algebra (with its main involution) and $M$ is Hermitian, we
may choose a basis for our space such that $M$ is diagonal, with entries from
$\mathbb{F}$. In this case we just have
\begin{lem}
If $M$ has dimension 1 then $GSp_{B}(M)=B^{\times}$ and $Sp_{B}(M)=B^{1}$.
\label{Sp1B}
\end{lem}

\begin{proof}
Indeed, $M$ is just a scalar from $\mathbb{F}^{\times}$, which may be taken to
be 1. Hence the $GSp$ relation just states that
$x\overline{x}\in\mathbb{F}^{\times}$ (which poses no restriction on the element
$x \in GL_{1}(B)=B^{\times}$), and the $Sp$ condition means $x\overline{x}=1$.
This proves the lemma.
\end{proof}

In resemblance with the classical case, we shall use $Sp_{2n}(R)$ and
$GSp_{2n}(R)$ for the case where $M$ is the anti-Hermitian matrix $\binom{0\ \
-I}{I\ \ \ \ 0}$. As usual, for a subgroup $H\subseteq\mathbb{F}^{\times}$ we
define $GSp_{R}^{H}(M)$ and $GSp_{2n}^{H}(R)$ to be the subgroup of $GSp_{R}(M)$
and $GSp_{2n}(R)$ in which the multiplier comes from $H$. The proof of Lemma
\ref{Sp1B} shows that if $B$ is a quaternion algebra with its main involution
and $M$ is 1-dimensional and Hermitian then the first group is just $B^{H}$.

\section{Dimension $\leq3$ \label{Dim123}}

In dimension 1 we have only one quadratic space (up to isomorphism), namely
$\mathbb{F}$ itself, and the bilinear form is determined by the norm of 1 (which
is most naturally normalized to be 1). Proposition \ref{CDT} shows that
$O(\mathbb{F})$ is generated by the only reflection $-Id_{\mathbb{F}}$, so that
it equals $\{\pm1\}$ and $SO(\mathbb{F})=\{1\}$. The spinor norm is just 1 on
$SO(\mathbb{F})$ (i.e., $SO^{1}(\mathbb{F})=SO(\mathbb{F})$). The Gspin group is
thus $\mathbb{F}^{\times}$, and the spin group is $\{\pm1\}$, both mapping to
the trivial group $SO(\mathbb{F})=\{1\}$.

\medskip

For dimension 2 we define $\mathbb{E}$ to be $\mathbb{F}(\sqrt{d})$, where $d$
is the discriminant of the space. Our space is described by the following
\begin{lem}
Any 2-dimensional quadratic space of discriminant $d$ is isometric to
$\mathbb{E}$, with a rescaling of the quadratic form
$N^{\mathbb{E}}_{\mathbb{F}}$. Rescaled appropriately, we get $2\langle z,w
\rangle=Tr^{\mathbb{E}}_{\mathbb{F}}(zw^{\rho})$ for $z$ and $w$ in
$\mathbb{E}$. \label{sp2}
\end{lem}

\begin{proof}
The second equality for $\mathbb{E}$ with
$|z|^{2}=N^{\mathbb{E}}_{\mathbb{F}}(z)$ follows from Lemmas \ref{BEpol} and
\ref{Vpol}. Hence $1\in\mathbb{E}$ satisfies $|1|^{2}=1$, an element from
$\mathbb{E}_{0}$ has vector norm $-d$ (up to $(\mathbb{F}^{\times})^{2}$), and
they are orthogonal. Now, rescaling our original space space such that some
anisotropic vector has vector norm 1, we find the orthogonal complement must be
spanned by a vector whose vector norm is the determinant $-d$, just like in
$\mathbb{E}$. This proves the lemma.
\end{proof}

Multiplication from $\mathbb{E}^{1}$ preserves the vector norms, which defines a
map $\mathbb{E}^{1} \to O(\mathbb{E})$, which is clearly injective. However, in
order to define the spin and Gspin group and be in the same spirit as the
constructions for higher dimensions, we shall use
\begin{lem}
The action $g:z \mapsto gzg^{-\rho}$ defines a map $\mathbb{E}^{\times} \to
O(\mathbb{E})$, with kernel $\mathbb{F}^{\times}$. The semi-direct product of
$Gal^{\mathbb{E}}_{\mathbb{F}}=\{Id_{\mathbb{E}},\rho\}$ with
$\mathbb{E}^{\times}$ also maps to $O(\mathbb{E})$. \label{ac2}
\end{lem}

\begin{proof}
As $N^{\mathbb{E}}_{\mathbb{F}}(g^{\rho})=N^{\mathbb{E}}_{\mathbb{F}}(g)$ and
the norm is multiplicative, we have the equalities $|gzg^{-\rho}|^{2}=|z|^{2}$
as well as $|z^{\rho}|^{2}=|z|^{2}$. Hence both $\mathbb{E}^{\times}$ and $\rho$
map to $O(\mathbb{E})$. The kernel of the map from $\mathbb{E}^{\times}$
consists of those elements of $\mathbb{E}^{\times}$ such that $g=g^{\rho}$,
which is $\mathbb{F}^{\times}$. The equality
$(gzg^{-\rho})^{\rho}=g^{\rho}z^{\rho}g^{-1}$ shows that the map from the
semi-direct product is also a homomorphism. This proves the lemma.
\end{proof}

In fact, the map from $\mathbb{E}^{\times}$ defined in Lemma \ref{ac2} and the
map defined above it have the same image, by Hilbert's Theorem 90. The next step
is
\begin{lem}
Fix some $0 \neq h\in\mathbb{E}^{0}$. Then for every $g\in\mathbb{E}^{\times}$
the map taking $z\in\mathbb{E}$ to $(gh)z^{\rho}(gh)^{-\rho}$ is the reflection
inverting $g$. \label{ref2}
\end{lem}

\begin{proof}
As $h^{\rho}=-h$, this map takes $z$ to $-gz^{\rho}g^{-\rho}$. It is clear that
$g$ is inverted (as $g^{\rho}=z^{\rho}$ cancels with $g^{-\rho}$). Now,
elements which are perpendicular to $g$ are those from $g\mathbb{E}_{0}$ (so
that multiplying by $g^{\rho}$ yields an element of
$N^{\mathbb{E}}_{\mathbb{F}}(g)\mathbb{E}_{0}=\mathbb{E}_{0}$, on which
$Tr^{\mathbb{E}}_{\mathbb{F}}$ vanishes). They are all multiples of $gh$. As
$(gh)^{\rho}=-hg^{\rho}$, a similar calculation shows that $gh$ is invariant
under this operation. This proves the lemma.
\end{proof}

Using all this, we can now establish
\begin{thm}
A special orthogonal group of a 1-dimensional space is a one-element trivial
group. For a 2-dimensional space of discriminant $d$, the Gspin group is
$\mathbb{E}^{\times}$, and the spin and special orthogonal groups are isomorphic
to $\mathbb{E}^{1}$. \label{dim12}
\end{thm}

\begin{proof}
The 1-dimensional part was already proven. Lemma \ref{ref2} and Proposition
\ref{CDT} show that the map the semi-direct product defined in Lemma \ref{ac2}
surjects onto $O(\mathbb{E})$. As $\rho$ represents an element of $O(\mathbb{E})
\setminus SO(\mathbb{E})$ (it inverts $\mathbb{E}_{0}$ and leaves $\mathbb{F}$
invariant), Lemma \ref{ref2} shows that $\mathbb{E}^{\times}$ maps to
$SO(\mathbb{E})$. As $\mathbb{E}^{\times}$ has the same index 2 in the
semi-direct product as $SO(\mathbb{E})$ has in $O(\mathbb{E})$, this map is also
surjective, with kernel $\mathbb{F}^{\times}$. Hence $\mathbb{E}^{\times}$ is
$Gspin(\mathbb{E})$, and $SO(\mathbb{E})$ is isomorphic to $\mathbb{E}^{1}$.
Now, the fact that $\rho$ has spinor norm $-d$ (as this is the vector norm of
non-zero elements of $\mathbb{E}_{0}$, up to $(\mathbb{F}^{\times})^{2}$),
implies that the image of $g\in\mathbb{E}^{\times}$ in $SO(\mathbb{E})$ has
spinor norm $N^{\mathbb{E}}_{\mathbb{F}}(g)$: Indeed, Lemma \ref{ref2} shows
that its composition with $\rho$ inverts an element of norm
$-dN^{\mathbb{E}}_{\mathbb{F}}(g)$, and the spinor norm is a group homomorphism.
Thus $SO^{1}(\mathbb{E})$ is the image of elements
$g\in\mathbb{E}^{(\mathbb{F}^{\times})^{2}}$, and as we may divide by elements
of the kernel $\mathbb{F}^{\times}$ of $\mathbb{E}^{\times} \to SO(\mathbb{E})$,
it suffices to consider $g\in\mathbb{E}^{1}$. Hence the map $\mathbb{E}^{1} \to
SO^{1}(\mathbb{E})$, which is just $g \mapsto g^{2}$ since $g^{-\rho}=g$ for
$g\in\mathbb{E}^{1}$, is surjective, and the kernel is just
$\mathbb{E}^{1}\cap\mathbb{F}^{\times}=\{\pm1\}$. It follows that
$spin(\mathbb{E})=\mathbb{E}^{1}$ as well, and $SO^{1}(\mathbb{E})$ is the group
$(\mathbb{E}^{1})^{2}$ of squares of elements from $\mathbb{E}^{1}$. This proves
the proposition.
\end{proof}

Another way to write the groups $Gspin(\mathbb{E})$ and $spin(\mathbb{E})$ are
as $GU_{\mathbb{E},\rho}(1)$ and $U_{\mathbb{E},\rho}(1)$ respectively, by Lemma
\ref{uniEB}. We remark that when $SO(\mathbb{E})$ is given in terms of
multiplication from $\mathbb{E}^{1}$, the spinor norm of $u\in\mathbb{E}^{1}$,
can be evaluated as $N^{\mathbb{E}}_{\mathbb{F}}(1+u)$ for $u\neq-1$ and $d$ for
$u=-1$ (note that the latter represents $-Id_{\mathbb{E}}$). To see this, write
$u=\frac{g}{g^{\rho}}$ for $g\in\mathbb{E}^{\times}$, so that $1+u$ equals
$\frac{Tr^{\mathbb{E}}_{\mathbb{F}}(g)}{N^{\mathbb{E}}_{\mathbb{F}}(g)}g$, and
$N^{\mathbb{E}}_{\mathbb{F}}(1+u) \in
N^{\mathbb{E}}_{\mathbb{F}}(g)(\mathbb{F}^{\times})^{2}$ since
$g\not\in\mathbb{E}_{0}$ for $u\neq-1$. Comparing these models we find that an
element $u\in\mathbb{E}^{1}$ (other than $-1$) lies in $(\mathbb{E}^{1})^{2}$ if
and only if $N^{\mathbb{E}}_{\mathbb{F}}(1+u)\in\mathbb{F}^{\times}$, a fact
which may also be verified directly. The remaining element $-1$ lies in
$SO^{1}(\mathbb{E})$ if and only if its spinor norm is a square, i.e., if
$\mathbb{E}=\mathbb{F}(\sqrt{-1})$. It is easy to verify that
$-1\in(\mathbb{E}^{1})^{2}$ precisely then this is indeed the case.

As a special case of Theorem \ref{dim12} we obtain
\begin{cor}
A 2-dimensional quadratic space is isotropic if and only if it has a trivial
discriminant. In this case the Gspin group is to
$\mathbb{F}^{\times}\times\mathbb{F}^{\times}$, while the spin and special
orthogonal groups are isomorphic to $\mathbb{F}^{\times}$. \label{iso2}
\end{cor}

\begin{proof}
By Lemma \ref{sp2}, an isotropic 2-dimensional quadratic space comes, up to
rescaling, from a quadratic algebra contains non-zero norm 0 elements. But such
an algebra cannot be a field, which is equivalent to $d$ being a square. In this
case $\mathbb{E}=\mathbb{F}\times\mathbb{F}$, so that
$Gspin(\mathbb{F}\times\mathbb{F})=\mathbb{E}^{\times}=\mathbb{F}^{\times}
\times\mathbb{F}^{\times}$. The group $\mathbb{E}^{1}$ (which is the spin group)
consists of the pairs $\big(r,\frac{1}{r}\big)$ with $r\in\mathbb{F}^{\times}$,
so it is isomorphic to $\mathbb{F}^{\times}$. As $\mathbb{F}^{\times}$ is
embedded in $\mathbb{E}^{\times}$ diagonally, the quotient
$SO(\mathbb{F}\times\mathbb{F})$ is the isomorphic image of the subgroup
$\{(r,1)|r\in\mathbb{F}^{\times}\}$, which is also a copy of
$\mathbb{F}^{\times}$. This proves the corollary.
\end{proof}
We remark that in the case presented in Corollary \ref{iso2}, the element
$\big(r,\frac{1}{r}\big)$ of $\mathbb{E}^{1}$ (considered as
$SO(\mathbb{F}\times\mathbb{F})$ now) is $\frac{g}{g^{\rho}}$ for $g=(r,1)$, so
that its spinor norm is just $r$. This value coincides with the norm of
$\big(1+r,1+\frac{1}{r}\big)$ for $r\neq-1$ and with $d=-1$ for $r=-1$. Hence
$SO(\mathbb{F}\times\mathbb{F})=\mathbb{F}^{\times}$ modulo
$(\mathbb{F}^{\times})^{2}$. The group $SO^{1}(\mathbb{F}\times\mathbb{F})$ is
just $(\mathbb{F}^{\times})^{2}$, given as the quotient of
$Spin(\mathbb{F}\times\mathbb{F})=\mathbb{E}^{1}\cong\mathbb{F}^{\times}$ modulo
$\{\pm1\}$.

The space appearing in Corollary \ref{iso2} is called a \emph{hyperbolic plane}.
It may also be generated by two isotropic vectors with non-zero pairing, so that
it is isometric to all its rescalings. In fact, every isotropic quadratic space
contains a hyperbolic plane, and the complement is uniquely determined up to
isomorphism by the Witt Cancelation Theorem.

\medskip

In dimension 3 we can assume, by rescaling , that our quadratic form has
determinant in $(\mathbb{F}^{\times})^{2}$ (hence discriminant $-1$). Such a
vector space (a 3-dimensional quadratic space with discriminant 1) will be
called a \emph{traceless quaternionic space} over $\mathbb{F}$, for a reason to
be explained by the following
\begin{lem}
If $B$ is a quaternion algebra over $\mathbb{F}$, then the space $B_{0}$ with
the vector norm $|x|^{2}=N^{B}_{\mathbb{F}}(x)$ is a traceless quaternionic
space. Every traceless quaternionic space is isometric to a space which obtained
in this way from some quaternion algebra $B$. The pairing on such a space is
given by $2\langle x,y \rangle=Tr^{\mathbb{E}}_{\mathbb{F}}(x\overline{y})$ for
$x$ and $y$ in $B_{0}$. \label{sp3}
\end{lem}

\begin{proof}
The latter formula for the pairing on $B_{0}$ is a consequence of Lemmas
\ref{BEpol} and \ref{Vpol} (in fact, the same formula holds for $x$ and $y$ in
$B$ with this quadratic form). Now, two elements of $B_{0}$ are orthogonal if
and only if they anti-commute. Writing $B$ is
$\big(\frac{\alpha,\beta}{\mathbb{F}}\big)$, we get two such elements having
norms $-\alpha$ and $-\beta$, and as their product is orthogonal to both of them
of squares to $-\alpha\beta$ (hence with norm $+\alpha\beta$), we get the
required determinant in $(\mathbb{F}^{\times})^{2}$. Conversely, if two
orthogonal elements of a traceless quaternionic space have norms $-\alpha$ and
$-\beta$ respectively, then the determinant condition shows that a generator for
their orthogonal complement can be normalized to have norm $+\alpha\beta$, and
this space is isometric to $B_{0}$ for
$\big(\frac{\alpha,\beta}{\mathbb{F}}\big)$. This proves the lemma.
\end{proof}

Regarding the group acting here we have
\begin{lem}
The group $B^{\times}$ is mapped into $O(B_{0})$ via $g:u \mapsto
gu\overline{g}/N^{B}_{\mathbb{F}}(g)$, with kernel $\mathbb{F}^{\times}$.
Letting $-1$ operate as $-Id_{B_{0}}$ yields a map from the direct product
$B^{\times}\times\{\pm1\}$ into $O(B_{0})$. \label{ac3}
\end{lem}

\begin{proof}
Since $N^{B}_{\mathbb{F}}(g)=\overline{g}g$, $g$ maps $u$ to $gug^{-1}$, and the
multiplicativity of the norm shows that $|gug^{-1}|^{2}=|u|^{2}$. An element
lies in the kernel if and only if it is central (since with the complement
$\mathbb{F}$ of $B_{0}$ in $B$ it does commute), so that the kernel is indeed
$\mathbb{F}^{\times}$. The centrality of $-Id_{B_{0}}$ in $O(B_{0})$ yields the
last assertion. This proves the lemma.
\end{proof}

We remark that the operation of $-1$ coincides with the operation of the main
involution of $B$. the analysis of the orthogonal group begins with the
following
\begin{lem}
If $g \in B_{0}$ has non-zero vector norm, then the orthogonal transformation
taking $u \in B_{0}$ to $-\frac{gu\overline{g}}{N^{B}_{\mathbb{F}}(g)}$ is the
reflection inverting $g$. \label{ref3}
\end{lem}

\begin{proof}
The proof of Lemma \ref{sp3} shows that conjugation by $g$ inverts $g^{\perp}$,
and it clearly leaves $g$ invariant. Composing with the central map
$-Id_{B_{0}}$, we establish the lemma.
\end{proof}

We can now prove
\begin{thm}
The Gspin group of a 3-dimensional $\mathbb{F}$-space is the group $B^{\times}$
for some quaternion algebra $B$ over $\mathbb{F}$, which is generated by
invertible elements of $B_{0}$. The spin is the subgroup $B^{1}$ arising from
this quaternion algebra $B$. \label{dim3}
\end{thm}

\begin{proof}
Proposition \ref{CDT} and Lemma \ref{ref3}, the operation in which is the
action, from Lemma \ref{ac3}, of $g$ on $-u$, imply that the map
$B^{\times}\times\{\pm1\} \to O(B_{0})$ is surjective. As reflections and
$-Id_{B_{0}}$ has determinant $-1$, the image of $B^{\times}$ lies in
$SO(B_{0})$, and index considerations show that the map $B^{\times} \to
SO(B_{0})$ is surjective. As the kernel is $\mathbb{F}^{\times}$, we find that
$Gspin(B_{0})=B^{\times}$. Taking out the action of the central element
$-Id_{B_{0}}$, Lemma \ref{ref3} shows that $B_{0} \cap B^{\times}$ indeed
generates $B^{\times}$ (a fact which in this case is easily verified directly),
since the full kernel $\mathbb{F}^{\times}$ is clearly generated by this set:
$t=(tg)g^{-1}$ for $t\in\mathbb{F}^{\times}$. As for spinor norms, we first
observe that under our normalization $-Id$ has spinor norm 1. Hence Lemma
\ref{ref3} and the fact that $|g|^{2}=N^{B}_{\mathbb{F}}(g)$ imply that the
spinor norm of any $g \in B_{0} \cap B^{\times}$ is $N^{B}_{\mathbb{F}}(g)$.
Since such elements were seen to generate $B^{\times}$, the spinor norm of any
$g \in B^{\times}$ is $N^{B}_{\mathbb{F}}(g)$. The group $SO^{1}(B_{0})$ is thus
the image of elements having reduced norms in $(\mathbb{F}^{\times})^{2}$, and
by appropriate scalar multiplication we may restrict to elements from $B^{1}$.
As the only scalars in $B^{1}$ are $\pm1$, we find that $B^{1}$ is indeed
$spin(B_{0})$. This proves the proposition.
\end{proof}
Lemmas \ref{uniEB} and \ref{Sp1B} shows that the Gspin group from Theorem
\ref{dim3} can also be described as $GSp_{B}(1)$ and as
$GSU_{\mathbb{K},\sigma}\binom{-\varepsilon\ \ 0}{\ \ 0\ \ 1}$, in case
$\mathbb{K}=\mathbb{F}(\eta)$ is a quadratic extension of $\mathbb{F}$ (with
Galois automorphism $\sigma$) which splits $B$, and $\varepsilon\in\mathbb{F}$
is such that $B\cong\big(\frac{\eta,\varepsilon}{\mathbb{F}}\big)$. They also
imply that the spin group in question is isomorphic to $Sp_{B}(1)$, as well as
to $SU_{\mathbb{K},\sigma}\binom{-\varepsilon\ \ 0}{\ \ 0\ \ 1}$ for such
$\mathbb{K}$, $\sigma$, and $\varepsilon$.

The isotropic case in dimension 3 is given in
\begin{cor}
A quadratic space of dimension 3 is isotropic if and only it is related to the
split quaternion algebra $B=M_{2}(\mathbb{F})$. The Gspin group
$Gspin\big(M_{2}(\mathbb{F})_{0}\big)$ is then $GL_{2}(\mathbb{F})$, and
$spin\big(M_{2}(\mathbb{F})_{0}\big)=SL_{2}(\mathbb{F})$. We also have
$SO\big(M_{2}(\mathbb{F})_{0}\big)=PGL_{2}(\mathbb{F})$ and
$SO^{1}\big(M_{2}(\mathbb{F})_{0}\big)=PSL_{2}(\mathbb{F})$. \label{iso3}
\end{cor}

\begin{proof}
If $B_{0}$ is isotropic then $B$ cannot be a division algebra, and the space
$M_{2}(\mathbb{F})_{0}$ does split. The Gspin and spin groups are determined by
Theorem \ref{dim3} with $B=M_{2}(\mathbb{F})$, and dividing the former by
$\mathbb{F}^{\times}$ and the latter by $\{\pm1\}$ yields the asserted
projective groups. This proves the corollary.
\end{proof}
Lemma \ref{uniEB} allows us to write the Gspin and spin groups from Corollary
\ref{iso3} also as $GSU_{\mathbb{K},\sigma}\binom{-1\ \ 0}{\ \ 0\ \ 1}$ and
$SU_{\mathbb{K},\sigma}\binom{-1\ \ 0}{\ \ 0\ \ 1}$ respectively, for any
quadratic extension $\mathbb{K}$ of $\mathbb{F}$ (with Galois automorphism
$\sigma$), since every such $\mathbb{K}$ splits $M_{2}(\mathbb{F})$. In this
case matrix transposition is also an element of
$O\big(M_{2}(\mathbb{F})_{0}\big) \setminus SO\big(M_{2}(\mathbb{F})_{0}\big)$,
and this element arises as the composition of $-Id$ and conjugation by
$\binom{0\ \ -1}{1\ \ \ \ 0}$ (see Lemma \ref{Sadjt}).

In this split case there is an additional assertion, which is given by
\begin{cor}
The groups $GL_{2}(\mathbb{F})$, $SL_{2}(\mathbb{F})$, $PGL_{2}(\mathbb{F})$,
and $PSL_{2}(\mathbb{F})$ are the Gspin, spin, special orthogonal, and spinor
norm kernel groups of the space $M_{2}^{sym}(\mathbb{F})$ of symmetric
$2\times2$ matrices over $\mathbb{F}$, on which they all operate via
$g:X\mapsto\frac{gXg^{t}}{\det g}$. \label{alt3}
\end{cor}

\begin{proof}
As right multiplication by $\binom{\ \ 0\ \ 1}{-1\ \ 0}$ takes
$M_{2}(\mathbb{F})_{0}$ to $M_{2}^{sym}(\mathbb{F})$ and preserves determinants,
the corollary follows from Corollary \ref{iso3} and Lemma \ref{Sadjt}.
\end{proof}
Note that the operation from Corollary \ref{alt3} replaces the symplectic main
involution on $B=M_{2}(\mathbb{F})$ by an orthogonal one, and the usual space
$B^{-}$ by an $M_{2}(\mathbb{F})^{+}$ space. The possible generators of
$O\big(M_{2}^{sym}(\mathbb{F})\big)/SO\big(M_{2}^{sym}(\mathbb{F})\big)$ are
again $-Id$ (which is now $X\mapsto-X^{t}$) and the adjoint involution.

\section{Dimension 4 \label{Dim4}}

Given a quadratic space of dimension 4 and discriminant $d$ over $\mathbb{F}$,
we define $\mathbb{E}$ to be $\mathbb{F}(\sqrt{d})$, with automorphism $\rho$.
Given a quaternion algebra $B$ over $\mathbb{F}$, the tensor product
$B_{\mathbb{E}}$ comes endowed with the (unitary) involution
$\iota\otimes\rho:x\mapsto\overline{x}^{\rho}$. Our space is given by
\begin{lem}
The space $B_{\mathbb{E}}^{-}$ of this involution becomes, when endowed with
the quadratic form $|x|^{2}=N^{B_{\mathbb{E}}}_{\mathbb{E}}(x)$, a quadratic
space over $\mathbb{F}$ with discriminant $d$, in which $2\langle x,y
\rangle=Tr^{B_{\mathbb{E}}}_{\mathbb{E}}(x\overline{y})$ holds for every $x$
and $y$. Every such space is obtained, up to rescaling in this way. \label{sp4}
\end{lem}

\begin{proof}
$B_{\mathbb{E}}^{-}$ is contained in the quadratic space $B_{\mathbb{E}}$ over
$\mathbb{E}$ with the same vector norm, and the combination of Lemmas
\ref{BEpol} and \ref{Vpol} shows that the formula for the pairing holds (in
$\mathbb{E}$) for any two elements of the larger space. Now,
$B_{\mathbb{E}}^{-}$ is $\mathbb{E}_{0} \oplus B_{0}$ inside $B_{\mathbb{E}}$,
and the direct sum is thus orthogonal inside there. Since
$N^{B_{\mathbb{E}}}_{\mathbb{E}}$ coincides with $N^{B}_{\mathbb{F}}$ on $B_{0}$
and is the square map on $\mathbb{E}_{0}$ (both $\mathbb{F}$-valued), we find
that $|x|^{2}\in\mathbb{F}$ for every $x \in B_{\mathbb{E}}^{-}$. Moreover,
$B_{0}$ has determinant 1 and $\mathbb{E}_{0}$ is spanned by an element $h$ with
$|h|^{2}=h^{2}=d$, so that the determinant and discriminant of such a space are
$d$. Conversely, given a quadratic space of discriminant $d$, we may rescale our
it such that a anisotropic element $v$ of our choice has vector norm $d$. The
subspace $v^{\perp}$ is a traceless quaternionic space, so that by Lemma
\ref{sp3} it can be presented as $B_{0}$ for a quaternion algebra $B$ over
$\mathbb{F}$. Hence we found a presentation of our space as $\mathbb{E}_{0}
\oplus B_{0}$, namely $B_{\mathbb{E}}^{-}$, with
$|x|^{2}=N^{B_{\mathbb{E}}}_{\mathbb{E}}(x)$. This proves the lemma.
\end{proof}

Note that the proof of Lemma \ref{sp4} involved a choice of a vector, and
choosing another vector (with the appropriate rescaling) may lead to other
quaternion algebras which are not isomorphic to $B$. However, $B_{\mathbb{E}}$
remains the same algebra, but with a different unitary involution. The
correspondence between quaternion algebras over $\mathbb{F}$ which are
contained in $B_{\mathbb{E}}$ and generate it over $\mathbb{E}$ and unitary
involutions on $B_{\mathbb{E}}$ is given in Proposition 2.22 of \cite{[KMRT]}.
However, we shall consider, $B$, $B_{\mathbb{E}}$, and the involution as fixed.
As $\iota$ is canonical, this gives an interpretation of $\rho$ on
$B_{\mathbb{E}}$ as well. We also remark that the complementary space
$B_{\mathbb{E}}^{+}$ is obtained from $B_{\mathbb{E}}^{-}$ via multiplication by
an element of $\mathbb{E}_{0}$, so that it is isometric to a rescaling of
$B_{\mathbb{E}}^{-}$ as a quadratic space over $\mathbb{F}$.

$B_{\mathbb{E}}$ operates on itself by $g:x \mapsto gx\overline{g}^{\rho}$, and
this action preserves the eigenspaces $B_{\mathbb{E}}^{\pm}$. Moreover, the two
spaces are invariant under $\iota$. Using this, we now prove
\begin{lem}
The group $B_{\mathbb{E}}^{\mathbb{F}^{\times}}$ maps into
$O(B_{\mathbb{E}}^{-})$ if an element $g$ operates as
$x\mapsto\frac{gx\overline{g}^{\rho}}{N^{B_{\mathbb{E}}}_{\mathbb{E}}(g)}$. The
kernel of this map is $\mathbb{F}^{\times}$. Let $\tilde{\iota}$ be the
non-trivial element of a cyclic group of order 2, which operates on
$B_{\mathbb{E}}^{\mathbb{F}^{\times}}$ as $\rho$. Sending $\tilde{\iota}$ to
operate as $\iota$ defines a group homomorphism from the semi-direct product
of $\{1,\tilde{\iota}\}$ and $B_{\mathbb{E}}^{\mathbb{F}^{\times}}$ to
$O(B_{\mathbb{E}}^{-})$. \label{ac4}
\end{lem}

\begin{proof}
For $g \in B_{\mathbb{E}}^{\mathbb{F}^{\times}}$ we may replace
$\frac{\overline{g}^{\rho}}{N^{B_{\mathbb{E}}}_{\mathbb{E}}(g)}$ by just
$g^{-\rho}$, and the fact that $N^{B_{\mathbb{E}}}_{\mathbb{E}}$ is
multiplicative implies that the equality
$\Big|\frac{gx\overline{g}^{\rho}}{N^{B_{\mathbb{E}}}_{\mathbb{E}}(g)}\Big|^{2}
=|x|^{2}$ holds for every $x \in B_{\mathbb{E}}^{-}$. An element $g$ is in the
kernel of this action if and only if it is central and satisfies
$gg^{\rho}=N^{B_{\mathbb{E}}}_{\mathbb{E}}(g)=g^{2}$, which is equivalent to $g$
being in $\mathbb{F}^{\times}$. We also know that $\iota$ preserves
$B_{\mathbb{E}}^{-}$ and $|\overline{x}|^{2}=|x|^{2}$ for elements of that
space. The equality
$\overline{gx\overline{g}^{\rho}}=g^{\rho}\overline{x}\overline{g}$ shows that
the map to $O(B_{\mathbb{E}}^{-})$ respects the product rule of the semi-direct
product, which completes the proof of the lemma.
\end{proof}

We remark that other elements of $B_{\mathbb{E}}^{\times}$ do not increase the
image of the map $B_{\mathbb{E}}^{\mathbb{F}^{\times}} \to
O(B_{\mathbb{E}}^{-})$ from Lemma \ref{ac4}. Indeed, if $g \in
B_{\mathbb{E}}^{\times}$ and $t\in\mathbb{E}^{\times}$ are such that
$x\mapsto\frac{gx\overline{g}^{\rho}}{t}$ preserves $B_{\mathbb{E}}^{-}$ and is
orthogonal, then $t\in\mathbb{F}^{\times}$, the number
$\frac{N^{B_{\mathbb{E}}}_{\mathbb{E}}(g)}{t}$ lies in $\mathbb{E}^{1}$ hence
equals $\frac{s^{\rho}}{s}$ for some $s\in\mathbb{E}^{\times}$ by Hilbert's
Theorem 90. But then $sg \in B_{\mathbb{E}}^{\mathbb{F}^{\times}}$ and operates
like $x\mapsto\frac{gx\overline{g}^{\rho}}{t}$. This assertion will also follow
from Theorem \ref{dim4} below. As for the kernel, note that a non-zero element
$r\in\mathbb{E}_{0}$ is also central and lies
$B_{\mathbb{E}}^{\mathbb{F}^{\times}}$, but as $rr^{\rho}=-r^{2}$ and
$N^{B_{\mathbb{E}}}_{\mathbb{E}}(r)=r^{2}$, such elements operate as $-Id$
rather than trivially.

The properties of the larger homomorphism from Lemma \ref{ac4} will follow from
\begin{lem}
For an element $g \in B_{\mathbb{E}}^{-} \cap B_{\mathbb{E}}^{\times}$, the
reflection inverting $g$ lies in the image of the map from Lemma \ref{ac4},
being
$x\mapsto\frac{g\overline{x}\cdot\overline{g}^{\rho}}{N^{B_{\mathbb{E}}}_{
\mathbb{E}}(g)}$.
\label{ref4}
\end{lem}

\begin{proof}
The proof of Lemma \ref{sp4} shows that every such $g$ is in
$B_{\mathbb{E}}^{\mathbb{F}^{\times}}$, so that the latter transformation comes
from the semi-direct product appearin in Lemma \ref{ac4}. Now,
$\overline{g}^{\rho}=-g$ for $g \in B_{\mathbb{E}}^{-}$, and
$N^{B_{\mathbb{E}}}_{\mathbb{E}}(g)=\overline{g}g$. Thus, for $x=g$ the result
of the action is $-g$, while if $x \in g^{\perp}$ the equality from Lemma
\ref{ac4} allows us to replace $-g\overline{x}$ by $+x\overline{g}$, so that the
total expression is just $x$. This proves the lemma.
\end{proof}

The groups obtained in dimension 4 are now given in the following
\begin{thm}
The quadratic space $B_{\mathbb{E}}^{-}$ has Gspin group
$B_{\mathbb{E}}^{\mathbb{F}^{\times}}$, and it is generated by
$B_{\mathbb{E}}^{-} \cap B_{\mathbb{E}}^{\times}$. The spin group is
$B_{\mathbb{E}}^{1}$. \label{dim4}
\end{thm}

\begin{proof}
The surjectivity of the map from the semi-direct product from Lemma \ref{ac4}
onto $O(B_{\mathbb{E}}^{-})$ follows from Lemma \ref{ref4}, and Proposition
\ref{CDT}. The fact that $\iota$ has determinant $-1$ in
$O(B_{\mathbb{E}}^{-})$ and index considerations show that
$B_{\mathbb{E}}^{\mathbb{F}^{\times}}$ maps onto $SO(B_{\mathbb{E}}^{-})$. The
kernel of latter map being $\mathbb{F}^{\times}$ by Lemma \ref{ac4}, we find
that $Gspin(B_{\mathbb{E}}^{-})=B_{\mathbb{E}}^{\mathbb{F}^{\times}}$. The
semi-direct product structure shows, with Lemma \ref{ref4}, that
$B_{\mathbb{E}}^{-} \cap B_{\mathbb{E}}^{\times}$ generates
$B_{\mathbb{E}}^{\mathbb{F}^{\times}}$ (since again the kernel
$\mathbb{F}^{\times}$ is not a problem), a fact which also in this case may
still be verified directly. Now, $\iota$ reflects the traceless quaternionic
space $B_{0}$, of determinant 1, so that its spinor norm is 1. Lemma \ref{ref4}
thus implies that the spinor norm of $g \in B_{\mathbb{E}}^{-} \cap
B_{\mathbb{E}}^{\times}$ is $|g|^{2}=N^{B_{\mathbb{E}}}_{\mathbb{E}}(g)$. As
such elements were seen to generate $B_{\mathbb{E}}^{\mathbb{F}^{\times}}$, the
multiplicativity of the norm implies that the spinor norm of any $g \in
B_{\mathbb{E}}^{\mathbb{F}^{\times}}$ is $N^{B_{\mathbb{E}}}_{\mathbb{E}}(g)$.
Elements whose (spinor, hence algebra) norms lie in $(\mathbb{F}^{\times})^{2}$
can be dividing by scalars from the kernel $\mathbb{F}^{\times}$, and land in
$B_{\mathbb{E}}^{1}$. Thus $B_{\mathbb{E}}^{1}$ maps surjectively onto
$SO^{1}(B_{\mathbb{E}}^{-})$, and the kernel consists of those scalars whose
norm (hence square) is 1. As these are just $\pm1$, $B_{\mathbb{E}}^{1}$ is the
spin group, which completes the proof of the proposition.
\end{proof}

In addition to $\iota$, $\rho$ also represents an element of
$O(B_{\mathbb{E}}^{-})$. Its spinor norm is $d$, being the composition of
$\iota$ and $-Id$ as well as being the reflection in a generator of
$\mathbb{E}_{0}$. However, $\iota$ is a more canonical representative of
$O(B_{\mathbb{E}}^{-})/SO(B_{\mathbb{E}}^{-})$, being independent of the
$\mathbb{F}$-structure on $B_{\mathbb{E}}$. Note that Theorems \ref{dim3} and
\ref{dim4} imply that if $d$ is not a square then any spin group arising from a
4-dimensional quadratic space over $\mathbb{F}$ with discriminant $d$ is
isomorphic to the spin group of a suitable 3-dimensional quadratic space over
$\mathbb{E}=\mathbb{F}(\sqrt{d})$. The converse does not hold though, as we need
the quaternion algebra over $\mathbb{E}$ to come from one over $\mathbb{F}$.

Recall now that every quaternion algebra $B$ over $\mathbb{F}$ becomes, with
$|x|^{2}=N^{B}_{\mathbb{F}}$, a quadratic space of discriminant 1. Lemma
\ref{sp4} and Theorem \ref{dim4} has the following
\begin{cor}
Any 4-dimensional space of discriminant 1 is isometric to some quaternion
algebra $B$ over $\mathbb{F}$ with its reduced norm (perhaps rescaled). The
Gspin group $Gspin(B)$ consists of those pairs in $B^{\times} \times B^{\times}$
having the same reduced norm, operating via left multiplication and inverted
right multiplication, and $spin(B)=B^{1} \times B^{1}$ with the same operation.
\label{dim4d1}
\end{cor}

\begin{proof}
The condition $d\in(\mathbb{F}^{\times})^{2}$ implies that
$\mathbb{E}=\mathbb{F}\times\mathbb{F}$. Hence $B_{\mathbb{E}}=B \times B$, and
$B_{\mathbb{E}}^{-}$ is the subspace $\{(x,-\overline{x})|x \in B\}$ of $B
\times B$. As
$N^{B_{\mathbb{E}}}_{\mathbb{E}}(x,-\overline{x})=N^{B}_{\mathbb{F}}(x)$, the
first assertion follows from Lemma \ref{sp4}. Now,
$B_{\mathbb{E}}^{\mathbb{F}^{\times}}$ consists, by definition, of those pairs
$(g,h) \in B^{\times} \times B^{\times}$ such that
$N^{B}_{\mathbb{F}}(g)=N^{B}_{\mathbb{F}}(h)$, sending $(x,-\overline{x})$ to
$(gx\overline{h},-h\overline{x}\cdot\overline{g})$ (the second entry being also
$-\overline{gx\overline{h}}$), divided by the common norm of $g$ and $h$. In
terms of the operation with $(g,h)^{-\rho}$ from the proof of Lemma \ref{sp4},
the action on $B$ is given by $(g,h):x \mapsto gxh^{-1}$. Restricting to
elements of norm 1, the spin group is seen to be $B^{1} \times B^{1}$. This
proves the corollary.
\end{proof}

The spin group from Corollary \ref{dim4d1} is the product of two copies of a
spin group of a 3-dimensional space over $\mathbb{F}$, which complements the
relation to spaces over the quadratic extension $\mathbb{E}$ in the other
discriminants. Lemmas \ref{uniEB} and \ref{Sp1B} allow us to write such spin
groups as $Sp_{B}(1) \times Sp_{B}(1)$, as well as
$SU_{\mathbb{K},\sigma}\binom{-\varepsilon\ \ 0}{\ \ 0\ \ 1} \times
SU_{\mathbb{K},\sigma}\binom{-\varepsilon\ \ 0}{\ \ 0\ \ 1}$ in case $B$ is
isomorphic to $\big(\frac{\eta,\varepsilon}{\mathbb{F}}\big)$ and
$\mathbb{K}=\mathbb{F}(\eta)$ with Galois automorphism $\sigma$.

We note that the $\mathbb{F}$-structure on $B_{\mathbb{E}}$, i.e., the
quaternion algebra over $\mathbb{F}$ which yields $B_{\mathbb{E}}$ after
tensoring with $\mathbb{E}$, is equivalent to the choice of the automorphism on
$B_{\mathbb{E}}$ which we denoted here also by $\rho$, or equivalently, by the
involution $x\mapsto\overline{x}^{\rho}$ on $B_{\mathbb{E}}$ (see Proposition
2.22 of \cite{[KMRT]}). Proposition 2.18 shows that all these involutions are
related, and are in one-to-one correspondence with the space
$(B_{\mathbb{E}}^{+} \cap B_{\mathbb{E}}^{\times})/\mathbb{F}^{\times}$, since
every such involution is $x \mapsto b\overline{x}^{\rho}b^{-1}$ for invertible
$b \in B_{\mathbb{E}}^{+}$ which is uniquely determined up to multiplication
from $\mathbb{F}^{\times}$. However, we shall not need these results in what
follows.

As for isotropy, here we have
\begin{cor}
The space $B_{\mathbb{E}}^{-}$ is isotropic if and only if $\mathbb{E}$ splits
$B$. We may then take $B=M_{2}(\mathbb{F})$.
$Gspin\big(M_{2}(\mathbb{E})^{-}\big)$ is
$GL_{2}^{\mathbb{F}^{\times}}(\mathbb{E})$, the spin group
$spin\big(M_{2}(\mathbb{E})^{-}\big)$ as $SL_{2}(\mathbb{E})$, and
$SO^{1}\big(M_{2}(\mathbb{E})^{-}\big)$ as $PSL_{2}(\mathbb{E})$.
\label{iso4}
\end{cor}

\begin{proof}
If $B_{\mathbb{E}}^{-}$ is isotropic then $B_{\mathbb{E}}$ cannot be a division
algebra. Conversely, if $\mathbb{E}$ splits $B$ then there is an embedding
$i:\mathbb{E} \to B$, and if $r\in\mathbb{E}_{0}$ then $r+i(r)$ belongs to
$B_{\mathbb{E}}^{-}$ and is a zero-divisor (hence isotropic). By splitting a
hyperbolic plane and rescaling a vector $v$ which is perpendicular to this
hyperbolic plane to have vector norm $d$. Then $v^{\perp}$ is isotropic, and
the corresponding quaternion algebra $B$ splits by Corollary \ref{iso3}. The
Gspin and spin groups are given in Theorem \ref{dim4} (written in terms of a
split algebra), and the assertion about $SO^{1}\big(M_{2}(\mathbb{E})^{-}\big)$,
which is $spin\big(M_{2}(\mathbb{E})^{-}\big)/\{\pm1\}$, is immediate. This
proves the corollary.
\end{proof}

Also here we can consider matrix conjugation an element of
$O\big(M_{2}(\mathbb{E})^{-}\big)$ of determinant 1, and it comes as the
composition of the main involution (adjoint) and conjugation by $\binom{\ \ 0\ \
1}{-1\ \ 0}$ by Lemma \ref{Sadjt}. In this case we can get an equivalent
quadratic space, as is given in the following
\begin{cor}
We can consider the groups $GL_{2}^{\mathbb{F}^{\times}}(\mathbb{E})$ and
$SL_{2}(\mathbb{E})$ as the Gspin and spin groups of the space
$M_{2}^{Her}(\mathbb{E},\rho)$ of $2\times2$ matrices of $\mathbb{E}$ which are
Hermitian with respect to $\rho$, with the vector norm being the determinant.
The operation is via $g:X\mapsto\frac{gXg^{t\rho}}{\det g}$. In the case of
trivial discriminant, we may consider the subgroup of
$GL_{2}(\mathbb{F}) \times GL_{2}(\mathbb{F})$ consisting of pairs of matrices
of the same determinant and $SL_{2}(\mathbb{F}) \times SL_{2}(\mathbb{F})$ as
the Gspin and the spin groups of $M_{2}(\mathbb{F})$ with the determinant also
via the action $(g,h):M \mapsto gMh^{t}$ divided by the common norm of $g$ an
$h$ (which is trivial on the latter group). \label{alt4}
\end{cor}

\begin{proof}
The first assertion follows directly from Corollary \ref{iso3} and Lemma
\ref{Sadjt}, since right multiplication by $\binom{\ \ 0\ \ 1}{-1\ \ 0}$ takes
$M_{2}(\mathbb{E})^{-}$ to $M_{2}^{Her}(\mathbb{E})$. The second assertion is
obtained from the same considerations together with Corollary \ref{dim4d1}.
\end{proof}

The quotient group
$O\big(M_{2}^{Her}(\mathbb{E},\rho)\big)/SO\big(M_{2}^{Her}(\mathbb{E},
\rho)\big)$, as well as the group
$O\big(M_{2}(\mathbb{F})\big)/SO\big(M_{2}(\mathbb{F})\big)$, is again
generated by adjoint or transposition, but here $\rho$ coincides with the latter
transformation.

\section{Dimension 6, Discriminant 1 \label{Dim6d1}}

Our presentation of 6-dimensional spaces of discriminant 1 is based on
presentations of bi-quaternion algebras over $\mathbb{F}$ as tensor products of
two quaternion algebras, as in
\begin{lem}
For two quaternion algebras over $\mathbb{F}$, $B$ and $C$ say, the subspace
$(B_{0}\otimes1)\oplus(1 \otimes C_{0})$ of the bi-quaternion algebra $A=B
\otimes C$ is a quadratic space of dimension 6 and discriminant 1 if we define
$|x\otimes1+1 \otimes y|^{2}=N^{B}_{\mathbb{F}}(x)-N^{C}_{\mathbb{F}}(y)$.
Every 6-dimensional quadratic space of discriminant 1 may be obtained, up to
isometries and rescalings, in this way. \label{sp6d1}
\end{lem}

\begin{proof}
Lemma \ref{sp3} shows that $(B_{0}\otimes1)\oplus(1 \otimes C_{0})$ is the
direct sum of two 3-dimensional spaces of determinants 1 and $-1$, so that the
total determinant is $-1$ and the discriminant is 1. Conversely, given a
quadratic space $V$ of dimension 6 and discriminant 1, choose a 3-dimensional
non-degenerate subspace of $V$, and rescale $V$ such that the chosen space has
determinant 1. The discriminant 1 condition implies that the orhogonal
complement is also a traceless quaternionic space but rescaled by $-1$, so that
the assertion follows from Lemma \ref{sp3}. This proves the lemma.
\end{proof}

The space from Lemma \ref{sp6d1} is called the \emph{Albert form} of the two
quaternion algebras $B$ and $C$. The bi-quaternion algebra $A=B \otimes C$
comes with the involution $\iota_{B}\otimes\iota_{C}$, which depends on the
presentation of $A$ as $B \otimes C$. When we decompose $A$ as $A^{+} \oplus
A^{-}$ according to this involution, then $A^{+}$ is the 10-dimensional space
$B_{0} \otimes C_{0}\oplus\mathbb{F}(1\otimes1)$, and $A^{-}$ is the space
$(B_{0}\otimes1)\oplus(1 \otimes C_{0})$ from Lemma \ref{sp6d1}. In particular,
$\iota_{B}\otimes\iota_{C}$ is orthogonal---compare Proposition 2.23(1) of
\cite{[KMRT]}. Note that the product structure on our original space depends on
the choice of the 3-dimensional space which we normalize to be
$B_{0}\otimes1$: Observe that $B_{0}\otimes1$ and $1 \otimes C_{0}$ are
precisely those elements of $A^{-}$ whose algebra square lies in $\mathbb{F}$.

Several remarks are in order here. First, there are many involutions, orthogonal
and symplectic, on $A$, and to each of them one may associate a quadratic
6-dimensional space of discriminant 1 (an Albert form). Two Albert forms are
isometric up to rescaling if and only if they come from isomorphic
bi-quaternion algebras, by a result of \cite{[J]} (Lemma \ref{sp6d1} already
shows that this construction gives all the Albert form, up to rescaling).
However, not all the involutions on $A$, and not even all the orthogonal
involutions on $A$, come from a presentation of $A$ as $B \otimes C$, though
the mere existence of an involution of the first kind on a degree 4 central
simple algebra $A$ implies that $A$ is the tensor product of two quaternion
algebras over $\mathbb{F}$ (a theorem of Albert---see Theorem 16.1 of
\cite{[KMRT]}). However, we stick to one fixed involution, which does arise
in this way, and our results are independent of all the results mentioned in
this paragraph. We also remark that Section 16 of \cite{[KMRT]} presents
results which are very similar to ours, but some of the calculations do not
appear there (in particular, the main calculation required for Proposition
\ref{NAFg} below), and we concentrate on a case where many technical aspects
become simpler.

Our analysis is based on the following
\begin{lem}
If $\theta:A^{-} \to A^{-}$ takes $u=x\otimes1+1 \otimes y$ to
$\tilde{u}=-x\otimes1+1 \otimes y$ then $u\tilde{u}=\tilde{u}u=|u|^{2}$ in $A$,
and $2\langle v,w \rangle=v\tilde{w}+w\tilde{v}=\tilde{v}w+\tilde{w}v$.
\label{vnorm6}
\end{lem}

\begin{proof}
The first assertion follows from a simple and direct calculation. The second
equality then follows from Lemma \ref{Vpol}. This proves the lemma.
\end{proof}

Note that $\theta$ might be considered as the restriction of the map $\iota_{B}
\otimes Id_{C}$ on $A$ to $A^{-}$, but as the latter map behaves badly with
respect to products (it neither preserves them nor inverts them), we consider
only the restriction $\theta$ to $A^{-}$. We also remark that interchanging the
roles of $B$ and $C$ just means replacing $\theta$ by $-\theta$ and inverting
the sign of the Albert form.

The reduced norm $N^{A}_{\mathbb{F}}$ is a degree 4 form on $A$ which takes any
tensor product $b \otimes c$ to
$N^{B}_{\mathbb{F}}(b)^{2}N^{C}_{\mathbb{F}}(c)^{2}$. In certain calculations
we shall need to evaluate it, which is done, under some assumptions, in the
following
\begin{prop}
If $A=M_{2}(B)$ for some quaternion algebra $B$ then the equality
\[N^{A}_{\mathbb{F}}\binom{a\ \ b}{c\ \
d}=N^{B}_{\mathbb{F}}(a)N^{B}_{\mathbb{F}}(d)+N^{B}_{\mathbb{F}}(b)N^{B}_{
\mathbb{F}}(c)-Tr^{B}_{\mathbb{F}}(\overline{a}b\overline{d}c)\] holds for
every $\binom{a\ \ b}{c\ \ d} \in A$, i.e., every $a$, $b$, $c$, and $d$ from
$B$. \label{NAexp}
\end{prop}

\begin{proof}
First, the assertion holds in case one of $a$, $b$, $c$, or $d$ is 0. This
follows from evaluation of $4\times4$ determinants in case $B=M_{2}(\mathbb{F})$
and $A=M_{4}(\mathbb{F})$, hence holds in general since $B$ may be considered as
a subalgebra of matrices over a splitting field and the assertions are invariant
under scalar extensions. Now, if one entry is invertible then we can determine
the reduced norm by right multiplication with a matrix of the sort $\binom{1\ \
x}{0\ \ 1}$ (of reduced norm 1): For example, if $a$ is invertible then we take
$x=-a^{-1}b$ and find that that $N^{A}_{\mathbb{F}}\binom{a\ \ b}{c\ \ d}$
equals $N^{B}_{\mathbb{F}}(a)N^{B}_{\mathbb{F}}(d-ca^{-1}b)$ by
multiplicativity, and then using Lemma  \ref{BEpol} we get the asserted value.
Similar considerations cover the cases where $b$, $c$, or $d$ are invertible,
which completes the proof in the case where $B$ is a division algebra. In case
$B=M_{2}(\mathbb{F})$ and all of $a$, $b$, $c$, and $d$ are non-zero and not
invertible, we may conjugate everything in $B$ (an operation leaving our
expression invariant) such that $a$ becomes the matrix $\binom{t\ \ 0}{0\ \ 0}$
for some $t\in\mathbb{F}^{\times}$. Observe that left multiplication by
$\binom{1\ \ x}{0\ \ 1}$ takes the upper left entry of our matrix to $a+xc$.
Recall that $\det c=0$ but $c\neq0$. Hence if $c$ has right column 0 then we may
choose $x$ such that $a+xc=0$, while if the right column of $c$ is not 0 we may
choose $x$ such that $x$ such that $xc$ has upper row 0 and lower row non-zero,
so that $a+xc$ becomes invertible. As our expression is invariant under left
multiplication by $\binom{1\ \ x}{0\ \ 1}$, this completes the proof of the
Proposition.
\end{proof}

We can use Proposition \ref{NAexp} to get an explicit expression for the
$N^{A}_{\mathbb{F}}$ in general. Writing an element of $A$ as $\alpha+\delta
i+\beta j+\gamma ij$ where $i$ and $j$ generate $C$, anti-commute, and square to
$\eta$ and $\varepsilon$ respectively, we may embed $C$ into
$M_{2}(\mathbb{K})$, where $\mathbb{K}=\mathbb{F}(\sqrt{\eta})$ with Galois
automorphism $\sigma$, as the algebra $(\mathbb{K},\sigma,\delta)$, and then our
element becomes an element of $A_{\mathbb{K}}=M_{2}(B_{\mathbb{K}})$ for which
we may apply Proposition \ref{NAexp}. The result is
\[N^{A}_{\mathbb{F}}(\alpha+\delta i+\beta j+\gamma
ij)=N^{B}_{\mathbb{F}}(\alpha)^{2}+\eta^{2}N^{B}_{\mathbb{F}}(\delta)^{2}
+\varepsilon^{2}N^{B}_{\mathbb{F}}(\beta)^{2}+\eta^{2}\varepsilon^{2}N^{B}_{
\mathbb{F}}(\gamma)^{2}+\] \[-\eta
Tr^{B}_{\mathbb{F}}\big((\overline{\alpha}\delta)^{2}\big)-\varepsilon
Tr^{B}_{\mathbb{F}}\big((\overline{\alpha}\beta)^{2}\big)+\eta\varepsilon
Tr^{B}_{\mathbb{F}}\big((\overline{\alpha}\gamma)^{2}\big)+\eta\varepsilon
Tr^{B}_{\mathbb{F}}\big((\overline{\delta}\beta)^{2}\big)+\]
\[-\eta^{2}\varepsilon
Tr^{B}_{\mathbb{F}}\big((\overline{\delta}\gamma)^{2}\big)-\eta\varepsilon^{2}
Tr^{B}_{\mathbb{F}}\big((\overline{\beta}\gamma)^{2}\big)-2\eta\varepsilon
Tr^{B}_{\mathbb{F}}\big(\overline{\alpha}\beta\overline{\delta}
\gamma)+2\eta\varepsilon
Tr^{B}_{\mathbb{F}}\big(\overline{\alpha}\gamma\overline{\delta}\beta),\] and it
is seen to reduce to $N^{B}_{\mathbb{F}}(b)^{2}N^{C}_{\mathbb{F}}(c)^{2}$ when
our element is a single tensor $b \otimes c$ (by choosing $i$ to be the
traceless part of $c$ if it is non-zero). More importantly, we have
\begin{cor}
The equality $N^{A}_{\mathbb{F}}(u)=|u|^{4}$ holds for every $u \in A^{-}$.
\label{NAvn2}
\end{cor}

\begin{proof}
One way to see it is by writing $u=x\otimes1+1 \otimes y$ and choosing $i=y$ in
the basis for $C$ for the latter formula (if it does not vanish).
Alternatively, we consider $A=M_{2}(B)$ (with $C=M_{2}(\mathbb{F})$) first,
where elements of $A^{-}=M_{2}(B)^{-}$ take the form $u=\binom{\lambda\ \ -r}{s\
\ -\overline{\lambda}}$, in which $\lambda \in B$ and $r$ and $s$ are from
$\mathbb{F}$. For such an element we have
$|u|^{2}=N^{B}_{\mathbb{F}}(\lambda)-rs$ (this expression resembles the
Moore determinant of Hermitian matrices---see Corollary \ref{iso6d1} below), and
the expression from Proposition \ref{NAexp} indeed yields the square of the
latter expression. For the general case we embed $C$ in $M_{2}(\mathbb{K})$ and
$A$ in $M_{2}(B_{\mathbb{K}})$ as above and use extension of scalars. This
proves the corollary.
\end{proof}

Corollary \ref{NAvn2} emphasizes the fact that an element $u \in A^{-}$ is
invertible if and only if $|u|^{2}\neq0$. We also note that for every $a \in A$
we have $N^{A}_{\mathbb{F}}(a)=N^{A}_{\mathbb{F}}(\overline{a})$, either by a
direct evaluation using our formulae (and the fact that the same assertion holds
for the reduced norms of the quaternion algebras $B$ and $C$) or by Corollary
2.2 of \cite{[KMRT]}.

\smallskip

The bi-quaternion algebra $A$ operates on itself via $g:M \mapsto
gM\overline{g}$, and this action preserves the subspaces $A^{\pm}$. The
properties for this action of $A$ on the 6-dimensional space $A^{-}$ underlying
the Albert form which will be useful for our purposes are given in the
following
\begin{prop}
The action of $g \in A$ multiplies the norm $|u|^{2}$ of the element $u \in
A^{-}$ by $N^{A}_{\mathbb{F}}(g)$. The only invertible elements whose action on
$A^{-}$ is a global scalar multiplication are scalars from
$\mathbb{F}^{\times}$. \label{NAFg}
\end{prop}

\begin{proof}
We first consider the case where $C=M_{2}(\mathbb{F})$ and $A=M_{2}(B)$. The
element $u$ takes the form $\binom{\lambda\ \ -r}{s\ \ -\overline{\lambda}}$,
with $|u|^{2}=N^{B}_{\mathbb{F}}(\lambda)-rs$, as in the proof of Corollary
\ref{NAexp}. Now, for $g=\binom{a\ \ b}{c\ \ d} \in M_{2}(B)$ we have
$\overline{g}=\binom{\ \ \overline{d}\ \ -\overline{b}}{-\overline{c}\ \ \ \
\overline{a}}$, and then $gu\overline{g} \in M_{2}(B)^{-}$ involves the
quaternion
$a\lambda\overline{d}+sb\overline{d}+ra\overline{c}+b\overline{\lambda}\overline
{c}$ and the numbers
$rN^{B}_{\mathbb{F}}(a)+Tr^{B}_{\mathbb{F}}(a\lambda\overline{b})+sN^{B}_{
\mathbb{F}}(b)$ and
$sN^{B}_{\mathbb{F}}(d)+Tr^{B}_{\mathbb{F}}(c\lambda\overline{d})+rN^{B}_{
\mathbb{F}}(c)$ in the places of $r$ and $s$ respectively. Evaluating the norm
of the element involving these parameters yields the desired expression
\[\big[N^{B}_{\mathbb{F}}(a)N^{B}_{\mathbb{F}}(d)+N^{B}_{\mathbb{F}}(b)N^{B}_{
\mathbb{F}}(c)-Tr^{B}_{\mathbb{F}}(\overline{a}b\overline{d}c)\big]
\big(N^{B}_{\mathbb{F}}(\lambda)-rs\big).\] In addition, if $g$ is such that
these parameters are multiples of the original ones by a global scalar, then
$b\overline{d}=a\overline{c}=0$ and
$N^{B}_{\mathbb{F}}(b)=N^{B}_{\mathbb{F}}(c)=0$, and the fact that the trace
form on $B$ is non-degenerate implies also $\overline{b}a=\overline{d}c=0$.
Since invertible elements of $A$ have non-zero norm, we find that
$N^{B}_{\mathbb{F}}(a)$ and $N^{B}_{\mathbb{F}}(d)$ must not vanish, hence $a$
and $d$ are invertible and thus $b=c=0$. The global scalar property now implies
$N^{B}_{\mathbb{F}}(a)=N^{B}_{\mathbb{F}}(d)$, examining the effect on $\lambda$
being $d$ or $\overline{a}$ implies $a=d$, and the scalar multiplication
property shows that this element $a=d$ is central in $A$. Hence
$g\in\mathbb{F}^{\times}$ as desired, and the converse direction is trivial. .
In the case where $A$ is a division algebra, we can extend scalars to a
splitting field $\mathbb{K}$ of $C$, apply the arguments over $\mathbb{K}$, and
then return to our space over $\mathbb{F}$. This works since the multiplier lies
in $\mathbb{F}$, it is the reduced norm of an element of a bi-quaternion algebra
over $\mathbb{F}$, and the only elements in $A$ which become scalars in
$A_{\mathbb{K}}$ come from $\mathbb{F}$. This proves the proposition.
\end{proof}

Before we define the group which will become the Gspin group in this case, we
need another
\begin{lem}
The equality
$\widetilde{gu\overline{g}}=N^{A}_{\mathbb{F}}(g)\overline{g}^{-1}\tilde{u}g^{-1
}$ holds for any $g \in A^{\times}$ and $u \in A^{-}$. \label{AxA-rel}
\end{lem}

\begin{proof}
From Lemma \ref{vnorm6} and Proposition \ref{NAFg} we get
\[gu\overline{g}\cdot\widetilde{gu\overline{g}}=|\widetilde{gu\overline{g}}|^{2}
=N^{A}_{\mathbb{F}}(g)|u|^{2}=N^{A}_{\mathbb{F}}(g)g \cdot u\tilde{u} \cdot
g^{-1}=N^{A}_{\mathbb{F}}(g)gu\overline{g}\cdot\overline{g}^{-1}\tilde{u}g^{-1},
\] since $|u|^{2}=u\tilde{u}$ is a scalar (hence central in $A$). The assertion
now follows for every $u \in A^{-} \cap A^{\times}$, and extends to all of
$A^{-}$ by linearity and the fact that $A^{-}$ admits a basis of (orthogonal)
vectors which lie in $A^{\times}$. This proves the lemma.
\end{proof}

We consider the subgroup $A^{(\mathbb{F}^{\times})^{2}}$ of $A^{\times}$. This
group contains every single tensor $b \otimes c$ with $b \in B^{\times}$ and $c
\in C^{\times}$ (and in particular any scalar from $\mathbb{F}^{\times}$), and
Corollary \ref{NAvn2} shows that $A^{-} \cap A^{\times} \subseteq
A^{(\mathbb{F}^{\times})^{2}}$. Let $\widetilde{A}^{(\mathbb{F}^{\times})^{2}}$
denote the ``metaplectic-like'' double cover of $A^{(\mathbb{F}^{\times})^{2}}$
consisting of pairs $(g,t) \in
A^{(\mathbb{F}^{\times})^{2}}\times\mathbb{F}^{\times}$ such that
$N^{A}_{\mathbb{F}}(g)=t^{2}$ (with coordinate-wise product). This group appears
in Section 17 of \cite{[KMRT]} in relation with the orthogonal group of $O(V)$,
and we prove this connection using very simple means.
\begin{lem}
The operation $(g,t):u\mapsto\frac{gu\overline{g}}{t}$ defines a map
$\widetilde{A}^{(\mathbb{F}^{\times})^{2}} \to O(A^{-})$, with kernel consisting
of the elements $(r,r^{2})$ for $r\in\mathbb{F}^{\times}$. The automorphism
$g\mapsto\overline{g}^{-1}$ of $A^{\times}$ preserves
$A^{(\mathbb{F}^{\times})^{2}}$, and
$(g,t)\mapsto(t\overline{g}^{-1},t)$ is an automorphism of
$\widetilde{A}^{(\mathbb{F}^{\times})^{2}}$ of order 2. If $\tilde{\theta}$
generates a cyclic group of order 2 acting on
$\widetilde{A}^{(\mathbb{F}^{\times})^{2}}$ by this automorphism then sending it
to $\theta$ yields a map from the associated semi-direct product to $O(A^{-})$.
\label{ac6d1}
\end{lem}

\begin{proof}
The fact that the image of $\widetilde{A}^{(\mathbb{F}^{\times})^{2}}$ lies in
$O(A^{-})$ follows from Proposition \ref{NAFg} and the definition of
$\widetilde{A}^{(\mathbb{F}^{\times})^{2}}$. As for the kernel, Proposition
\ref{NAFg} implies that elements of the kernel must lie over
$\mathbb{F}^{\times}$, and the fact that they take the form $(r,r^{2})$ (and
that these elements indeed lie in $\widetilde{A}^{(\mathbb{F}^{\times})^{2}}$)
is immediate. The fact that $N^{A}_{\mathbb{F}}(\overline{g}^{\ -1})$ is the
reciprocal of $N^{A}_{\mathbb{F}}(g)$ shows that the automorphism preserves
$A^{(\mathbb{F}^{\times})^{2}}$, and the fact that
$N^{A}_{\mathbb{F}}(t)=t^{4}$ and $t\overline{t}^{\ -1}=1$ for
$t\in\mathbb{F}^{\times}$ proves the assertions about the automorphism of
$\widetilde{A}^{(\mathbb{F}^{\times})^{2}}$. As $\theta \in O(A^{-})$ (clear),
the result about the semi-direct product follows from Lemma \ref{AxA-rel} and
the fact that $t^{2}=N^{A}_{\mathbb{F}}(g)$ for
$(g,t)\in\widetilde{A}^{(\mathbb{F}^{\times})^{2}}$ . This proves the lemma.
\end{proof}
The multiplication by $t$ in the automorphism of
$\widetilde{A}^{(\mathbb{F}^{\times})^{2}}$ is not necessary. However, we put
it there for certain maps below (see Lemmas \ref{Qtheta} and \ref{Qhatpsi}) to
take a neater form.

As in the previous cases, we now prove
\begin{lem}
The reflection in the vector $g \in A^{-} \cap A^{\times}$ lies in the image of
the semi-direct product from Lemma \ref{ac6d1}, and takes the form
$u\mapsto\frac{g\tilde{u}\overline{g}}{|g|^{2}}$. \label{ref6d1}
\end{lem}

\begin{proof}
Corollary \ref{NAvn2} shows that
$(g,|g|^{2})\in\widetilde{A}^{(\mathbb{F}^{\times})^{2}}$, so that the asserted
operation lies in the image of the map from Lemma \ref{ac6d1}. Now,
$\overline{g}=-g$ since $g \in A^{-}$, and Lemma \ref{vnorm6} shows that the
denominator is $g\tilde{g}$. In the same manner as in the previous cases,
replacing $-g\tilde{u}$ by $u\tilde{g}$ for $u \in g^{\perp}$ and just
substituting $g$ for $u=g$ shows that this expression yields $u$ for $u \in
g^{\perp}$ case and $-g=-u$ for $u=g$. This proves the lemma.
\end{proof}

We remark that the operation from Lemma \ref{ref6d1} is indeed invariant under
interchanging the roles of $B$ and $C$, since then both $\theta$ and $|g|^{2}$
are being inverted.

We now come to prove
\begin{thm}
The Gspin group of $A^{-}$ is $\widetilde{A}^{(\mathbb{F}^{\times})^{2}}$, and
it is generated by the elements $(g,|g|^{2})$ with $g \in A^{-} \cap
A^{\times}$. The spin group is a subgroup mapping bijectively onto $A^{1}$.
\label{dim6d1}
\end{thm}

\begin{proof}
As in the previous cases, the fact that $\theta$ has determinant $-1$ (it
inverts a 3-dimensional subspace), Lemma \ref{ref6d1}, and Proposition \ref{CDT}
show that the map from the semi-direct product in Lemma \ref{ac6d1} to
$O(A^{-})$ is surjective, and its restriction to
$\widetilde{A}^{(\mathbb{F}^{\times})^{2}}$ maps surjectively onto $SO(A^{-})$,
with kernel (isomorphic to) $\mathbb{F}^{\times}$. Hence
$Gspin(A^{-})=\widetilde{A}^{(\mathbb{F}^{\times})^{2}}$, and from the structure
of the semi-direct product we deduce the generation of the latter group by
elements arising from $A^{-} \cap A^{\times}$. Note that here the latter
assertion is not so easy to establish in a direct manner. As $\theta$ has spinor
norm 1 (the 3-dimensional space which it inverts has determinant 1), Lemma
\ref{ref6d1} implies that the spinor norm of the element $(g,|g|^{2}=t)$ is
$t=|g|^{2}$. As these elements generate
$\widetilde{A}^{(\mathbb{F}^{\times})^{2}}$, multiplicativity shows that the
spinor norm of every every element
$(g,t)\in\widetilde{A}^{(\mathbb{F}^{\times})^{2}}$ is just $t$ (in particular,
the element $(1,-1)$, operating as $-Id_{A^{-}}$, has the required spinor norm
$-1$ like the determinant of $A^{-}$). Thus, the elements lying over
$SO^{1}(A^{-})$ come with second coordinate from $(\mathbb{F}^{\times})^{2}$,
and as $(r,r^{2})$ operates trivially for every $r\in\mathbb{F}^{\times}$, we
may divide by it and consider just elements of
$\widetilde{A}^{(\mathbb{F}^{\times})^{2}}$ of the form $(g,1)$. As this group
maps isomorphically onto $A^{1}$ in the natural projection
$\widetilde{A}^{(\mathbb{F}^{\times})^{2}} \to A^{(\mathbb{F}^{\times})^{2}}$,
and its intersection with the image of the kernel $\mathbb{F}^{\times}$ is just
$\pm1$, we find that the spin group $spin(A^{-})$ is just (this isomorphic image
of) $A^{1}$. This proves the theorem.
\end{proof}

Observe that while the space depends on the choice of the involution on $A$
(i.e., the decomposition as a tensor product of quaternion algebras), the groups
$SO(A^{-})$ and $SO^{1}(A^{-})$ indeed do not.

For the isotropic case, we have
\begin{cor}
If $A^{-}$ is isotropic then there is a quaternion algebra $B$ such that the
Gspin group is the double cover of the group
$\widetilde{GL}_{2}^{(\mathbb{F}^{\times})^{2}}(B)$ consisting of those matrices
in which the expression from Proposition \ref{NAFg} is a square, where the
double cover involves choosing a square root for it. The spin group is just the
subgroup $GL_{2}^{1}(B)$ in which this expression equals 1. If $A^{-}$ splits
another hyperbolic plane then it is the sum of three such planes, the Gspin
group is $\widetilde{GL}_{4}^{(\mathbb{F}^{\times})^{2}}(\mathbb{F})$, and the
spin group is $SL_{4}(\mathbb{F})$. \label{iso6d1}
\end{cor}

\begin{proof}
The fact that $A$ cannot be a division algebra follows immediately from the fact
that some elements of $A^{-}$ are zero-divisors when $A^{-}$ is isotropic.
Moreover, in this case $A^{-}$ splits a hyperbolic plane, so that by choosing
the subspace to normalize as the $B_{0}$ part to be a subspace of the orthogonal
complement of this hyperbolic plane we may assume (see Corollary \ref{iso3})
that $C=M_{2}(\mathbb{F})$. Hence $A=M_{2}(B)$, the space $A^{-}$ consists of
the matrices $\binom{\lambda\ \ -r}{s\ \ -\overline{\lambda}}$ from the proof of
Corollary \ref{NAvn2}, and
$A^{(\mathbb{F}^{\times})^{2}}=GL_{2}^{(\mathbb{F}^{\times})^{2}}(B)$ is defined
according to the expression from Proposition \ref{NAFg}. Thus, Theorem
\ref{dim6d1} shows that $Gspin\big(M_{2}(B)^{-}\big)$ is the double cover
$\widetilde{GL}_{2}^{(\mathbb{F}^{\times})^{2}}(B)$, and
$spin\big(M_{2}(B)^{-}\big)$ is $A^{1}=GL_{2}^{1}(B)$ according to the same
expression. Observe that the complement of a hyperbolic plane (the $r$ and $s$
coordinates) is the space $B$ from Corollary \ref{dim4d1}. Hence $A^{-}$ splits
another hyperbolic plane if and only if $B$ also splits (see Corollary
\ref{iso4}), and as in this case the complement of these two planes has
discriminant 1, it is again a hyperbolic plane by Corollary \ref{iso2}. Hence
$A=M_{4}(\mathbb{F})$, the reduced norm is the determinant, and the Gspin and
spin groups from Theorem \ref{dim6d1} become indeed
$\widetilde{GL}_{4}^{(\mathbb{F}^{\times})^{2}}(\mathbb{F})$ and
$SL_{4}(\mathbb{F})$ respectively. This proves the corollary.
\end{proof}

Theorem 16.5 of \cite{[KMRT]} relates the index of $A$ (4 in case of a division
algebra, 1 in the split case, 2 in the middle) to the number of hyperbolic
planes one can split from the Albert form of $A$ (0, 3, and 1 respectively).
The hardest direction is to prove that if the Albert form is anisotropic then
$A$ is a division algebra, yielding the inverse direction of Corollary
\ref{iso6d1}. On the other hand, all the remaining assertions are either trivial
or follow from the proof of Corollary \ref{iso6d1}.

We also have
\begin{cor}
Given a quaternion algebra $B$, the groups
$\widetilde{GL}_{2}^{(\mathbb{F}^{\times})^{2}}(B)$ and $GL_{2}^{1}(B)$ are also
the Gspin and spin groups of the space $M_{2}^{Her}(B)$ of Hermitian $2\times2$
matrices over $B$, which is quadratic with the vector norm being (minus) the
Moore determinant, through the action via
$g:X\mapsto\frac{gX\iota_{B}(g)^{t}}{t}$. In addition, the groups
$\widetilde{GL}_{4}^{(\mathbb{F}^{\times})^{2}}(\mathbb{F})$ and
$SL_{4}(\mathbb{F})$ are the Gspin and spin groups of $M_{4}^{as}(\mathbb{F})$
the space of anti-symmetric matrices over $\mathbb{F}$ (with minus the Pfaffian
as the quadratic form) on which they operate as $g:T\mapsto\frac{gTg^{t}}{t}$,
as well as on another space, with the operation $\big(g=\binom{a\ \ b}{c\ \
d},t\big):S\mapsto\binom{a\ \ b}{c\ \ d}S\binom{\ \ d^{t}\ \ -b^{t}}{-c^{t}\ \ \
\ a^{t}}/t$. \label{alt6d1}
\end{cor}

\begin{proof}
The first assertion follows, as in Corollaries \ref{alt3} and \ref{alt4}, from
Corollary \ref{iso6d1}, Lemma \ref{Sadjt}, and the fact that $M_{2}^{Her}(B)$ is
the image of $M_{2}(B)^{-}$ under right multiplication by $\binom{\ \ 0\ \
1}{-1\ \ 0}$. Indeed, the expression $N^{B}_{\mathbb{F}}(\lambda)-rs$ becomes
minus the Moore determinant of $\binom{r\ \ \lambda}{\overline{\lambda}\ \ s}$.
For the second assertion, we may apply Lemma \ref{Sadjt} again, but now for the
elements of $B=M_{2}(\mathbb{F})$ which are entries of matrices in $A=M_{2}(B)$.
The representation just described becomes $M_{4}^{as}(\mathbb{F})$ under this
operation: The diagonal blocks $rI$ and $sI$ are taken to multiples of the
anti-symmetric matrix $\binom{\ \ 0\ \ 1}{-1\ \ 0}$, the off-diagonal blocks are
taken care of by Lemma \ref{Sadjt} and the fact that $\binom{\ \ 0\ \ 1}{-1\ \
0}$ is anti-symmetric, and the Moore determinant becomes the Pfaffian. The
fourth representation is obtained from the one considered in Theorem
\ref{dim6d1} by this right multiplication. This proves the corollary.
\end{proof}

The representation on $M_{4}^{as}(\mathbb{F})$ uses an orthogonal involution,
while the two others use symplectic ones (compare with Proposition 2.23(1) of
\cite{[KMRT]} again). The generator $\theta$ of $O(A^{-})/SO(A^{-})$ corresponds
to the adjoint involution of $2\times2$ matrices in Theorem \ref{dim6d1} in case
$A=M_{2}(B)$, and to $X\mapsto-X^{t}$ on the space $M_{2}^{Her}(B)$ from
Corollary \ref{alt6d1}. On the other hand, conjugating by $\binom{\ \ 0\ \
1}{-1\ \ 0}$ sends the latter involution to $X\mapsto-adjX$ on $M_{2}^{Her}(B)$,
and the product of $X$ and its image under this involution yields $|X|^{2}$
again. $\theta$ yields involutions on $M_{4}^{as}(\mathbb{F})$ and on the other
space as well, and after appropriate conjugations (by $\binom{\ \ 0\ \ 1}{-1\ \
0}$ on the entries as $2\times2$ matrices, or on both), we see that these spaces
also admit involutions with the property that the product of a vector with its
image yields the vector norm. The latter involution on $M_{4}^{as}(\mathbb{F})$
will be denoted $T\mapsto\hat{T}$. 

\section{Dimension 5 \label{Dim5}}

The spaces we get in dimension 5 are those given in
\begin{lem}
Every 5-dimensional space is isometric, up to a scalar multiple, to the
orthogonal complement of an anisotropic vector $Q$ in the space $A^{-}$ arising
from a bi-quaternion algebra $A$ presented as $B \otimes C$. \label{sp5}
\end{lem}

\begin{proof}
First we remark that for anisotropic $Q \in A^{-}$, the subspace $Q^{\perp}$ of
$A^{-}$ is 5-dimensional and non-degenerate. Conversely, a 5-dimensional
quadratic space can be can be extended, uniquely up to the choice of a
generator, to a 6-dimensional space of discriminant 1: If the space has some
discriminant (or determinant) $d$, this is done by adding an element $Q$ with
$|Q|^{2}=-d$. The lemma now follows directly from Lemma \ref{sp6d1}.
\end{proof}

In fact, by choosing the 3-dimensional subspace from the proof of Lemma
\ref{sp6d1} to be contained in our original space $Q^{\perp}$ we make sure that
$Q \in 1 \otimes C_{0}$ (see Proposition \ref{NQF2NA} below for a more precise
statement), but we shall not need this fact. In any case, we shall write our
vector space as $Q^{\perp} \subseteq A^{-}$. The next step is
\begin{lem}
The groups $SO(Q^{\perp} \subseteq A^{-})$, $SO^{1}(Q^{\perp} \subseteq A^{-})$,
$Gspin(Q^{\perp} \subseteq A^{-})$, and $spin(Q^{\perp} \subseteq A^{-})$ are
the subgroups of $SO(A^{-})$, $SO^{1}(A^{-})$, $Gspin(A^{-})$, and $spin(A^{-})$
respectively, consisting of those elements whose action stabilizes $Q$.
\label{ac5}
\end{lem}

\begin{proof}
The Witt Cancelation Theorem implies that any element of $O(Q^{\perp} \subseteq
A^{-})$ comes from an element of $O(A^{-})$ under which $Q^{\perp}$ is
invariant, and $\mathbb{F}Q$ must therefore also be invariant under such
extensions. As $O(\mathbb{F})=\{\pm1\}$, and $-1$ has determinant $-1$ there,
the assertion about $SO(Q^{\perp} \subseteq A^{-})$ follows. Considering
extensions of $O(Q^{\perp} \subseteq A^{-})$ to $O(A^{-})$ by taking only $+1$
on $\mathbb{F}Q$, Proposition \ref{CDT} shows that any element there has the
same spinor norm in both $O(Q^{\perp} \subseteq A^{-})$ and $O(A^{-})$ (since
this is true for reflections in vectors from $Q^{\perp}$). The assertion for
$SO^{1}(Q^{\perp} \subseteq A^{-})$ follows, and for the remaining two groups we
just use the maps into $O(A^{-})$. This proves the lemma.
\end{proof}

For a simpler description, consider the group $A^{\times}_{\mathbb{F}Q}$ of
those $g \in A^{\times}$ such that $gQ\overline{g}\in\mathbb{F}Q$ (i.e., those
which preserve the 1-dimensional space $\mathbb{F}Q$), which comes with a group
homomorphism $m:A^{\times}_{\mathbb{F}Q}\to\mathbb{F}^{\times}$ defined
by $gQ\overline{g}=t(g)Q$. We remark that $A^{\times}_{\mathbb{F}Q}$ can be seen
as the group of similitudes of $A$ with respect to the involution $x \mapsto
Q\overline{x}Q^{-1}$ (which is symplectic since $Q \in A^{-}$). We now have
\begin{lem}
The group $A^{\times}_{\mathbb{F}Q}$ is a subgroup of
$A^{(\mathbb{F}^{\times})^{2}}$. The double cover
$\widetilde{A}^{(\mathbb{F}^{\times})^{2}}$ splits over
$A^{\times}_{\mathbb{F}Q}$, and the splitting group contains the group
$\{(r,r^{2})|r\in\mathbb{F}^{\times}\}$. \label{AFQAF2}
\end{lem}

\begin{proof}
The equality $gQ\overline{g}=t(g)Q$ holding for $g \in A^{\times}_{\mathbb{F}Q}$
implies, by Proposition \ref{NAFg} and the fact that $|Q|^{2}\neq0$, the
equality $N^{A}_{\mathbb{F}}(g)=t(g)^{2}$. This establishes the first assertion,
as well as introducing the splitting homomorphism $g\mapsto\big(g,t(g)\big)$
from $A^{\times}_{\mathbb{F}Q}$ into $\widetilde{A}^{(\mathbb{F}^{\times})}$.
The fact that $rQ\overline{r}=r^{2}Q$ for $r\in\mathbb{F}^{\times}$ completes
the proof of the  lemma.
\end{proof}

We can now prove
\begin{thm}
We have $Gspin(Q^{\perp} \subseteq A^{-})=A^{\times}_{\mathbb{F}Q}$, and
$spin(Q^{\perp} \subseteq A^{-})$ is the group $A^{1}_{Q}$ of elements of $g
\in A^{1}$ such that $gQ\overline{g}=Q$. \label{dim5}
\end{thm}

\begin{proof}
For the Gspin group, Lemma \ref{ac5} shows that it suffices to prove that the
image of $A^{\times}_{\mathbb{F}Q}$ under the splitting map from Lemma
\ref{AFQAF2} is the stabilizer of $Q$ under the action of
$\widetilde{A}^{(\mathbb{F}^{\times})}$. But it is clear from the definition
that the image of this splitting map stabilizes $Q$, and that elements of this
stabilizer must come from $A^{\times}_{\mathbb{F}Q}$. As replacing $m$ by $-m$
yields an element taking $Q$ to $-Q$, the assertion for the Gspin group follows.
The statement for the spin group follows directly from Lemma \ref{ac5}, since no
scalar appears in the action of $(g,1)$ for $A^{1}$. Note that for $g \in
A^{1}_{Q}$ we have $t(g)=1$, so that the corresponding element is indeed of the
form $(g,1)$ for $A^{1}$. This proves the proposition.
\end{proof}

Note that our construction is based on the choice of $Q \in A^{-}$. However, the
only parameter which is required to know the isomorphism class of
$A^{\times}_{\mathbb{F}Q}$ and $A^{1}_{Q}$ is given in the following
\begin{prop}
If an element $R \in A^{-}$ satisfies
$|R|^{2}\in|Q|^{2}(\mathbb{F}^{\times})^{2}N^{A}_{\mathbb{F}}(A^{\times})$ then
$A^{\times}_{\mathbb{F}Q} \cong A^{\times}_{\mathbb{F}R}$ and $A^{1}_{Q} \cong
A^{1}_{R}$. \label{NQF2NA}
\end{prop}

\begin{proof}
First note that $A^{\times}_{\mathbb{F}rQ}=A^{\times}_{\mathbb{F}Q}$ and
$A^{1}_{rQ}=A^{1}_{Q}$ for every $r\in\mathbb{F}^{\times}$, so once only
(non-zero) vector norms are involved, we can consider them modulo
$(\mathbb{F}^{\times})^{2}$. Now, for any $a \in A^{\times}$, conjugation by $a$
takes $A^{\times}_{\mathbb{F}Q}$ to $A^{\times}_{\mathbb{F}aQ\overline{a}}$ and
$A^{1}_{Q}$ to $A^{1}_{aQ\overline{a}}$. Applying this to $(a,t)$ with $a \in
A^{(\mathbb{F}^{\times})^{2}}$ and using the transitivity of the action of
$SO(A^{-})$ on elements of the same vector norm (Witt Cancelation again)
establishes the assertion in case $|R|^{2}=|Q|^{2}$ (hence also if
$|R|^{2}\in|Q|^{2}(\mathbb{F}^{\times})^{2}$), and then doing so for general $a
\in A^{\times}$ allows us to divide also by $N^{A}_{\mathbb{F}}(A^{\times})$.
This proves the proposition.
\end{proof}

When we consider the case where $A^{-}$ is isotropic, the relation from
Proposition \ref{NQF2NA} becomes simpler by the following
\begin{lem}
If $A=M_{2}(B)$ then $N^{A}_{\mathbb{F}}(A)=N^{B}_{\mathbb{F}}(B)$ and
$N^{A}_{\mathbb{F}}(A^{\times})=N^{B}_{\mathbb{F}}(B^{\times})$. \label{NM2B}
\end{lem}

\begin{proof}
For any $b \in B$ we get $N^{B}_{\mathbb{F}}(b)$ as the reduced norm of the
element $\binom{b\ \ 0}{0\ \ 1}$ of $A$. For the other direction, the proof of
Proposition \ref{NAFg} shows that if a matrix in $M_{2}(B)$ has an invertible
entry then its reduced norm is the product of two norms from $B$ (e.g.,
$N^{B}_{\mathbb{F}}(a)N^{B}_{\mathbb{F}}(d-ca^{-1}b)$ if $a \in B^{\times}$),
which completes the proof for division algebras $B$ since $N^{B}_{\mathbb{F}}$
is multiplicative. As for $B=M_{2}(\mathbb{F})$ and $A=M_{4}(\mathbb{F})$, every
element of $\mathbb{F}$ is both a $2\times2$ determinant and a $2\times2$
determinant, the proof of the lemma is complete.
\end{proof}

Rather than isotropy, in this case we have
\begin{cor}
If the 5-dimensional space of discriminant $d$ represents $-d$ then its Gspin
and spin groups are isomorphic to $GSp_{B}\binom{-\delta\ \ 0}{\ \ 0\ \ 1}$ and
$Sp_{B}\binom{-\delta\ \ 0}{\ \ 0\ \ 1}$ respectively, for some quaternion
algebra $B$ and some $\delta\in\mathbb{F}^{\times}$ which may be taken from a
set of representatives for $F^{\times}/N^{B}_{\mathbb{F}}(B^{\times})$. In case
this space splits two hyperbolic planes, these groups are isomorphic to the
classical groups $GSp_{4}(\mathbb{F})$ and $Sp_{4}(\mathbb{F})$ respectively.
\label{iso5}
\end{cor}

\begin{proof}
In this case (which covers the case of an isotropic space since a hyperbolic
plane represents every element of $\mathbb{F}$) the space $A^{-}$ from the
proof of Lemma \ref{sp5} will be isotropic. Hence we can assume that
$A=M_{2}(B)$ (with the associated involution) by Corollary \ref{iso6d1}, and the
Gspin and spin groups from Theorem \ref{dim5} are $GL_{2}(B)_{\mathbb{F}Q}$ and
$GL_{2}^{1}(B)_{Q}$ respectively. Lemma \ref{NQF2NA} shows that the only
parameter required for determining the isomorphism classes of these groups is
$|Q|^{2}$ up to $(\mathbb{F}^{\times})^{2}N^{A}_{\mathbb{F}}(A^{\times})$, and
by Lemma \ref{NM2B} and the fact that $(\mathbb{F}^{\times})^{2} \subseteq
N^{B}_{\mathbb{F}}(B^{\times})$ for quaternion algebras (as
$N^{B}_{\mathbb{F}}(r)=r^{2}$), we may take the vector norm $\delta$ from the
required set of representatives and any $Q$ with $|Q|^{2}=\delta$ will do. We
choose $Q$ to be the element $\binom{0\ \ \delta}{1\ \ 0}$ of
$M_{2}(\mathbb{F})_{0}=1 \otimes C_{0} \subseteq M_{2}(B)^{-}$, with
$|Q|^{2}=-\det Q=\delta$. Using Corollary \ref{alt6d1}, the fact that
multiplying $Q$ by $\binom{\ \ 0\ \ 1}{-1\ \ 0}$ from the right gives this
diagonal matrix shows that our groups are indeed $GSp_{B}\binom{-\delta\ \ 0}{\
\ 0\ \ 1}$ and $Sp_{B}\binom{-\delta\ \ 0}{\ \ 0\ \ 1}$.

In case our space contains two hyperbolic planes, Corollary \ref{iso6d1} shows
that $B$ also splits. As $N^{B}_{\mathbb{F}}(B^{\times})=\mathbb{F}^{\times}$
in this case, we have, up to isomorphism, the same group for every choice of $Q$
(Lemma \ref{NQF2NA}). Instead of taking $Q$ as above, we recall from the proof
of Corollary \ref{NAvn2} that the space $A^{-}$ (now for split $B$ and $C$)
consists of matrices of the form $\binom{\ Y\ \ \ -rI\ \ }{sI\ \ -adjY}$ with
$r$ and $s$ from $\mathbb{F}$ and $Y \in M_{2}(\mathbb{F})$, and we take $Q$ to
be the element with $r=s=0$ and $Y=\binom{\ \ 0\ \ 1}{-1\ \ 0}$. Going over to
the representation from Corollary \ref{iso6d1} in which
$\widetilde{GL}_{4}^{(\mathbb{F}^{\times})^{2}}(\mathbb{F})$ operates by
$g:T\mapsto\frac{gTg^{t}}{t}$, our groups are easily seen to be the classical
$GSp_{4}(\mathbb{F})$ and $Sp_{4}(\mathbb{F})$. This proves the corollary.
\end{proof}

As for equivalent representations, we get
\begin{cor}
The groups $GSp_{B}\binom{-\delta\ \ 0}{\ \ 0\ \ 1}$ and $Sp_{B}\binom{-\delta\
\ 0}{\ \ 0\ \ 1}$ are the Gspin and spin groups of the orthogonal complement of
$\binom{-\delta\ \ 0}{\ \ 0\ \ 1}$ inside $M_{2}^{Her}(B)$ from Corollary
\ref{alt6d1}. The classical groups $GSp_{4}(\mathbb{F})$ and
$Sp_{4}(\mathbb{F})$ can be seen as the Gspin and spin group of either the
complement of $\binom{0\ \ -I}{I\ \ \ \ 0}$ in $M_{4}^{as}(\mathbb{F})$ or of
the orthogonal complement of the adjoint representation on the Lie algebra
$\mathfrak{gsp}_{4}(\mathbb{F})$ of $GSp_{4}(\mathbb{F})$ inside
$M_{4}(\mathbb{F})=\mathfrak{gl}_{4}(\mathbb{F})$.
\label{alt5}
\end{cor}

\begin{proof}
This follows directly by restricting the representations from Corollary
\ref{alt6d1} to our groups. Note that the operation $g=\binom{a\ \ b}{c\ \
d}:S\mapsto\binom{a\ \ b}{c\ \ d}S\binom{\ \ d^{t}\ \ -b^{t}}{-c^{t}\ \ \ \
a^{t}}$ is conjugation for $Sp_{4}(\mathbb{F})$ (and conjugation tensored with
the determinant for $GSp_{4}(\mathbb{F})$), leaving the Lie algebra
$\mathfrak{sp}_{4}(\mathbb{F})$ invariant. The additional invariant vector is
$I$, adding which to $\mathfrak{sp}_{4}(\mathbb{F})$ yields
$\mathfrak{gsp}_{4}(\mathbb{F})$. This proves the corollary.
\end{proof}

We remark that the simplicity of $Sp_{4}(\mathbb{F})$ as an algebraic group
implies that the action on $\mathfrak{sp}_{4}(\mathbb{F})$ is an irreducible
representation of $Sp_{4}(\mathbb{F})$. Hence Corollary \ref{alt5} yields a
complete reduction of the representation $\mathfrak{gl}_{4}(\mathbb{F})$ of
these groups as the direct sum of $\mathfrak{sp}_{4}(\mathbb{F})$,
$\mathbb{F}I$ (as the determinant), and the 5-dimensional representation of
$Sp_{4}(\mathbb{F})$ as an $SO$ or $SO^{1}$ group. The adjoint representation
appears for the group $Sp_{B}\binom{-\delta\ \ 0}{\ \ 0\ \ 1}$ in Corollary
\ref{alt5} only if $\delta=-1$.

\section{Dimension 6, General Discriminant \label{Dim6gen}}

If our 6-dimensional quadratic space over $\mathbb{F}$ now has some discriminant
$d$, take $\mathbb{E}=\mathbb{F}(\sqrt{d}))$ with Galois automorphism $\rho$ as
before. Our object of interest will be bi-quaternion algebras over $\mathbb{E}$,
which take the form $A_{\mathbb{E}}$ for some bi-quaternion algebra $A$ over
$\mathbb{F}$. Given a presentation of $A$ as $B \otimes C$ as in Lemma
\ref{sp6d1}, both the orthogonal involution $\iota_{B}\otimes\iota_{C} \otimes
Id_{\mathbb{E}}:x\mapsto\overline{x}$ and the unitary involution
$\iota_{B}\otimes\iota_{C}\otimes\rho:x\mapsto\overline{x}^{\rho}$ are defined
on $A_{\mathbb{E}}$. The space $A_{\mathbb{E}}^{-}$ from Lemma \ref{sp6d1} is
defined as a quadratic space over $\mathbb{E}$, and we shall consider the vector
space over $\mathbb{F}$ which is defined in the following
\begin{lem}
Take some anisotropic $Q \in A^{-}$. The set of elements $u \in
A_{\mathbb{E}}^{-}$ satisfying the equality
$u^{\rho}=-\frac{Q\tilde{u}Q}{|Q|^{2}}$ form a 6-dimensional quadratic space of
discriminant $d$ over $\mathbb{F}$. Every quadratic space of dimension 6 and
discriminant $d$ over $\mathbb{F}$ is isomorphic, up to rescaling, to a space
which is obtained in this way. \label{sp6gen}
\end{lem}

\begin{proof}
The proof of Lemma \ref{ref6d1} shows that the expression, which we require to
equal $u^{\rho}$, is just $u$ if $u \in Q^{\perp} \subseteq
A_{\mathbb{E}}^{-}$, and is $-cQ$ in case $u=cQ$. As $Q^{\rho}=Q$, it follows
that the elements $u$ which satisfy this property is precisely $(Q^{\perp}
\subseteq A^{-})\oplus\mathbb{E}_{0}Q$, which is the quadratic space over
$\mathbb{F}$ which is obtained from $A^{-}$ by rescaling the vector norm of $Q$
by $d$. It is thus 6-dimensional and of the required discriminant. Conversely,
given a space of dimension 6 and discriminant $d$, choose an arbitrary isotropic
vector, and divide its vector norm by $d$. The resulting space has discriminant
1, so that by Lemma \ref{sp6d1} it takes the form $A^{-}$ for some bi-quaternion
algebra $A=B \otimes C$ over $\mathbb{F}$. But then we may extend scalars to
$\mathbb{E}$, and multiplying $|Q|^{2}$ back by $d$ (to get its original value)
means precisely replacing the $\mathbb{F}$-subspace $A^{-}$ of
$A_{\mathbb{E}}^{-}$ by the space $(Q^{\perp} \subseteq
A^{-})\oplus\mathbb{E}_{0}Q$ considered above. This completes the proof of the
lemma.
\end{proof}

We denote the space from Lemma \ref{sp6gen} by $(A_{\mathbb{E}}^{-})_{\rho,Q}$,
and observe that extending its scalars to $\mathbb{E}$ also gives
$A_{\mathbb{E}}^{-}$. The formulae from Lemma \ref{vnorm6} hold for $u \in
(A_{\mathbb{E}}^{-})_{\rho,Q}$, since such elements lie in $A_{\mathbb{E}}^{-}$
and the expressions are $\mathbb{F}$-valued on $(A_{\mathbb{E}}^{-})_{\rho,Q}$.
Note that $\rho$ acts on $(A_{\mathbb{E}}^{-})_{\rho,Q}$ like the reflection in
a generator of $\mathbb{E}_{0}Q$.

\smallskip

The group $A_{\mathbb{E}}^{\times}$ operates on $A$ either through the map from
Lemma \ref{ac6d1} (preserving $A_{\mathbb{E}}^{-}$), or via $g:M \mapsto
gM\overline{g}^{\rho}$. We define $A_{\mathbb{E},\rho,\mathbb{F}Q}^{\times}$ to
be the subgroup of $A_{\mathbb{E}}^{\times}$ which stabilizes the 1-dimensional
vector space $\mathbb{F}Q$ under the latter action (i.e., multiplies $Q$ by a
scalar), and let
$t:A_{\mathbb{E},\rho,\mathbb{F}Q}^{\times}\to\mathbb{F}^{\times}$ be the group
homomorphism defined via $gQ\overline{g}^{\rho}=t(g)Q$. Note that as
$\overline{Q}^{\rho}=-Q$ and $gQ\overline{g}^{\rho}$ lies in the same eigenspace
of $\iota_{B}\otimes\iota_{C}\otimes\rho$, allowing $t(g)$ to be in
$\mathbb{E}^{\times}$ does not produce a larger group. This group can be seen as
the group of similitudes of the unitary involution $x \mapsto
Q\overline{x}^{\rho}Q^{-1}$. It is stable under $\rho$ (just apply $\rho$ to
the defining equation), and its automorphism $\rho$ commutes with the map $t$
into $\mathbb{F}^{\times}$. The group which we shall consider here is the one
appearing in the following
\begin{lem}
The group $A_{\mathbb{E},\rho,\mathbb{F}Q}^{t^{2}}$ consisting of those elements
$g \in A_{\mathbb{E},\rho,\mathbb{F}Q}^{\times}$ such that
$N^{A_{\mathbb{E}}}_{\mathbb{E}}(g)=t(g)^{2}$ has index either 1 or 2 in
$A_{\mathbb{E},\rho,\mathbb{F}Q}^{\times} \cap
A_{\mathbb{E}}^{\mathbb{F}^{\times}}$. The former group is stable under $\rho$,
and it coincides with $A_{\mathbb{E},\rho,\mathbb{F}Q}^{\times} \cap
A_{\mathbb{E}}^{(\mathbb{F}^{\times})^{2}}$ unless
$-1\in(\mathbb{F}^{\times})^{2}$ and the above index is 2. \label{ind2int}
\end{lem}

\begin{proof}
Given $g \in A_{\mathbb{E},\rho,\mathbb{F}Q}^{\times}$, the norm
$N^{A_{\mathbb{E}}}_{\mathbb{E}}(g)$ lies in $t(g)^{2}\mathbb{E}^{1}$. For
$g \in A_{\mathbb{E}}^{\mathbb{F}^{\times}}$ the second multiplier comes from
$\mathbb{E}^{1} \cap \mathbb{F}=\{\pm1\}$, so that the group
$A_{\mathbb{E},\rho,\mathbb{F}Q}^{t^{2}}$ equals
$A_{\mathbb{E},\rho,\mathbb{F}Q}^{\times} \cap
A_{\mathbb{E}}^{\mathbb{F}^{\times}}$, unless
$A_{\mathbb{E},\rho,\mathbb{F}Q}^{\times}$ contains elements $g$ with
$N^{A_{\mathbb{E}}}_{\mathbb{E}}(g)=-t(g)^{2}$, a case in which
$A_{\mathbb{E},\rho,\mathbb{F}Q}^{t^{2}}$ has index 2 in this intersection. The
stability under $\rho$ follows from the fact that for $g \in
A_{\mathbb{E},\rho,\mathbb{F}Q}^{t^{2}}$ we have $t(g)=t(g^{\rho})$ and
$N^{A}_{\mathbb{F}}(g)=t(g)^{2}$ lies in $\mathbb{F}$. The assertion about the
intersection with $A_{\mathbb{E}}^{(\mathbb{F}^{\times})^{2}}$ follows
immediately from the considerations concerning
$A_{\mathbb{E}}^{\mathbb{F}^{\times}}$. This proves the lemma.
\end{proof}

Lemma \ref{ind2int} is related to some delicate facts about the definition of
$GSU$ groups---see also Corollary \ref{iso6gen} below. Note that as
$N^{A_{\mathbb{E}}}_{\mathbb{E}}$ is of degree 4, Hilbert's Theorem 90 cannot
help us in getting a more accurate result, like the one following Lemma
\ref{ac4}. In any case, we have $A_{\mathbb{E},\rho,\mathbb{F}Q}^{t^{2}}
\subseteq A_{\mathbb{E}}^{(\mathbb{F}^{\times})^{2}} \subseteq
A_{\mathbb{E}}^{(\mathbb{E}^{\times})^{2}}$ (by Lemma \ref{ind2int}), and the
double cover $\widetilde{A}_{\mathbb{E}}^{(\mathbb{E}^{\times})^{2}}$ splits
over $A_{\mathbb{E},\rho,\mathbb{F}Q}^{t^{2}}$---indeed, the map
$g\mapsto\big(g,t(g)\big)$ is, by definition, a splitting map.

Note that unless $Q^{2}\in\mathbb{F}$ (i.e., $\tilde{Q}\in\mathbb{F}Q$, which
means that $Q$ lies in either $B_{\mathbb{E}}\otimes1$ or $1 \otimes
C_{\mathbb{E}}$), the map $\theta$ from Lemma \ref{vnorm6} does not preserve
$(A_{\mathbb{E}}^{-})_{\rho,Q}$ (see Lemma \ref{QthetaQrel} below). Moreover,
the automorphism from Lemma
\ref{ac6d1} does not preserve $A_{\mathbb{E},\rho,\mathbb{F}Q}^{t^{2}}$, since
inverting $Q$ may change the 1-dimensional space it generates. However, in order
to carry out the usual considerations also in this case, we have the following
\begin{lem}
The element $(Q,|Q|^{2})\tilde{\theta}$ squares to $(-|Q|^{2},|Q|^{4})$ in the
semi-direct product of $\{1,\tilde{\theta}\}$ and
$\widetilde{A}_{\mathbb{E}}^{^{(\mathbb{E}^{\times})^{2}}}$. It operates on
$(A_{\mathbb{E}}^{-})_{\rho,Q}$ as $\rho$ (hence preserves this space), and it
takes $(g,t) \in \widetilde{A}_{\mathbb{E}}^{^{(\mathbb{E}^{\times})^{2}}}$
to $(tQ\overline{g}^{-1}Q^{-1},t)$ by conjugation. This is an order 2
automorphism of $\widetilde{A}_{\mathbb{E}}^{^{(\mathbb{E}^{\times})^{2}}}$,
which preserves the image of $A_{\mathbb{E},\rho,\mathbb{F}Q}^{t^{2}}$ under the
splitting map, and coincides with the action of $\rho$ on the latter group.
\label{Qtheta}
\end{lem}

\begin{proof}
The product rule of the semi-direct product from Lemma \ref{ac6d1} implies that
the square of the element in question is
$(Q\cdot|Q|^{2}\overline{Q}^{-1},|Q|^{2}\cdot|Q|^{2})$, which equals
$(-|Q|^{2},|Q|^{4})$ since $Q \in A_{\mathbb{E}}^{-}$. This element sends $u \in
A_{\mathbb{E}}^{-}$ to $-\frac{Q\tilde{u}Q}{|Q|^{2}}$ (use the fact that $Q \in
A_{\mathbb{E}}^{-}$ again), which for $u \in (A_{\mathbb{E}}^{-})_{\rho,Q}$
coincides with $u^{\rho}\in(A_{\mathbb{E}}^{-})_{\rho,Q}$ according to Lemma
\ref{sp6gen}. The formula for the conjugation follows directly from Lemma
\ref{ac6d1}, and the fact that it is an automorphism of order 2 either follows
from the centrality of $(-1,1)$ or can be easily verified directly using the
fact that $Q \in A_{\mathbb{E}}^{-}$ once more. Now, multiplying the equation
stating that an element $g$ lies  in $A_{\mathbb{E},\rho,\mathbb{F}Q}^{\times}$
by $(gQ)^{-1}$ from the left yields the equality
$\overline{g}^{\rho}=t(g)Q^{-1}g^{-1}Q$, and after applying
$\iota_{B}\otimes\iota_{C}$ we get the equality
$g^{\rho}=t(g)Q\overline{g}^{\ -1}Q^{-1}$ (since $Q \in A_{\mathbb{E}}^{-}$).
This shows that conjugation by $(Q,|Q|^{2})\tilde{\theta}$ operates on
$A_{\mathbb{E},\rho,\mathbb{F}Q}^{t^{2}}$ as $\rho$, and as the latter group was
seen to be preserved by $\rho$ in Lemma \ref{ind2int}, this completes the proof
of the lemma.
\end{proof}

We can now proceed with
\begin{lem}
The image of the embedding of $A_{\mathbb{E},\rho,\mathbb{F}Q}^{t^{2}}$ into
$\widetilde{A}_{\mathbb{E}}^{^{(\mathbb{E}^{\times})^{2}}}$ by
$g\mapsto\big(g,t(g)\big)$ preserves the $\mathbb{F}$-subspace
$(A_{\mathbb{E}}^{-})_{\rho,Q}$ in the action of
$\widetilde{A}_{\mathbb{E}}^{^{(\mathbb{E}^{\times})^{2}}}$ on
$A_{\mathbb{E}}^{-}$. Adding the element $(Q,|Q|^{2})\tilde{\theta}$ from Lemma
\ref{Qtheta} gives a group, which contains
$A_{\mathbb{E},\rho,\mathbb{F}Q}^{t^{2}}$ as a subgroup of index 2, and the
larger group also maps to $O\big((A_{\mathbb{E}}^{-})_{\rho,Q}\big)$.
\label{ac6gen}
\end{lem}

\begin{proof}
An element $u \in (A_{\mathbb{E}}^{-})_{\rho,Q}$ satisfies the condition of
Lemma \ref{sp6gen}, and for $g \in A_{\mathbb{E},\rho,\mathbb{F}Q}^{\times}$ we
have the formulae for $g^{\rho}$ and $\overline{g}^{\rho}$ from the proof of
Lemma \ref{Qtheta}. We evaluate \[(gu\overline{g})^{\rho}=t(g)Q\overline{g}^{\
-1}Q^{-1}\cdot\frac{-Q\tilde{u}Q}{|Q|^{2}}t(g)Q^{-1}g^{-1}Q=-\frac{QN^{A}_{
\mathbb{F}}(g)\overline{g}^{\ -1}\tilde{u}g^{-1}Q}{|Q|^{2}}\] (recall the
$A_{\mathbb{E},\rho,\mathbb{F}Q}^{t^{2}}$ condition), and this equals
$\frac{Q\widetilde{gu\overline{g}}Q}{|Q|^{2}}$ by Lemma \ref{AxA-rel} (over
$\mathbb{E}$). This proves the first assertion. The remaining assertions follow
directly from Lemma \ref{Qtheta} and Lemma \ref{ac6d1} over $\mathbb{E}$, as the
quadratic structure on $(A_{\mathbb{E}}^{-})_{\rho,Q}$ is defined as a subset of
$A_{\mathbb{E}}^{-}$. This proves the lemma.
\end{proof}

The issue with the reflections is a bit tricky here:
\begin{lem}
Let $h\in\mathbb{E}^{0}$ be some non-zero element (so that $d=h^{2}$). For any
$g\in(A_{\mathbb{E}}^{-})_{\rho,Q} \cap A_{\mathbb{E}}^{\times}$, the element
$ghQ^{-1}$ of $A_{\mathbb{E}}^{\times}$ lies in
$A_{\mathbb{E},\rho,\mathbb{F}Q}^{t^{2}}$. The combination of this element with
the element from Lemma \ref{Qtheta} operates on $(A_{\mathbb{E}}^{-})_{\rho,Q}$
as the reflection in $g$. \label{ref6gen}
\end{lem}

\begin{proof}
The fact that $Q \in A_{\mathbb{E}}^{-}$ satisfies $Q^{\rho}=Q$,
$g\in(A_{\mathbb{E}}^{-})_{\rho,Q}$, and $h^{\rho}=-h$ allows us to write
\[(ghQ^{-1})Q\overline{(ghQ^{-1})}^{\rho}=-h^{2}gQ^{-1}g^{\rho}=d\frac{g\tilde{g
}Q}{|Q|^{2}}=\frac{d|g|^{2}}{|Q|^{2}}Q.\] As
$N^{A}_{\mathbb{F}}(ghQ^{-1})=\frac{h^{4}N^{A}_{\mathbb{F}}(g)}{N^{A}_{\mathbb{F
}}(Q)}$ equals $\big(\frac{d|g|^{2}}{|Q|^{2}}\big)^{2}$ by Corollary
\ref{NAvn2}, $ghQ^{-1}$ indeed lies in $A_{\mathbb{E},\rho,\mathbb{F}Q}^{t^{2}}$
with $t(ghQ^{-1})=\frac{d|g|^{2}}{|Q|^{2}}$. We may decompose the resulting
element of $\widetilde{A}_{\mathbb{E}}^{^{(\mathbb{E}^{\times})^{2}}}$ as the
product of $(g,|g|^{2})$, $(h,d)=(h,h^{2})$ (which acts trivially on
$A_{\mathbb{E}}^{-}$), and $\big(Q^{-1},\frac{1}{|Q|^{2}}\big)$. In the product
with the element from Lemma \ref{Qtheta} the two terms involving $Q$ cancel,
and the total action (on $A_{\mathbb{E}}^{-}$) is as
$u\mapsto\frac{g\tilde{u}\overline{g}}{|g|^{2}}$, which was seen in Lemma
\ref{ref6d1} to be the reflection in $g$. As $(A_{\mathbb{E}}^{-})_{\rho,Q}$
inherits its quadratic structure from $A_{\mathbb{E}}^{-}$ and $g$ lies in the
smaller space, this total operation is the reflection in $g$ also on
$(A_{\mathbb{E}}^{-})_{\rho,Q}$. This proves the lemma.
\end{proof}

We are now in position to prove
\begin{thm}
The Gspin group $Gspin\big((A_{\mathbb{E}}^{-})_{\rho,Q}\big)$ is
$A_{\mathbb{E},\rho,\mathbb{F}Q}^{t^{2}}$. The spin group
$spin\big((A_{\mathbb{E}}^{-})_{\rho,Q}\big)$ is the subgroup
$A_{\mathbb{E},\rho,Q}^{1}$ consisting of those elements $g \in
A_{\mathbb{E}}^{1}$ satisfying $gQ\overline{g}^{\rho}=Q$. \label{dim6gen}
\end{thm}

\begin{proof}
Theorem \ref{dim6d1} shows that the semi-direct product of
$\{1,\tilde{\theta}\}$ with
$\widetilde{A}_{\mathbb{E}}^{(\mathbb{E}^{\times})^{2}}$, which is also
generated by $(Q,|Q|^{2})\tilde{\theta}$ and the latter group, maps surjectively
onto $O(A_{\mathbb{E}}^{-})$, with kernel $\mathbb{E}^{\times}$, and such that
$\widetilde{A}_{\mathbb{E}}^{(\mathbb{E}^{\times})^{2}}$ is the inverse image of
$O(A_{\mathbb{E}}^{-})$. Lemma \ref{ac6gen} implies that the subgroup generated
by $(Q,|Q|^{2})\tilde{\theta}$ and the image of
$A_{\mathbb{E},\rho,\mathbb{F}Q}^{t^{2}}$ under $g\mapsto\big(g,t(g)\big)$ lies
in the inverse image of the subgroup $O\big((A_{\mathbb{E}}^{-})_{\rho,Q}\big)$
of $O(A_{\mathbb{E}}^{-})$ (respecting the $\mathbb{F}$-structure), and Theorem
\ref{dim6d1} and the fact the elements of
$O\big((A_{\mathbb{E}}^{-})_{\rho,Q}\big)$ have the same determinant in
$O\big((A_{\mathbb{E}}^{-})_{\rho,Q}\big)$ and in $O(A_{\mathbb{E}}^{-})$ (this
is just extension of scalars) show that the inverse image of
$SO\big((A_{\mathbb{E}}^{-})_{\rho,Q}\big)$ in the smaller group is
$A_{\mathbb{E},\rho,\mathbb{F}Q}^{t^{2}}$. Lemma \ref{ref6gen} and Proposition
\ref{CDT} imply that the map from the group generated by
$A_{\mathbb{E},\rho,\mathbb{F}Q}^{t^{2}}$ and $(Q,|Q|^{2})\tilde{\theta}$ to
$O\big((A_{\mathbb{E}}^{-})_{\rho,Q}\big)$ is surjective. Hence the map
$A_{\mathbb{E},\rho,\mathbb{F}Q}^{t^{2}} \to
SO\big((A_{\mathbb{E}}^{-})_{\rho,Q}\big)$ also surjects, and its kernel,
which consists of those scalars $r\in\mathbb{E}^{\times}$ such that
$rQ\overline{r}^{\rho}=r^{2}Q$, is precisely $\mathbb{F}^{\times}$. Note that
elements $r\in\mathbb{E}_{0}$ also lie in
$A_{\mathbb{E},\rho,\mathbb{F}Q}^{t^{2}}$, but as they satisfy $t(r)=-r^{2}$
they not in the kernel (they operate as $-Id_{(A_{\mathbb{E}}^{-})_{\rho,Q}}$).
It follows that $Gspin\big((A_{\mathbb{E}}^{-})_{\rho,Q}\big)$ is indeed
$A_{\mathbb{E},\rho,\mathbb{F}Q}^{t^{2}}$. Now, the proof of Theorem
\ref{dim6d1} shows that the spinor norm map factors through the projection
$(g,t) \mapsto t$ before going over to any quotient
($\mathbb{E}^{\times}/(\mathbb{E}^{\times})^{2}$, or anything else). Thus
elements of $SO^{1}\big((A_{\mathbb{E}}^{-})_{\rho,Q}\big)$ are obtained from
$g \in A_{\mathbb{E},\rho,\mathbb{F}Q}^{t^{2}}$ with
$t(g)\in(\mathbb{F}^{\times})^{2}$, and the usual interplay with scalars
reduces to those $g$ with $t(g)=1$. But such $g$ preserve $Q$ in the twisted
operation, and must lie in $A_{\mathbb{E}}^{1}$ by the
$A_{\mathbb{E},\rho,\mathbb{F}Q}^{t^{2}}$ condition, which proves the second
assertion since the scalars which square to 1 in $\mathbb{F}$ (which again form
the kernel) are just $\pm1$. This completes the proof of the theorem.
\end{proof}

The case $d=-1$ of Theorem \ref{dim6gen} gives back Theorem \ref{dim6d1} in the
following way. We have $\mathbb{E}=\mathbb{F}\times\mathbb{F}$, so that
$A_{\mathbb{E}}=A \times A$, and $Q$ lies in the original space (no
re-normalization). As
$A_{\mathbb{E}}^{(\mathbb{E}^{\times})^{2}}=A^{(\mathbb{F}^{\times})^{2}} \times
A^{(\mathbb{F}^{\times})^{2}}$, taking $g \in A^{(\mathbb{F}^{\times})^{2}}$ to
be the first coordinate and choosing $t$ such that
$t^{2}=N^{A}_{\mathbb{F}}(g)$, we find that the condition
$(g,h)Q\overline{(g,h)}^{\rho}=tQ$ is equivalent to each one of the conditions
$gQ\overline{h}=tQ$ and $hQ\overline{g}=tQ$, hence to $h$ being
$tQ\overline{g}^{\ -1}Q^{-1}$ (indeed, with the same reduced norm as $g$). Hence
$A_{\mathbb{E},\rho,\mathbb{F}Q}^{t^{2}} \cong
\widetilde{A}^{(\mathbb{F}^{\times})^{2}}$, through the map $(g,t) \mapsto
(g,tQ\overline{g}^{\ -1}Q^{-1})$ (in the opposite direction), so that the second
coordinate associated with $(g,t)$ depends on $Q$, but the isomorphism class of
the group does not (see also Proposition \ref{NQF2NAE} below for the general
case). A similar assertion holds for $A_{\mathbb{E},\rho,Q}^{1}$.

\smallskip

For the independence property, we observe that apart from the choice of the
element $Q$, there is also the choice of the $\mathbb{F}$-structure on
$A_{\mathbb{E}}$, i.e., of the interpretation of $x^{\rho}$ for $x \in
A_{\mathbb{E}}$. We thus replace this notation by defining $\sigma$ to be a ring
automorphism of order 2 whose restriction to $\mathbb{E}$ is $\rho$, and which
satisfies $(\overline{x})^{\sigma}=\overline{x^{\sigma}}$ for every $x \in A$.
By replacing $\sigma$ with another such automorphism, say $\tau$, we get a
presentation of $A_{\mathbb{F}}$ as the tensor product of another bi-quaternion
algebra over $\mathbb{F}$ with $\mathbb{E}$. Any element $y \in
A_{\mathbb{E}}^{\times}$ such that $\overline{y}y\in\mathbb{E}^{\times}$
produces such an automorphism by defining $x^{\tau}=yx^{\sigma}y^{-1}$, where
the similitude condition is equivalent to the condition
$(\overline{x})^{\tau}=\overline{x^{\tau}}$ for every $x \in A$: Compare
$(\overline{x})^{\tau}=y\overline{x}^{\sigma}y^{-1}$ with
$\overline{x^{\tau}}=\overline{y}^{-1}\overline{x}^{\sigma}\overline{y}$ and use
the centrality of $A_{\mathbb{E}}$. The relations we have are contained in the
following
\begin{prop}
$(i)$ Let $Q$ and $R$ be $\sigma$-invariant elements of $A_{\mathbb{E}}^{-}$
(i.e., elements of $A^{-}$) such that
$|R|^{2}\in|Q|^{2}(\mathbb{F}^{\times})^{2}N^{A}_{\mathbb{F}}(A^{\times})$. Then
the groups $A_{\mathbb{E},\sigma,\mathbb{F}Q}^{t^{2}}$ and
$A_{\mathbb{E},\sigma,\mathbb{F}R}^{t^{2}}$ as well as
$A_{\mathbb{E},\sigma,Q}^{1}$ and $A_{\mathbb{E},\sigma,R}^{1}$ are isomorphic.
$(ii)$ Assume that $e \in A_{\mathbb{E}}^{\mathbb{F}^{\times}}$ is such that
$ee^{-\sigma}$ is invariant under $x\mapsto\overline{x}^{\sigma}$, and $S \in
A_{\mathbb{E}}^{-}$ is invariant under $\tau:x \mapsto
ee^{-\sigma}x^{\sigma}e^{\sigma}e^{-1}$ and satisfies
$|S|^{2}\in|Q|^{2}N^{A_{\mathbb{E}}}_{\mathbb{E}}(e)(\mathbb{F}^{\times})^{2}N^{
A}_{\mathbb{F}}(A^{\times})$. In this case the groups
$A_{\mathbb{E},\sigma,\mathbb{F}Q}^{t^{2}}$ and
$A_{\mathbb{E},\tau,\mathbb{F}S}^{t^{2}}$ and the groups
$A_{\mathbb{E},\sigma,Q}^{1}$ and $A_{\mathbb{E},\tau,R}^{1}$ are isomorphic.
$(iii)$ Let $b \in A_{\mathbb{E}}^{\times}$ be such that
$\overline{b}^{\sigma}=b$ and $b\overline{b}\in\mathbb{E}^{\times}$, and define
$\eta:x \mapsto bx^{\sigma}b^{-1}$. If the element $T=Qb^{-1}$ lies in
$A_{\mathbb{E}}^{-}$ then the groups $A_{\mathbb{E},\sigma,\mathbb{F}Q}^{t^{2}}$
and $A_{\mathbb{E},\eta,\mathbb{F}T}^{t^{2}}$ coincide, and same assertion
holds for $A_{\mathbb{E},\sigma,Q}^{1}$ and $A_{\mathbb{E},\eta,T}^{1}$.
\label{NQF2NAE}
\end{prop}

\begin{proof}
Considering $A^{-}$ as a space of discriminant 1 again, we find (as in the
proof of Proposition \ref{NQF2NA}) that we can write
$R=rcQ\overline{c}=rcQ\overline{c}^{\sigma}$ with $r\in\mathbb{F}^{\times}$ and
$c \in A^{\times}$ (hence $c=c^{\sigma}$). It follows that conjugation by $c$
takes the group $A_{\mathbb{E},\sigma,\mathbb{F}Q}^{t^{2}}$ to
$A_{\mathbb{E},\sigma,\mathbb{F}R}^{t^{2}}$ (as conjugation preserves reduced
norms, the $N^{A}_{\mathbb{E}}=t^{2}$ condition is also preserved), in a manner
which commutes with the multiplier maps. This proves part $(i)$, since the
groups $A_{\mathbb{E},\sigma,Q}^{1}$ and $A_{\mathbb{E},\sigma,R}^{1}$ are
defined by the condition $t=1$ on the larger groups
$A_{\mathbb{E},\sigma,\mathbb{F}Q}^{t^{2}}$ to
$A_{\mathbb{E},\sigma,\mathbb{F}R}^{t^{2}}$. For part $(ii)$, first note that
the $\sigma$-image of $ee^{-\sigma}$ is its inverse, so that the condition on
that element implies that it is a similitude and
$(\overline{x})^{\tau}=\overline{x^{\tau}}$ holds for every $x \in
A_{\mathbb{E}}$. Consider the element $eQ\overline{e}$ of $A_{\mathbb{E}}^{-}$,
whose vector norm is
$N^{A_{\mathbb{E}}}_{\mathbb{E}}(e)|Q|^{2}\in\mathbb{F}^{\times}$ by Proposition
\ref{NAFg} and the assumption on $e$. The relation between $\tau$ and $\sigma$
and the assumption on $Q$ and $e$ show that this element is $\tau$-invariant,
and as the properties of $e$ and $\tau$ imply that for $g \in
A_{\mathbb{E}}^{\times}$ we have
\[\overline{ege^{-1}}^{\tau}=\overline{e}^{-1}\overline{e}^{\sigma}
\cdot\overline{ege^{-1}}^{\sigma}\cdot\overline{e}^{-\sigma}\overline{e}
=\overline{e}^{-1}\overline{g}^{\sigma}\overline{e},\] it follows that
conjugation by $e$ sends $A_{\mathbb{E},\sigma,\mathbb{F}Q}^{t^{2}}$ to
$A_{\mathbb{E},\tau,\mathbb{F}eQ\overline{e}}^{t^{2}}$ as well as
$A_{\mathbb{E},\sigma,Q}^{1}$ to $A_{\mathbb{E},\tau,eQ\overline{e}}^{1}$. As
$S$ is related to $eQ\overline{e}$ in the same manner as $Q$ and $R$, part
$(ii)$ now follows from part $(i)$, using conjugation by an appropriate
$\tau$-invariant element $d \in A_{\mathbb{E}}^{\times}$. Now, the fact that
$Q \in A^{-}$ and the element $T$ from part $(iii)$ lies in $A_{\mathbb{E}}^{-}$
imply
$T^{\eta}=-\overline{T}^{\eta}=-b\overline{Qb^{-1}}^{\sigma}b^{-1}=+b\overline{b
}^{-\sigma}Qb^{-1}$, and this gives $Qb^{-1}=T$ again by our assumption on $b$.
Furthermore, for $g \in A$ we have
$gT\overline{g}^{\eta}=gQ\overline{g}^{\sigma}b^{-1}$ by the definitions of
$T$ and $\eta$, so that the $A_{\mathbb{E},\sigma,\mathbb{F}Q}^{t^{2}}$ and
$A_{\mathbb{E},\eta,\mathbb{F}S}^{t^{2}}$ conditions coincide, with the same
multiplier $t$. This proves part $(iii)$ (for both the Gspin and spin groups),
and completes the proof of the proposition.
\end{proof}

Note that part $(i)$ of Proposition \ref{NQF2NAE} is a special case of part
$(ii)$ there, but it is required for the proof of the latter part. We also
remark that for $\tau$ and $\sigma$ which are connected as in part $(ii)$ of
Proposition \ref{NQF2NAE}, the subring which $\tau$ stabilizes is isomorphic to
$A$. Indeed, if $x=x^{\tau}=aa^{-\sigma}x^{\sigma}a^{\sigma}a^{-1}$ then
$a^{-1}xa$ is $\sigma$-invariant, so that conjugation by $a$ maps $A$ onto this
subring (in fact, the Skolem--Noether Theorem implies that every automorphism on
$A_{\mathbb{E}}$ which stabilizes a subring which is isomorphic to $A$ must
arise in this way). This is related to the fact that the vector norm relation is
based on $N^{A}_{\mathbb{F}}(A^{\times})$ also when $\tau$ is involved. On the
other hand, $\eta$ from part $(iii)$ of that proposition may arise from a
non-isomorphic bi-quaternion algebra over $\mathbb{F}$. E.g., assuming that $Q
\in C_{0}$, any element $b\in\mathbb{F}\oplus\mathbb{E}_{0}C_{0}$ such that
$Tr^{C_{\mathbb{E}}}_{\mathbb{E}}(Qb^{-1})=0$ satisfies all the assumptions
(including $Qb^{-1} \in A_{\mathbb{E}}^{-}$), and we have seen in the paragraph
preceding Corollary \ref{iso4} that up to letting $Q \in C_{0}$ vary, this
covers all the bi-quaternion algebras which arise as the tensor product of $B$
with some quaternion algebra over $\mathbb{F}$ which becomes $C_{\mathbb{E}}$
over $\mathbb{E}$. We may also apply this construction with $b \in B_{0} \otimes
Q$, yielding an operation on both $B$ in $C$. We remark that by Proposition 2.18
of \cite{[KMRT]}, any two involutions $x\mapsto\overline{x}^{\sigma}$ and
$x\mapsto\overline{x}^{\tau}$ must be related through conjugation by some
element $b$ satisfying $\overline{b}^{\sigma}=b$, with the similitude condition
to preserve commutativity with $\iota_{B}\otimes\iota_{C} \otimes
Id_{\mathbb{E}}$. It seems likely that the similitude condition implies the
existence of $Q \in A^{-} \cap A^{\times}$ such that $Qb^{-1} \in
A_{\mathbb{E}}^{-}$, but we have not checked this in detail. Was this the case,
Proposition \ref{NQF2NAE} would relate all the involutions which commute with
$\iota_{B}\otimes\iota_{C} \otimes Id_{\mathbb{E}}$ to one another, yielding an
$\mathbb{F}$-structure invariance result.

\smallskip

Fixing $Q$ and $\sigma$ back again (and writing $\rho$ for the automorphism of
$A_{\mathbb{E}}$ as before), the assertion involving isotropy in this case is
\begin{cor}
If the space becomes isotropic when one extends scalars to $\mathbb{E}$ then
the Gspin group is a ``quaternionic $GSU$ group'', consisting of the
matrices $g \in M_{2}(B_{\mathbb{E}})$, for some quaternion algebra
$B_{\mathbb{E}}$ over $\mathbb{E}$ which comes from a quaternion algebra $B$
over $\mathbb{F}$, which multiply a diagonal matrix of the sort $\binom{-\delta\
\ 0}{\ \ 0\ \ 1} \in M_{2}(\mathbb{F})$ with $\delta$ determined up to
multiplication from $N^{B}_{\mathbb{F}}(B^{\times})$, through the action $g:X
\mapsto gX\overline{g}^{t}$. The $GSU$ condition means that the reduced norm of
these matrices equals the square of the multiplier of this element. The spin
group is the associated ``quaternionic special unitary group'', of matrices
which preserve this element and have reduced norm 1. In case the space splits
two hyperbolic planes over $\mathbb{E}$, the Gspin group becomes the
$GSU_{\mathbb{E},\rho}$ group of a unitary space of dimension 4 and determinant
(or discriminant) 1 over $\mathbb{E}$ (with $\rho$), and the spin group is the
associated special unitary group. \label{iso6gen}
\end{cor}

\begin{proof}
In the first case $A_{\mathbb{E}}$ is isomorphic to $M_{2}(B_{\mathbb{E}})$ for
some quaternion algebra $B$ over $\mathbb{F}$. We may normalize (using parts
$(ii)$ and $(iii)$ of Proposition \ref{NQF2NAE} if necessary) the involutions
such that $A=M_{2}(B)$, i.e., with $C=M_{2}(\mathbb{F})$, and as every element
of $\mathbb{F}^{\times}$ is a norm from $C_{0}$, part $(i)$ of Proposition
\ref{NQF2NAE} allows us to restrict attention to spaces
$\big(M_{2}(B_{\mathbb{E}})^{-}\big)_{\rho,Q}$ with $Q \in
M_{2}(\mathbb{F})_{0}$. Moreover, as in the proof of Corollary \ref{iso5}, we
may take $Q$ of the form $\binom{0\ \ \delta}{1\ \ 0}$, and part $(i)$ of
Proposition \ref{NQF2NAE} and Lemma \ref{NM2B} allows us to take $\delta$ from a
set of representatives for $\mathbb{F}^{\times}/N^{B}_{\mathbb{F}}(B^{\times})$.
By Lemma \ref{Sadjt}, the condition $g \in
A_{\mathbb{E},\sigma,\mathbb{F}Q}^{t^{2}}$ (or $g \in
GL_{2}^{t^{2}}(B_{\mathbb{E}})_{\sigma,\mathbb{F}Q}$) becomes $g\binom{-\delta\
\ 0}{\ \ 0\ \ 1}\iota_{B}(g)^{t\rho}=t\binom{-\delta\ \ 0}{\ \ 0\ \ 1}$ for some
$t=t(g)\in\mathbb{F}^{\times}$ and
$N^{M_{2}(B_{\mathbb{E}})}_{\mathbb{E}}(g)=t(g)^{2}$, while elements of
$A_{\mathbb{E},\sigma,Q}^{1}=GL_{2}^{1}(B_{\mathbb{E}})_{\sigma,Q}$ have reduced
norm 1 and their action leaves the element $\binom{-\delta\ \ 0}{\ \ 0\ \ 1}$
invariant. This proves the first assertion.

If after tensoring $(A_{\mathbb{E}}^{-})_{\rho,Q}$ with $\mathbb{E}$ we get a
2-dimensional isotropic subspace (and then we even get a 3-dimensional such
space) then $B_{\mathbb{E}}$ also splits, hence $B$ may be presented as
$\big(\frac{-d,\varepsilon}{\mathbb{F}}\big)$ for some
$\varepsilon\in\mathbb{F}^{\times}$ (which is unique up to
$(\mathbb{F}^{\times})^{2}$ and multiplication by $+d$). The model we thus take
for $B$ is the algebra $(\mathbb{E},\rho,\varepsilon)$, and by Lemma
\ref{KsplitB}, $M \mapsto M^{\rho}$ is $\binom{a\ \ b}{c\ \ d}\mapsto\binom{0\ \
\varepsilon}{1\ \ 0}\binom{a^{\rho}\ \ b^{\rho}}{c^{\rho}\ \ d^{\rho}}\binom{0\
\ \varepsilon}{1\ \ 0}/\varepsilon$ on $B_{\mathbb{E}}=M_{2}(\mathbb{E})$. Thus
the operation $g\mapsto\overline{g}^{\rho}$ or $g\mapsto\iota_{B}(g)^{t\rho}$
involves the map $M\mapsto\binom{0\ \
\varepsilon}{1\ \ 0}M^{\rho}\binom{0\ \ \varepsilon}{1\ \
0}/\varepsilon$ on the entries (in addition to adjoint or transposition of
$2\times2$ matrices over $B_{\mathbb{E}}$). It follows that each entry from
$B_{\mathbb{E}}$ in $\binom{0\ \ \delta}{1\ \ 0}$ or $\binom{-\delta\ \ 0}{\ \
0\ \ 1}$ appearing in the definition of the groups (which we write as $\binom{0\
\ \delta I}{I\ \ \ 0\ }$ or $\binom{-\delta I\ \ 0}{\ \ 0\ \ \ I}$ since
$B_{\mathbb{E}}$ is the matrix algebra $M_{2}(\mathbb{E})$) is replaced by the
corresponding multiple of $\binom{0\ \ \varepsilon}{1\ \ 0}$ when we apply
$\rho$ on the $\mathbb{E}$-entries of matrices in
$M_{2}(B_{\mathbb{E}})=M_{4}(\mathbb{E})$. Applying Lemma \ref{Sadjt} to the
$2\times2$-entries changes each such $\binom{0\ \ \varepsilon}{1\ \ 0}$ to
$\binom{-\varepsilon\ \ 0}{\ \ 0\ \ 1}$. Therefore, the
$A_{\mathbb{E},\sigma,\mathbb{F}Q}^{t^{2}}=GL_{2}^{t^{2}}(B_{\mathbb{E}})_{
\sigma,\mathbb{F}Q}$ condition just means being in the
$GSU_{\mathbb{E},\rho}$ group of a 4-dimensional space over $\mathbb{E}$ with a
sesqui-linear form having an orthogonal basis with norms $\delta\varepsilon$,
$-\varepsilon$, $-\delta$, and 1 (hence the determinant is 1 modulo
$N^{\mathbb{E}}_{\mathbb{F}}(\mathbb{E}^{\times})$), while the group
$A_{\mathbb{E},\sigma,Q}^{1}=GL_{2}^{1}(B_{\mathbb{E}})_{\sigma,Q}$ is the 
corresponding $SU_{\mathbb{E},\rho}$ group. This proves the corollary.
\end{proof}

We emphasize that only unitary spaces with determinant (or discriminant) 1
appear in this context.

In this isotropic case we have the following equivalent representations:
\begin{cor}
The groups $GL_{2}^{t^{2}}(B_{\mathbb{E}})_{\sigma,\mathbb{F}Q}$ and
$GL_{2}^{1}(B_{\mathbb{E}})_{\sigma,Q}$ with $Q=\binom{0\ \ \delta}{1\ \ 0}$ are
the Gspin and spin groups of the space $\mathbb{E} \oplus B$, embedded as
$(y,\lambda)\mapsto\binom{\lambda\ \ \ -\delta y}{y^{\rho}\ \
-\overline{\lambda}\ }$ or as $(y,\lambda)\mapsto\binom{\delta y\ \ \lambda\
}{\overline{\lambda}\ \ \ y^{\rho}}$, with the vector norm being
$N^{B}_{\mathbb{F}}(\lambda)-\delta N^{\mathbb{E}}_{\mathbb{F}}(y)$ and the
actions being $g:M \mapsto gM\overline{g}$ and $g:X \mapsto gX\iota{B}(g)^{t}$
respectively. The $GSU_{\mathbb{E},\rho}$ and $SU_{\mathbb{E},\rho}$ groups of a
4-dimensional unitary space of discriminant 1 are the Gspin and spin groups of
the direct sum of 3 copies of $\mathbb{E}$, with the vector norm of $(z,w,y)$
being $N^{\mathbb{E}}_{\mathbb{F}}(z)-\varepsilon
N^{\mathbb{E}}_{\mathbb{F}}(w)-\delta N^{\mathbb{E}}_{\mathbb{F}}(y)$. It may be
embedded in $M_{4}(\mathbb{E})$ by replacing $\lambda$ by $\binom{\ z\ \
\varepsilon w}{w^{\rho}\ \ z^{\rho}}$ with $z$ and $w$ from $\mathbb{E}$ in the
two representations from above. In addition, these groups act via the operation
$g:T \mapsto gTg^{t}$ on the subspace of $M_{4}^{as}(\mathbb{E})$ consisting of
those matrices $\binom{a\ \ b}{c\ \ d}$ in which the anti-symmetric matrices $a$
and $d=\binom{0\ \ -y}{y\ \ \ \ 0}$ satisfies $a=\delta d^{\rho}$ while
$b=-c^{t}$ takes the form $\binom{\varepsilon w\ \ -z\ }{z^{\rho}\ \
-w^{\rho}}$, as well as through $\binom{a\ \ b}{c\ \ d}:S\mapsto\binom{a\ \
b}{c\ \ d}S\binom{\ \ d^{t}\ \ -b^{t}}{-c^{t}\ \ \ \ a^{t}}$ on another
embedding of the three copies of $\mathbb{E}$ into the complement of
$\mathfrak{sp}_{4}(\mathbb{E})$ in $M_{4}(\mathbb{E})$. \label{alt6gen}
\end{cor}

\begin{proof}
Recall that $(A_{\mathbb{E}}^{-})_{\rho,Q}$ is $(Q^{\perp} \subseteq
A^{-})\oplus\mathbb{E}_{0}Q$, and the off-diagonal matrices which lie in
$Q^{\perp}$ are spanned by $\binom{0\ \ -\delta}{1\ \ \ \ 0}$. Adding
$\mathbb{E}_{0}Q$ yields the first representation, on which the action is via
$g:M\mapsto\frac{gM\overline{g}}{t(g)}$, and the second one, on which the groups
operated via $g:X\mapsto\frac{gX\iota_{B}(g)^{t}}{t(g)}$, arises from Lemma
\ref{Sadjt} as in the proofs of Corollaries \ref{alt6d1} and \ref{alt5} (with
the norm being some generalization of the Moore determinant). The two
representations of the $GSU_{\mathbb{E},\rho}$ and $SU_{\mathbb{E},\rho}$ groups
are also obtained by applying Lemma \ref{Sadjt} to the entries as $2\times2$
matrices, as the proofs of Corollaries \ref{alt6d1} and \ref{alt5} do, while
recalling that $B$ is embedded into $M_{2}(\mathbb{E})$ as
$(\mathbb{E},\rho,\varepsilon)$. This proves the corollary. 
\end{proof}

We remark that for the case where $A=M_{4}(\mathbb{F})$ we have used a different
choice of $Q$ in Corollaries \ref{iso5} and \ref{alt5} in order to obtain the
classical symplectic group. Choosing this $Q$ for Corollaries \ref{iso6gen} and
\ref{alt6gen} yields the $GSU_{\mathbb{E},\rho}$ and $SU_{\mathbb{E},\rho}$
conditions for an anti-diagonal symmetric matrix (which is explicitly $\binom{0\
\ E}{E\ \ 0}$ with $E=\binom{0\ \ 1}{1\ \ 0}$). On the other hand, here the
splitting of $A_{\mathbb{E}}$ might come from $A=M_{2}(B)$ with $B$ which is
non-split over $\mathbb{F}$, so that we have many groups in this case. In cany
case, we restrict attention to the classical unitary groups (of diagonal
matrices), as any unitary group is conjugate to a classical one.

\section{Relations with the Exterior Square \label{Wedge}}

For the groups arising from bi-quaternion algebras, hence related to $4\times4$
matrices, there are representations which are equivalent to those presented
here. Given a field $\mathbb{M}$, the group $GL_{4}(\mathbb{M})$ operates on
the 6-dimensional exterior square $\bigwedge^{2}\mathbb{M}^{4}$ of the natural
representation space $\mathbb{M}^{4}$, and if we denote the canonical
 basis for $\mathbb{M}^{4}$ by $e_{i}$, $1 \leq i\leq4$ for $\mathbb{M}^{4}$
then the 6 elements $e_{i} \wedge e_{j}$ with $1 \leq i<j\leq4$ form a basis for
$\bigwedge^{2}\mathbb{M}^{4}$. The map taking $u$ and $v$ from
$\bigwedge^{2}\mathbb{M}^{4}$ to the multiple of $e_{1} \wedge e_{2} \wedge
e_{3} \wedge e_{4}$ which equals $u \wedge w \in \bigwedge^{4}\mathbb{M}^{4}$ is
bilinear and symmetric (hence we denote this value by $\langle u,w \rangle$),
and in fact independent, up to rescaling, of the choice of basis. The connection
which we have arises from the following
\begin{lem}
For any $n$, the natural representation of $GL_{n}(\mathbb{M})$ on
$\bigwedge^{2}\mathbb{M}^{n}$ is equivalent to its representation on the space
$M_{n}^{as}(\mathbb{M})$ of anti-symmetric matrices via $g:T \mapsto gTg^{t}$.
For $n=4$ the equivalence preserves the bilinear forms, where on
$\bigwedge^{2}\mathbb{M}^{4}$ we take the bilinear form defined above and on
$M_{n}^{as}(\mathbb{M})$ we use the Pfaffian. \label{wedge2equiv}
\end{lem}

\begin{proof}
Using the standard basis for $\mathbb{M}^{n}$, consider the map which takes, for
all $1 \leq i<j \leq n$, the basis element $e_{i} \wedge e_{j}$ of
$\bigwedge^{2}\mathbb{M}^{n}$ to the anti-symmetric matrix $E_{ij}-E_{ji}$. The
fact that for $g=(a_{kl})$ we have
$\big(g(E_{ij}-E_{ji})g^{t}\big)_{kl}=a_{ki}a_{lj}-a_{kj}a_{li}$ shows that
this is indeed an equivalence of representations. For $n=4$ we find that $e_{1}
\wedge e_{2}$ and $e_{3} \wedge e_{4}$ as well as $e_{1} \wedge e_{4}$ and
$e_{2} \wedge e_{3}$ span hyperbolic planes, while $e_{1} \wedge e_{3}$ and
$e_{2} \wedge e_{4}$ span another hyperbolic plane but with the sign inverted.
As this is in correspondence with the values bilinear form arising from the
Pfaffian takes on the images of these vectors in $M_{n}^{as}(\mathbb{M})$, this
completes the proof of the lemma.
\end{proof}

Given $g \in GL_{4}(\mathbb{M})$ and $u$ and $v$ from
$\bigwedge^{2}\mathbb{M}^{4}$, we have $\langle gu,gw \rangle=\det g\langle u,w
\rangle$. This is so, since $gu \wedge gv$ is $\langle u,w \rangle$ times the
image of the generator $e_{1} \wedge e_{2} \wedge e_{3} \wedge e_{4}$ of
$\bigwedge^{4}\mathbb{M}^{4}$ via the $\bigwedge^{4}$-action of $g$, and the
latter action is multiplication by the determinant by definition. Lemma
\ref{wedge2equiv} and the isomorphism of representations appearing in the proof
of Corollary \ref{alt6d1} thus provide an alternative proof for Proposition
\ref{NAFg}, via appropriate restriction of scalars.

Assume now that $ch\mathbb{M}\neq2$. Then sums and differences of complementary
pairs form an orthogonal basis for $\mathbb{M}^{4}$. All the representations in
this Section will be variants of this one, using this basis. The first one is
\begin{prop}
$SL_{4}(\mathbb{F})$ is the spin group of the quadratic space
$\bigwedge^{2}\mathbb{M}^{4}$, and the Gspin group is the double cover
$\widetilde{GL}_{4}^{(\mathbb{F}^{\times})^{2}}(\mathbb{F})$, operating with
division by the chosen square root. The classical $Sp_{4}(\mathbb{F})$ and
$GSp_{4}(\mathbb{F})$ are the spin and Gspin groups of the subspace of
$\bigwedge^{2}\mathbb{M}^{4}$ which is spanned by $e_{1} \wedge e_{2}$, $e_{1}
\wedge e_{4}$, $e_{2} \wedge e_{3}$, $e_{3} \wedge e_{4}$, and $e_{1} \wedge
e_{3}-e_{2} \wedge e_{4}$.
\label{wedge2splitF}
\end{prop}

\begin{proof}
The first two assertions follow directly from Corollary \ref{alt6d1} by taking
$\mathbb{M}=\mathbb{F}$ in Lemma \ref{wedge2equiv}. For the last two assertions
we apply Corollary \ref{alt5}, observing that our element $Q$ is taken through
all these maps to $e_{1} \wedge e_{3}+e_{2} \wedge e_{4}$, the orthogonal
complement of which is spanned by the asserted vectors. This proves the
proposition.
\end{proof}

When we consider spaces of general discriminant $d$, let $\mathbb{E}$ be
$\mathbb{F}(\sqrt{d})$ with Galois automorphism $\rho$. Considering the
bi-quaternion algebra $A=M_{2}(B)$ with
$B=\big(\frac{d,\varepsilon}{\mathbb{F}}\big)$ (which splits over
$\mathbb{E}$), and choosing some $\delta\in\mathbb{F}^{\times}$, we find that
\begin{prop}
The $SU_{\mathbb{E},\rho}$ group of a space with an orthogonal basis having
norms $\delta\varepsilon$, $-\varepsilon$, $-\delta$, and 1 is the spin group of
the $\mathbb{F}$-subspace of $\bigwedge^{2}\mathbb{E}^{4}$ which is the direct
sum of $\mathbb{F}(e_{1} \wedge e_{4}-e_{2} \wedge e_{3})$,
$\mathbb{E}_{0}(e_{1} \wedge e_{4}+e_{2} \wedge e_{3})$, $\mathbb{F}(e_{2}
\wedge e_{4}-\varepsilon e_{1} \wedge e_{3})$, $\mathbb{E}_{0}(e_{2} \wedge
e_{4}+\varepsilon e_{1} \wedge e_{3})$, $\mathbb{F}(e_{3} \wedge e_{4}+\delta
e_{1} \wedge e_{2})$, and $\mathbb{E}_{0}(e_{3} \wedge e_{4}-\delta e_{1} \wedge
e_{2})$. The Gspin group is the associated $GSU_{\mathbb{E},\rho}$ group.
\label{wedge2splitE}
\end{prop}

\begin{proof}
Multiplying each entry of the elements of the second representation appearing
in Corollary \ref{iso6gen} by $\binom{\ \ 0\ \ 1}{-1\ \ 0}$ from the right
takes, if $\lambda \in B=\big(\frac{d,\varepsilon}{\mathbb{F}}\big)$, this
incarnation of $\mathbb{E} \oplus B$ to the asserted direct sum (the first two
generate the image of the diagonal matrices in $B$, the second two yield those
of the off-diagonal matrices, and the latter two give the image of
$\mathbb{E}$). The assertion now follows from Corollary \ref{alt6gen} and Lemma
\ref{wedge2equiv} for $\mathbb{M}=\mathbb{E}$. This proves the proposition.
\end{proof}

Note that unlike in Proposition \ref{wedge2splitF}, here the image of $Q$
became the element $e_{3} \wedge e_{4}-\delta e_{1} \wedge e_{2}$, appearing in
the 6th direct summand in Proposition \ref{wedge2splitE}.

\smallskip

We now consider the case where $A^{-}$ is isotropic, but does not necessarily
split more than one hyperbolic plane. Then $A$ is isomorphic to $M_{2}(B)$ for
some quaternion algebra $B$, and by taking $\mathbb{K}$ to be a quadratic
extension of $\mathbb{F}$, whose Galois automorphism we denote $\eta$, which
splits $B$, we may write $B\cong(\mathbb{K},\eta,\varepsilon)$ for some
$\varepsilon\in\mathbb{F}^{\times}$. In this case we have
\begin{prop}
The groups $\widetilde{GL}_{2}^{(\mathbb{F}^{\times})^{2}}(B)$ and
$GL_{2}^{1}(B)$ are the Gspin and spin groups of the $\mathbb{F}$-subspace of
$\bigwedge^{2}\mathbb{K}^{4}$ which is obtained as the direct sum of the four
spaces $\mathbb{F}(e_{1} \wedge e_{4}-e_{2} \wedge e_{3})$,
$\mathbb{K}_{0}(e_{1} \wedge e_{4}+e_{2} \wedge e_{3})$, $\mathbb{F}(e_{2}
\wedge e_{4}-\varepsilon e_{1} \wedge e_{3})$, $\mathbb{K}_{0}(e_{2} \wedge
e_{4}+\varepsilon e_{1} \wedge e_{3})$, and the hyperbolic plane
$\mathbb{F}e_{1} \wedge e_{2}\oplus\mathbb{F}e_{3} \wedge e_{4}$. The
quaternionic symplectic groups $Sp_{B}\binom{-\delta\ \ 0}{\ \ 0\ \ 1}$ and
$GSp_{B}\binom{-\delta\ \ 0}{\ \ 0\ \ 1}$ are the spin and Gspin groups of the
subspace which is the direct sum of the first four 1-dimensional spaces with
$\mathbb{F}(e_{3} \wedge e_{4}+\delta e_{1} \wedge e_{2})$. If $d$ is a
non-trivial discriminant and $\mathbb{E}=\mathbb{F}(\sqrt{d})$ has Galois
automorphism $\rho$ over $\mathbb{F}$, then
$GL_{2}^{t^{2}}(B_{\mathbb{E}})_{\rho,\mathbb{F}Q}$ and
$GL_{2}^{1}(B_{\mathbb{E}})_{\rho,Q}$ are the Gspin and spin groups of the
direct sum of the latter 5-dimensional space and $\mathbb{E}_{0}(e_{3} \wedge
e_{4}-\delta e_{1} \wedge e_{2})$. \label{wedge2iso}
\end{prop}

\begin{proof}
The first space is the image of $M_{2}(B)^{-}$ under the isomorphisms from
Corollary \ref{alt6d1} and Lemma \ref{wedge2equiv} with
$\mathbb{M}=\mathbb{K}$. For the second and third space we note that our vector
$Q$ is the same matrix from Proposition \ref{wedge2splitE}, hence has the same
image. This completes the proof of the proposition as we did for Propositions
\ref{wedge2splitF} and \ref{wedge2splitE}, where for the latter assertion we
need to apply Lemma \ref{wedge2equiv} with $\mathbb{M}=\mathbb{K}\mathbb{E}$.
\end{proof}

We remark that if $B$ is not split but $\mathbb{E}$ splits $B$ then taking
$\mathbb{K}=\mathbb{E}$ in Proposition \ref{wedge2iso} yields Proposition
\ref{wedge2splitE} again.

\smallskip

In the general case, where $A$ might be a division algebra, we write $A=B
\otimes C$ take $\mathbb{L}$ to be a quadratic extension of $\mathbb{F}$ (with
Galois automorphism $\omega$) which splits $C$. As we have seen that our choice
of $Q$ may be taken from $C_{0}$, we assume that
$C\cong(\mathbb{L},\omega,\delta)$ with the same $\delta$. Hence in general we
have
\begin{prop}
A space for which the groups $\widetilde{A}^{(\mathbb{F}^{\times})^{2}}$ and
$A^{1}$ appear as the Gspin and spin groups is the direct sum of the
$\mathbb{F}$-spaces $\mathbb{L}_{0}(e_{1} \wedge e_{4}-e_{2} \wedge e_{3})$,
$\mathbb{K}_{0}(e_{1} \wedge e_{4}+e_{2} \wedge e_{3})$, $\mathbb{F}(e_{2}
\wedge e_{4}-\varepsilon e_{1} \wedge e_{3})$, $\mathbb{K}_{0}(e_{2} \wedge
e_{4}+\varepsilon e_{1} \wedge e_{3})$, $\mathbb{L}_{0}(e_{3} \wedge
e_{4}+\delta e_{1} \wedge e_{2})$, and $\mathbb{F}(e_{3} \wedge e_{4}-\delta
e_{1} \wedge e_{2})$ inside $\bigwedge^{2}(\mathbb{K}\mathbb{L})^{4}$. The
stabilizers $A^{\times}_{\mathbb{F}Q}$ and $A^{\times}_{Q}$ are the Gspin and
spin groups of the direct sum of the first 5 spaces, if $Q \in A^{-}$ is chosen
to be $\binom{0\ \ \delta}{1\ \ 0} \in C_{0}$. With this choice of $Q$, a space
of discriminant $d$ for which $A_{\mathbb{E},\rho,\mathbb{F}Q}^{t^{2}}$ is the
Gspin group and $A_{\mathbb{E},\rho,Q}^{1}$ is the spin group, where
$\mathbb{E}=\mathbb{F}(\sqrt{d})$ and $\rho$ is its Galois automorphism, is
obtained by adding $\mathbb{E}_{0}(e_{3} \wedge e_{4}-\delta e_{1} \wedge
e_{2})$ to the latter 5-dimensional space. \label{wedge2gen}
\end{prop}

\begin{proof}
The embedding of $C$ inside $M_{2}(\mathbb{L})$ embeds $A$ into
$M_{2}(B_{\mathbb{L}})$, hence $A^{-}$ into $M_{2}(B_{\mathbb{L}})^{-}$. As the
two scalar entries represent scalar entries of $C \subseteq M_{2}(\mathbb{L})$,
they are related to one another via $\omega$ and multiplication by $\delta$,
yielding the two latter spaces under the maps from Corollary \ref{alt6d1} and
Lemma \ref{wedge2equiv} with $\mathbb{M}=\mathbb{K}\mathbb{L}$. As the
remaining entries come from $\lambda \in B_{\mathbb{L}}$ which in fact lies in
$B_{0}\oplus\mathbb{L}_{0}$, the off-diagonal entries of $\lambda \in
M_{2}(\mathbb{K}\mathbb{L})$ are related through $\eta$ and multiplication by
$\varepsilon$, yielding the middle two spaces via these identifications. The
diagonal entries must therefore be negated by $\eta\omega$, hence come from
$\mathbb{K}_{0}\oplus\mathbb{L}_{0}$, which gives the first two subspaces. This
proves the assertions about $\widetilde{A}^{(\mathbb{F}^{\times})^{2}}$ and
$A^{1}$, and the ones about the symplectic groups follow since we use the same
element $Q=\binom{0\ \ \delta}{1\ \ 0}$ as in Propositions \ref{wedge2splitE}
and \ref{wedge2iso} (note that this element belongs to $C$ as a subalgebra of
$M_{2}(\mathbb{L})$ in our normalizations). For the remaining assertions we use
again the same element $Q$, and we argue as in Propositions \ref{wedge2splitE}
and \ref{wedge2iso}, where $\mathbb{M}$ is now taken to be
$\mathbb{K}\mathbb{L}\mathbb{E}$ in Lemma \ref{wedge2equiv}. This proves the
proposition.
\end{proof}

We remark again that out choice of $Q$ in Proposition \ref{wedge2gen}, which
looks rather special, is entirely general when we normalize $B$ and $C$
appropriately. Note that taking $\mathbb{K}=\mathbb{L}$ in Proposition
\ref{wedge2gen} in case $A$ is not a division algebra, or
$\mathbb{L}=\mathbb{E}$ in case $A$ is a division algebra but $A_{\mathbb{E}}$
is not, does not give us Proposition \ref{wedge2iso} again. This is so, since
the split algebra $C$ is normalized as $M_{2}(\mathbb{F})$ in Proposition
\ref{wedge2iso}, but as some subalgebra $(\mathbb{L},\omega,\delta)$ of
$M_{2}(\mathbb{L})$, with the $\mathbb{F}$-structure from Lemma \ref{KsplitB},
in Proposition \ref{wedge2gen}.

\section{Dimension 8, Isotropic, Discriminant 1 \label{Dim8id1}}

The spaces we consider here are described in the following
\begin{lem}
The direct sum of a hyperbolic plane with an Albert form arising from a
presentation of a bi-quaternion algebra $A$ as the tensor product $B \otimes C$
of two quaternion algebras over $\mathbb{F}$ is isotropic of dimension 8 and
discriminant 1. Any isotropic 8-dimensional space of discriminant 1 is obtained,
up to rescaling and isometries, in this way. \label{sp8id1}
\end{lem}

\begin{proof}
Recall that a space is isotropic if and only if it has a hyperbolic plane as a
direct summand, and that a hyperbolic plane is isometric to its rescalings. The
lemma now follows directly from Lemma \ref{sp6d1}.
\end{proof}
We thus fix a bi-quaternion algebra $A=B \otimes C$ with the (orthogonal)
involution $\iota_{B}\otimes\iota_{C}:x\mapsto\overline{x}$ as above, and the
space from Lemma \ref{sp8id1} may be denoted $A^{-} \oplus H$ (where $H$ stands
for a hyperbolic plane). It will be useful to embed this space into $M_{2}(A)$
by taking the sum of $u \in A^{-}$, $-p$ times one generator of the hyperbolic
plane, and $q$ times the second generator, to the matrix $U=\binom{u\ \ -p}{q\ \
-\tilde{u}}$. We shall henceforth identify $A^{-} \oplus H$ with the space of
those matrices. The norm $|U|^{2}=u\tilde{u}-pq$ (see Lemma \ref{vnorm6})
resembles the ``Moore-like determinant'' of $M_{2}(B)^{-}$ from the proof of
Corollary \ref{NAvn2}, but now with a subset of a bi-quaternion algebra rather
than a usual quaternion algebra.

Let $\widehat{GSp}_{A}\binom{1\ \ \ \ 0}{0\ \ -1}$ be the group of (invertible)
matrices $\binom{a\ \ b}{c\ \ d} \in M_{2}(A)$ which preserve the space
$\mathbb{F}\binom{1\ \ \ \ 0}{0\ \ -1}$ under the operation $\binom{a\ \ b}{c\ \
d}:M\mapsto\binom{a\ \ b}{c\ \ d}M\binom{\ \ \overline{d}\ \
-\overline{b}}{-\overline{c}\ \ \ \ \overline{a}}$ on $M_{2}(A)$. An element
$\binom{a\ \ b}{c\ \ d} \in M_{2}(A)$ lies in $\widehat{GSp}_{A}\binom{1\ \ \ \
0}{0\ \ -1}$ if and only if $a\overline{b}$ and $c\overline{d}$ lie in $A^{-}$
and there exists an element $m\in\mathbb{F}^{\times}$ such that
$a\overline{d}+b\overline{c}=m$ (and equivalently
$d\overline{a}+c\overline{b}=m$). The fact that our matrix $\binom{1\ \ \ \
0}{0\ \ -1}$ equals its inverse implies that if $\binom{a\ \ b}{c\ \
d}\in\widehat{GSp}_{A}\binom{1\ \ \ \ 0}{0\ \ -1}$ then so is $\binom{\ \
\overline{d}\ \ -\overline{b}}{-\overline{c}\ \ \ \ \overline{a}}$ (with the
same multiplier $m$), so that $\overline{c}a$ and $\overline{d}b$ belong to
$A^{-}$ and $\overline{a}d+\overline{c}b=m$ as well as
$\overline{d}a+\overline{b}c=m$. We call these relations \emph{the $GSp$
relations}, and the map taking an element of $\widehat{GSp}_{A}\binom{1\ \ \ \
0}{0\ \ -1}$ to the scalar $m$ by which its action multiplies $\binom{1\ \ \ \
0}{0\ \ -1}$ is a group homomorphism $\widehat{GSp}_{A}\binom{1\ \ \ \
0}{0\ \ -1}\to\mathbb{F}^{\times}$. Now, any element $\binom{a\ \ b}{c\ \ d}$ of
$\widehat{GSp}_{A}\binom{1\ \ \ \ 0}{0\ \ -1}$, with multiplier $m$, satisfies
$N^{M_{2}(A)}_{\mathbb{F}}\binom{a\ \ b}{c\ \ d}^{2}=m^{8}$ (for the reduced
norm of the degree 8 algebra $M_{2}(A)$ over $\mathbb{F}$). The function
$\frac{N^{M_{2}(A)}_{\mathbb{F}}}{m^{4}}$ is thus a group homomorphism
$\widehat{GSp}_{A}\binom{1\ \ \ \ 0}{0\ \ -1}\to\{\pm1\}$, and we define
$GSp_{A}\binom{1\ \ \ \ 0}{0\ \ -1}$ to be it kernel. We shall see in Lemma
\ref{genie} below that unless $A=M_{4}(\mathbb{F})$, the latter homomorphism is
trivial and $\widehat{GSp}_{A}\binom{1\ \ \ \ 0}{0\ \ -1}$ and $GSp_{A}\binom{1\
\ \ \ 0}{0\ \ -1}$ coincide, so that the reduced norm condition becomes
redundant. The kernel of the restriction of $m$ to $\widehat{GSp}_{A}\binom{1\ \
\ \ 0}{0\ \ -1}$ is just the symplectic group $Sp_{A}\binom{1\ \ \ \ 0}{0\ \
-1}$ consisting of those matrices which preserve the element $\binom{1\ \ \ \
0}{0\ \ -1}$ itself under this operation, and have reduced norm 1.

We begin our analysis of this group with the following
\begin{lem}
Let $\binom{a\ \ b}{c\ \ d}$ be an element of $\widehat{GSp}_{A}\binom{1\ \ \ \
0}{0\ \ -1}$, with multiplier $m$. If $a \in A^{\times}$ then $\binom{a\ \ b}{c\
\ d}$ is the product $\binom{1\ \ 0}{\beta\ \ 1}\binom{a\ \ \ \ 0\ \ \ }{0\ \
m\overline{a}^{-1}}\binom{1\ \ \alpha}{0\ \ 1}$ with $\alpha$ and $\beta$ from
$A^{-}$. Moreover, these matrices lie in $GSp_{A}\binom{1\ \ \ \ 0}{0\ \ -1}$.
\label{invent}
\end{lem}

\begin{proof}
$b\overline{a}$ and $\overline{a}c$ lie in $A^{-}$ by the $GSp$ relations, and
since $a \in A^{\times}$ we find that
$\alpha=a^{-1}b=a^{-1}(b\overline{a})\overline{a}^{\ -1}$ and
$\beta=ca^{-1}=\overline{a}^{\ -1}(\overline{a}c)a^{-1}$ are also in $A^{-}$.
Hence $b=a\alpha$ and $c=\beta a$, and as $d\overline{a}+c\overline{b}=m$ by the
$GSp$ relations, it follows that $d=m\overline{a}^{-1}+\beta a\alpha$. But this
is easily seen to be the value of the asserted product. We claim that the
multipliers lie in $GSp_{A}\binom{1\ \ \ \ 0}{0\ \ -1}$, which will prove the
second assertion. In order to evaluate the reduced norms we may extend scalars
to a splitting field of $A$, and then we evaluate $8\times8$ determinants. Now,
the unipotent elements become unipotent $8\times8$ matrices, hence they have
determinant 1, in correspondence with their multipliers being 1. The diagonal
matrix is a block matrix, hence has reduced norm $N^{A}_{\mathbb{F}}(a) \cdot
N^{A}_{m\mathbb{F}}(\overline{a}^{-1})=m^{4}$, which proves the claim as the
multiplier of this element is $m$. This completes completes the proof of the
lemma.
\end{proof}

The following technical result will be very useful in what follows.
\begin{lem}
Given any element $\binom{a\ \ b}{c\ \ d} \in GSp_{A}\binom{1\ \ \ \ 0}{0\ \
-1}$, there exists $v \in A^{-}$ such that the result of multiplying it from the
left by a matrix of the sort product $\binom{1\ \ v}{0\ \ 1}$ has an invertible
upper left entry. In fact, this is the case for ``almost any $v \in A^{-}$''. In
addition, $GSp_{A}\binom{1\ \ \ \ 0}{0\ \ -1}$ has index 2 in
$\widehat{GSp}_{A}\binom{1\ \ \ \ 0}{0\ \ -1}$ if and only if $A$ is split.
\label{genie}
\end{lem}

\begin{proof}
The upper left entry of the product in question is $a+vc$. By fixing an
arbitrary non-zero $w \in A^{-}$, consider the expression
$N^{A}_{\mathbb{F}}(a+swc)$ as a function of $s\in\mathbb{F}$. It is a
polynomial of degree not exceeding 4 in $s$. Now, if $a \in A^{\times}$ then
this polynomial never vanishes at $s=0$. Hence for all but at most 4 (non-zero)
values of $s$ (the roots of this polynomial), $v=sw$ has the desired property
(in fact, this polynomial has at most two roots, as it decomposes as the
product of the global constant $N^{A}_{\mathbb{F}}(a)$ with the polynomial
$N^{A}_{\mathbb{F}}(1+swca^{-1})$ and the latter expression is a square by the
proof of Lemma \ref{invent} and Lemma \ref{normsq} below). In case $a=0$ we know
that $c$ must be invertible (for the $GSp$ relation
$a\overline{d}+b\overline{c}=mI$ to be possible), and then with any anisotropic
$v$ the element $a+vc=vc$ is invertible. This covers the case where $A$ is a
division algebra (since then either $a=0$ or $a \in A^{\times}$), so assume that
$A=M_{2}(B)$ for $B$ a quaternion algebra over $\mathbb{F}$. We only have to
consider the case where $a$ is non-zero singular $2\times2$ matrices over $B$.
Hence there exist $\sigma$ and $\tau$ in $B$, not both zero, such that
$a\binom{\sigma}{\tau}=0$ as a 2-vector over $B$. As $a\neq0$, this allows us to
construct some $y \in GL_{2}(B)$ such that $ay$ has right column 0. It is then
easy to find $x \in GL_{2}(B)$ such that $xay=\binom{1\ \ 0}{0\ \ 0}$, and by
multiplying $\binom{a\ \ b}{c\ \ d}$ from the left by $\binom{y\ \ \ 0\ \ }{0\ \
\overline{y}^{-1}}$ and from the right by $\binom{x\ \ \ 0\ \ }{0\ \
\overline{x}^{-1}}$ we may consider only elements with $a=\binom{1\ \ 0}{0\ \
0}$. Indeed, the right multiplication by $\binom{y\ \ \ 0\ \ }{0\ \
\overline{y}^{-1}}$ does not affect our assertions, and conjugating $\binom{1\ \
v}{0\ \ 1}$ by $\binom{x\ \ \ 0\ \ }{0\ \ \overline{x}^{-1}}$ yields just
$\binom{1\ \ xv\overline{x}}{0\ \ \ 1\ }$, with $xv\overline{x} \in A^{-}$.
Write $c$ as $\binom{\lambda\ \ \mu}{\kappa\ \ \nu}$, with entries from $B$. The
$GSp$ condition $\overline{c}a \in A^{-}$ implies $\nu=0$ and
$\kappa\in\mathbb{F}$, and as $\overline{d}a+\overline{b}c=mI$ is invertible, we
find that $\mu \in B^{\times}$. Choose $w$ of the sort $\binom{\sigma\ \ -r}{1\
\ -\overline{\sigma}}$, and consider the polynomial in $s$ defined above. As
$a+swc=\binom{1+s\sigma\lambda+sr\kappa\ \ s\sigma\mu}{\
s\lambda-s\kappa\overline{\sigma}\ \ \ \ \ s\mu}$, one may use Lemma \ref{Vpol}
to evaluate the coefficients of the powers of $s$ in the resulting expression
from the formula from Proposition \ref{NAexp}, which gives us
$s^{2}N^{B}_{\mathbb{F}}(\mu)+2s^{3}\kappa
N^{B}_{\mathbb{F}}(\mu)|w|^{2}+s^{4}N^{A}_{\mathbb{F}}(w)N^{A}_{\mathbb{F}}
(c)$. Evaluating $N^{A}_{\mathbb{F}}(c)$ as $\kappa^{2}N^{B}_{\mathbb{F}}(\mu)$
(by Proposition \ref{NAexp}) and using Corollary \ref{NAvn2}, this polynomial is
(up the the global scalar $N^{B}_{\mathbb{F}}(\mu)\in\mathbb{F}^{\times}$) just
$s^{2}(1+\kappa|w|^{2}s)^{2}$. Hence our upper entry is invertible for all
non-zero $s$ if $w$ is isotropic or $c$ is not invertible (i.e., $\kappa=0$),
and we also have to omit the value $s=-\frac{1}{\kappa|w|^{2}}$ otherwise. As
these are at most two values of $s$ and $ch\mathbb{F}\neq2$, there is at least
one multiple of $w$ which we yields an invertible $a+swc$. Note that we have
only assumed that our matrix lies in $\widehat{GSp}_{A}\binom{1\ \ \ \ 0}{0\ \
-1}$, so that the last assertion of Lemma \ref{invent} implies that
$\widehat{GSp}_{A}\binom{1\ \ \ \ 0}{0\ \ -1}=GSp_{A}\binom{1\ \ \ \ 0}{0\ \
-1}$ whenever $A$ is not split.

Consider now the case $B=M_{2}(\mathbb{F})$ and $A=M_{4}(\mathbb{F})$. By an
argument similar to the  one given above, we may restrict attention to the case
where $a$ is $\binom{I_{k}\ \ 0}{0\ \ \ 0}$ for some $1 \leq k\leq3$. If $k=2$
then $c=\binom{\alpha\ \ \beta}{\gamma\ \ \delta}$ once again (with entries from
$M_{2}(\mathbb{F})$), the $GSp$ condition imply $\delta=0$,
$\gamma\in\mathbb{F}I$, and $\beta \in GL_{2}(\mathbb{F})$, and the argument
from the case of $B$ division works equally well. In particular, all these
elements (as well as the elements containing invertible entries) lie in
$GSp_{A}\binom{1\ \ \ \ 0}{0\ \ -1}$. We now claim that the cases $k=3$ and
$k=1$ can only occur for elements of $\widehat{GSp}_{A}\binom{1\ \ \ \ 0}{0\ \
-1}$ which are not in $GSp_{A}\binom{1\ \ \ \ 0}{0\ \ -1}$. First we demonstrate
the existence of elements of $\widehat{GSp}_{A}\binom{1\ \ \ \ 0}{0\ \ -1}$
which are not in $GSp_{A}\binom{1\ \ \ \ 0}{0\ \ -1}$: One example is the matrix
$\binom{e\ \ f}{g\ \ h}$ in which $e=\binom{I_{3}\ \ 0}{0\ \ \ 0}$, $h=\binom{0\
\ 0\ }{0\ \ I_{3}}$, $f$ has only lower left entry 1 and other 15 entries
vanish, and $g$ has upper right entry 1 (and the rest 0), with multiplier 1 and
determinant $-1$. Consider now the case $k=3$, and write $\binom{\alpha\ \
\beta}{\gamma\ \ \delta}$ with the entries from $M_{2}(\mathbb{F})$. The
condition that $\overline{a}c=\binom{0\ \ 0\ }{0\ \ I_{3}}\binom{\alpha\ \
\beta}{\gamma\ \ \delta}$ lies in $A^{-}$ implies, in particular, that in the
rightmost column of $c$ the only entry which may be non-zero is the upper right
one, which we denote by $t$. $t$ may not vanish, for otherwise the matrix
$\overline{d}a+\overline{b}c=mI$ would have to be singular, a contradiction .
But now left multiplication from our representative of
$\widehat{GSp}_{A}\binom{1\ \ \ \ 0}{0\ \ -1}/GSp_{A}\binom{1\ \ \ \ 0}{0\ \
-1}$ would yield a matrix in $M_{2}(A)$ whose upper left entry is
$\binom{I_{3}\ \ 0}{*\ \ \ t}$ (where $(*\ \ t)$ is the upper row of $c$). As
this matrix is invertible, the resulting element lies in $GSp_{A}\binom{1\ \ \ \
0}{0\ \ -1}$, hence the original one $\binom{a\ \ b}{c\ \ d}$ does not. In the
case $k=1$, if none of the entries are invertible then $c$ must have rank 3
(again, for $\overline{d}a+\overline{b}c=mI$ to be non-singular), and again we
are in the case $k=3$. This completes the proof of the lemma.
\end{proof}

\begin{cor}
Any element of $GSp_{A}\binom{1\ \ \ \ 0}{0\ \ -1}$, with multiplier $m$, can be
written as $\binom{1\ \ v}{0\ \ 1}\binom{1\ \ 0}{\beta\ \ 1}\binom{a\ \ \ \ 0\ \
\ }{0\ \ m\overline{a}^{-1}}\binom{1\ \ \alpha}{0\ \ 1}$ for some $a \in
A^{\times}$, and $u$, $\alpha$, and $\beta$ from $A^{-}$. \label{genform}
\end{cor}

\begin{proof}
This follows directly from Lemmas \ref{invent} and \ref{genie}.
\end{proof}

Many arguments below shall make use of the following
\begin{lem}
For $\eta$ and $\omega$ in $A^{-}$, define
$D(\eta,\omega)=1+2\langle\eta,\tilde{\omega}\rangle+|\omega|^{2}|\eta|^{2}$.
Then the element $1+\eta\omega$ of $A^{-}$ has reduced norm
$D(\eta,\omega)^{2}$, and its product with $1+\tilde{\omega}\tilde{\eta}$ (from
either side) yields the scalar $D(\eta,\omega)$. \label{normsq}
\end{lem}

\begin{proof}
The fact that the products $(1+\eta\omega)(1+\tilde{\omega}\tilde{\eta})$ and
$(1+\tilde{\omega}\tilde{\eta})(1+\eta\omega)$ both yield $D(\eta,\omega)$
easily follow from Lemma \ref{vnorm6}. For the reduced norm, multiply our
element by $\tilde{\omega}$ from the right. The result is
$\tilde{\omega}+|\omega|^{2}\eta$ (Lemma \ref{vnorm6} again) and lies in
$A^{-}$, so that its reduced norm is $|\tilde{\omega}+|\omega|^{2}\eta|^{2}$ by
Corollary \ref{NAvn2}. Moreover, Lemma \ref{Vpol} evaluates this vector norm as
$|\omega|^{2}D(\eta,\omega)$, so that this proves the assertion for anisotropic
$\omega$. Assume now $|\omega|^{2}=0$. Take some $\xi \in A^{-}$ with
$|\xi|^{2}\neq0$, and consider the two expressions
$N^{A}_{\mathbb{F}}\big(1+\eta(\omega+s\xi)\big)$ and
$D(\eta,(\omega+s\xi))^{2}$ for $s\in\mathbb{F}$. Both are polynomials of degree
4 in $s$, which were seen to coincide wherever $|(\omega+s\xi)|^{2}\neq0$. The
latter assumption occurs for any $s$ other than $s=0$ or
$s=-\frac{2\langle\omega,\xi\rangle}{|\xi|^{2}}$ (Lemma \ref{Vpol} again and the
isotropy of $\omega$), so that we omit at most two values. By extending scalars
if necessary, we may assume that $\mathbb{F}$ has more than 6 elements. Then we
have two polynomials of degree 4 which coincide on more than 4 elements of
$\mathbb{F}$, hence they must be the same polynomial. Substituting $s=0$ in the
two equals polynomials verify the assertion for isotropic $\omega$ as well. This
completes the proof of the lemma.
\end{proof}

The freedom of choice we have in Lemma \ref{genie} shows that there are many
different choices of parameters to get the same element of $GSp_{A}\binom{1\ \ \
\ 0}{0\ \ -1}$ in the form of Corollary \ref{genform}. Hence some compatibility
assertions will be required wherever we use the form from that Corollary. these
will be based on the following
\begin{lem}
Assume that the expression using $a$, $v$, $\alpha$ and $\beta$ in Corollary
\ref{genform} yields the same element (of multiplier $m$) as the one arising
from $c$, $w$, $\gamma$, and $\delta$ respectively. Then we have the equalities
$(i)$ $c=\big(1-(w-v)\beta\big)a$, $(ii)$
$\delta=\frac{\beta-|\beta|^{2}(\tilde{w}-\tilde{v})}{D(\beta,w-v)}$, and
$(iii)$
$\gamma=\alpha-ma^{-1}\frac{w-v-|w-v|^{2}\tilde{\beta}}{D(\beta,w-v)}\overline{a
}^{\ -1}$. \label{comp}
\end{lem}

\begin{proof}
The product from Corollary \ref{genform} equals $\binom{(1+v\beta)a\ \
(1+v\beta)a\alpha+mv\overline{a}^{-1}}{\ \beta a\ \ \ \ \ \ \ \ \beta
a\alpha+m\overline{a}^{-1}\ }$ (matrix multiplication). When we compare this
matrix with the one arising from $c$, $w$, $\gamma$, and $\delta$, we first
find that $\beta a=\delta c$ and $(1+v\beta)a=(1+w\delta)c$. Multiplying the
first equality by $w$ from the left and subtracting the result from the second
equality establishes $(i)$. We write $\delta=\beta ac^{-1}$, substitute $c$ from
part $(i)$, use Lemma \ref{normsq} for the inverse of $1-(w-v)\beta$, and apply
Lemma \ref{vnorm6}, which proves part $(ii)$. We can now write $a\alpha$ as the
upper right entry of our common matrix minus $v$ times the lower right entry and
the same for $c\gamma$ (but with $w$). Comparing yields
$c\gamma=a\alpha-(w-v)(\beta a\alpha+m\overline{a}^{-1})$, which equals
$c\alpha-m(w-v)\overline{a}^{-1}$. Multiplying by $c^{-1}$ from the left and
using part $(i)$ and Lemmas \ref{normsq} and \ref{vnorm6} again yields part
$(iii)$. This proves the lemma.
\end{proof}

In addition, the description of the product in the parameters from Corollary
\ref{genform} is given in the following
\begin{lem}
Given two elements $g$ and $h$ of $GSp_{A}\binom{1\ \ \ \ 0}{0\ \ -1}$, with
multipliers $m$ and $n$ respectively, assume that $v \in A^{-}$ is such that
left multiplication of both $g$ and $gh$ by $\binom{1\ \ -v}{0\ \ \ \ 1}$ yields
matrices with invertible upper left entry. Then if $a$, $v$, $\alpha$ and
$\beta$ are the parameters thus obtained for $g$ as in Corollary \ref{genform}
and $e$, $z$, $\kappa$, and $\nu$ are parameters for $h$, then parameters for
$gh$ may be taken as $x$, $v$, $\xi$, and $\zeta$ with
$x=a\big(1+(\alpha+z)\nu\big)e$,
$\xi=\kappa+ne^{-1}\frac{\alpha+z+|\alpha+z|^{2}\tilde{\nu}}{D(\alpha+z,\nu)}
\overline{e}^{\ -1}$, and
$\zeta=\beta+m\overline{a}^{\
-1}\frac{\nu+|\nu|^{2}(\tilde{\alpha}+\tilde{z})}{D(\alpha+z,\nu)}a^{-1}$.
\label{GSpprod}
\end{lem}

\begin{proof}
Comparing the expressions for the product shows that the expression $\binom{a\ \
\ \ 0\ \ \ }{0\ \ m\overline{a}^{-1}}\binom{1\ \ \alpha+z}{0\ \ \ \ 1\ \
}\binom{1\ \ 0}{\nu\ \ 1}\binom{e\ \ \ \ 0\ \ \ }{0\ \ n\overline{e}^{-1}}$
equals $\binom{\ \ 1\ \ \ \ 0}{\zeta-\beta\ \ 1}\binom{x\ \ \ \ \ 0\ \ \
}{0\ \ mn\overline{x}^{-1}}\binom{1\ \ \xi-\kappa}{0\ \ \ \ 1\ \ }$. When we
consider the upper left entries of both sides we obtain the asserted value for
$x$. The values for $\xi$ (resp. $\zeta$) is obtained by comparing the upper
right (resp. lower left) entries, using the value of $x$, Lemma \ref{normsq},
and Lemma \ref{vnorm6}. This proves the lemma.
\end{proof}

Now, the group which will end up being the Gspin group is, in general, not the
full group $GSp_{A}\binom{1\ \ \ \ 0}{0\ \ -1}$, but a double cover of a certain
subgroup. This subgroup is defined using the following
\begin{prop}
The map $\varphi$ which takes an element of $GSp_{A}\binom{1\ \ \ \ 0}{0\ \
-1}$, decomposes it as in Corollary \ref{genform}, and sends it to the image of
$N^{A}_{\mathbb{F}}(a)$ in $\mathbb{F}^{\times}/(\mathbb{F}^{\times})^{2}$ is a
well-defined group homomorphism $\varphi:GSp_{A}\binom{1\ \ \ \ 0}{0\ \
-1}\to\mathbb{F}^{\times}/(\mathbb{F}^{\times})^{2}$. \label{GSphom}
\end{prop}

\begin{proof}
Given two decompositions of the same element of $GSp_{A}\binom{1\ \ \ \ 0}{0\ \
-1}$, part $(i)$ of Lemma \ref{comp} shows that the reduced norms used to
define $\varphi$ in these decompositions differ by the reduced norm of an
element of the form considered in Lemma \ref{normsq}. Hence, by that lemma, the
result in $\mathbb{F}^{\times}/(\mathbb{F}^{\times})^{2}$ is the same for both
decompositions. Hence $\varphi$ is well-defined. Now, given two elements of
$GSp_{A}\binom{1\ \ \ \ 0}{0\ \ -1}$, the level of freedom in Lemma \ref{genie}
allows us to find $v \in A^{-}$ for which left multiplication by $\binom{1\ \
-v}{0\ \ \ \ 1}$ renders the upper left entries of both $g$ and $gh$ invertible.
We then invoke Lemma \ref{GSpprod}, and using Lemma \ref{normsq} for the reduced
norm of the multiplier between $a$ and $e$ in $x$ there shows that $\varphi$ is
also multiplicative. This proves the lemma.
\end{proof}

Note that the choice of the upper left entry in Proposition \ref{GSphom} is
arbitrary, but does not affect the value of $\varphi$ in the sense that a
similar definition in terms of another entry would yield the same result.
Indeed, $\varphi$ is a group homomorphism (Proposition \ref{GSphom}) which
attains 1 on the unipotent matrices and $GSp_{A}\binom{1\ \ \ \ 0}{0\ \ -1}$
contains the element $\binom{0\ \ 1}{1\ \ 0}$, which has multiplier 1. As the
latter element may be obtained from the parameters $v=-u$, $a=u$, and
$\alpha=\beta=\frac{\tilde{u}}{|u|^{2}}$ in Corollary \ref{genform} for some
anisotropic $u$ (by Lemma \ref{vnorm6}), Corollary \ref{NAvn2} shows that
$\varphi\binom{0\ \ 1}{1\ \ 0}$ is also trivial. It thus suffices to compare the
reduced norms of the entries $b$, $c$, and $d$ in the matrix $\binom{a\ \ b}{c\
\ d}$ when given in the form of Lemma \ref{invent}, and see that if one of them
is invertible then it has the same reduced norm as $a$ up to
$(\mathbb{F}^{\times})^{2}$. For $b$ and $c$ the assertion follows directly from
Corollary \ref{NAvn2}, while $\frac{\overline{a}d}{m}$ (or
$\frac{d\overline{a}}{m}$) is an element of the form given in Lemma
\ref{normsq}. Hence $\varphi$ is more intrinsic than it seems at first sight. In
particular, in case an element of $GSp_{A}^{\mathbb{F}^{2}}\binom{1\ \ \ \ 0}{0\
\ -1}$ has any invertible entry, we may use the reduced norm of that entry in
order to evaluate the $\varphi$-image of that element.

\smallskip

Denote the kernel of the map $\varphi$ from Proposition \ref{GSphom} by
$GSp_{A}^{\mathbb{F}^{2}}\binom{1\ \ \ \ 0}{0\ \ -1}$. We now construct a
certain group automorphism of a double cover of the subgroup
$GSp_{A}^{\mathbb{F}^{2}}\binom{1\ \ \ \ 0}{0\ \ -1}$, which is again based on
some choice of square root. For the definition we shall use
\begin{lem}
Let $a$, $v$, $\alpha$, $\beta$, $c$, $w$, $\gamma$, and $\delta$ as in Lemma
\ref{comp}, and assume that the (common) element lies in
$GSp_{A}^{\mathbb{F}^{2}}\binom{1\ \ \ \ 0}{0\ \ -1}$. Let
$t\in\mathbb{F}^{\times}$ be such that $N^{A}_{\mathbb{F}}(a)=t^{2}$, and denote
$tD(\beta,w-v)$ by $s$. The following expressions remain invariant by replacing
$a$ by $c$, $v$ by $w$, $\alpha$ by $\gamma$, $\beta$ by $\delta$, and $t$ by
$s$: $(i)$ $t\tilde{\beta}\overline{a}^{\ -1}$. $(ii)$
$t(1+\tilde{v}\tilde{\beta})\overline{a}^{\ -1}$. $(iii)$
$t\tilde{\beta}\overline{a}^{\ -1}\tilde{\alpha}+\frac{m}{t}a$. $(iv)$
$t(1+\tilde{v}\tilde{\beta})\overline{a}^{\
-1}\tilde{\alpha}+\frac{m}{t}\tilde{v}a$. \label{psiwd}
\end{lem}

\begin{proof}
Lemma \ref{normsq} and part $(i)$ of Lemma \ref{comp} yield $\overline{c}^{\
-1}=\frac{1-(\tilde{w}-\tilde{v})\tilde{\beta}}{D(\beta,w-v)}\overline{a}^{\
-1}$. Part $(i)$ then follows from the definition of $s$, Lemma \ref{normsq},
and part $(ii)$ of Lemma \ref{comp}. Part $(ii)$ is now obtained from part
$(i)$, the latter equation, the definition of $s$, and simple algebra. Now, part
$(iii)$ of Lemma \ref{comp}, Lemma \ref{AxA-rel}, and the assumption on $t$
imply $s\overline{c}^{\
-1}\tilde{\gamma}=t\big(1-(\tilde{w}-\tilde{v})\tilde{\beta}\big)\overline{}^{\
-1}\tilde{\alpha}-\frac{m}{t}(\tilde{w}-\tilde{v})a$. Part $(iii)$ is
established using part $(ii)$ of Lemma \ref{comp}, Lemma \ref{normsq}, and the
definition of $s$. Part $(iv)$ now follows from the latter equality and part
$(iii)$. This completes the proof of the lemma.
\end{proof}

By definition, an element of $GSp_{A}\binom{1\ \ \ \ 0}{0\ \ -1}$ belongs to
$GSp_{A}^{\mathbb{F}^{2}}\binom{1\ \ \ \ 0}{0\ \ -1}$ if and only if when
decomposed as in Corollary \ref{genform}, the diagonal matrix has entries from
$A^{(\mathbb{F}^{\times})^{2}}$. Using the double cover
$\widetilde{A}^{(\mathbb{F}^{\times})^{2}}$ from Lemma \ref{ac6d1} we now have
\begin{thm}
The group $GSp_{A}^{\mathbb{F}^{2}}\binom{1\ \ \ \ 0}{0\ \ -1}$ admits a
well-defined double cover $\widetilde{GSp}_{A}^{\mathbb{F}^{2}}\binom{1\ \ \ \
0}{0\ \ -1}$, in which the parameter $a \in A^{(\mathbb{F}^{\times})^{2}}$ from
Corollary \ref{genform} is replaced by an element
$(a,t)\in\widetilde{A}^{(\mathbb{F}^{\times})^{2}}$ lying over it. Define a map
$\psi$ from $\widetilde{GSp}_{A}^{\mathbb{F}^{2}}\binom{1\ \ \ \ 0}{0\ \ -1}$
to itself by replacing the paramter $(a,t)$ by $\big(t\overline{a}^{\
-1},t\big)$ and sending the other paramters from Corollary \ref{genform} to
their $\theta$-images. Then $\psi$ is a well-defined group automorphism of order
2 of $\widetilde{GSp}_{A}^{\mathbb{F}^{2}}\binom{1\ \ \ \ 0}{0\ \ -1}$, which
commutes with the multiplier map to $\mathbb{F}^{\times}$. \label{GSppsi}
\end{thm}

\begin{proof}
Assume, as in Lemma \ref{comp}, that two sets of paramters in Corollary
\ref{genform}, say $a$, $v$, $\alpha$ and $\beta$ versus $c$, $w$, $\gamma$, and
$\delta$, describe the same element of $GSp_{A}^{\mathbb{F}^{2}}\binom{1\ \ \ \
0}{0\ \ -1}$. Lemma \ref{normsq} and part $(i)$ of Lemma \ref{comp} show that
the same element of the double cover
$\widetilde{GSp}_{A}^{\mathbb{F}^{2}}\binom{1\ \ \ \ 0}{0\ \ -1}$ is obtained
from $(a,t)$, $v$, $\alpha$ and $\beta$ and from $\big(c,D(v-w,\beta)t\big)$,
$w$, $\gamma$, and $\delta$. Hence this double cover is well-defined. Consider
now our element of $GSp_{A}^{\mathbb{F}^{2}}\binom{1\ \ \ \ 0}{0\ \ -1}$ written
in terms of the parameters $(a,t)$, $v$, $\alpha$ and $\beta$. Its matrix form
appears at the beginning of the proof of Lemma \ref{comp}, and one sees that its
$\psi$-image (using these parameters) has precisely the four entries which
appear in the various parts of Lemma \ref{psiwd}. But this lemma precisely shows
that taking the set of parameters $\big(c,D(v-w,\beta)t\big)$, $w$, $\gamma$,
and $\delta$ instead yields the same entries. This shows that $\psi$ is
well-defined, and its image is clearly in
$\widetilde{GSp}_{A}^{\mathbb{F}^{2}}\binom{1\ \ \ \ 0}{0\ \ -1}$ again and has
the same multiplier. It is a map of order 2 since so are $\theta$ and the map
$(a,t)\mapsto\big(t\overline{a}^{\ -1},t\big)$ on
$\widetilde{A}^{(\mathbb{F}^{\times})^{2}}$. Finally, let $(e,r)$, $z$,
$\kappa$, and $\nu$ be parameters of another element of
$\widetilde{GSp}_{A}^{\mathbb{F}^{2}}\binom{1\ \ \ \ 0}{0\ \ -1}$ and assume
that $v$ represents a parameter also for the product of these two elements.
Then Lemma \ref{GSpprod} provides expressions for the parameters $x$, $v$,
$\xi$, and $\zeta$ of the product in $GSp_{A}^{\mathbb{F}^{2}}\binom{1\ \ \ \
0}{0\ \ -1}$, and Lemma \ref{normsq} shows that we may replace $x$ by
$\big(x,tD(a+z,\nu)r\big)$ for the parameter of the product in the double cover
$\widetilde{GSp}_{A}^{\mathbb{F}^{2}}\binom{1\ \ \ \ 0}{0\ \ -1}$. Showing that
$\psi$ is multiplicative amounts to verifying that sending $(a,t)$ to
$\big(t\overline{a}^{\ -1},t\big)$, $(e,r)$ to
$\big(r\overline{e}^{\ -1},r\big)$, and the $A^{-}$-parameters $v$, $\alpha$,
$\beta$, $z$, $\kappa$, and $\nu$ to their $\theta$-images results in the same
effect on $\big(x,tD(a+z,\nu)r\big)$ and on $\xi$ and $\zeta$ (the parameter
$v$ of the product already appears). Now, for $x$ this follows directly from
Lemma \ref{normsq} and the multiplicativity of $y\mapsto\overline{y}^{\ -1}$ on
$A^{\times}$, and for $\xi$ and $\zeta$ this follows from Lemma \ref{AxA-rel},
the assumptions on $r$ and $t$, the preservation of multipliers, and the fact
that $\theta$ preserves the bilinear form on $A^{-}$ (this is relevant also for
the action on the denominators $D(\eta,\omega)$). This completes the proof of
the theorem.
\end{proof}

As with $\varphi$, the automorphism $\psi$ from Theorem \ref{GSppsi} may be
defined using the other entries, hence is a more intrinsic automorphism of
$\widetilde{GSp}_{A}^{\mathbb{F}^{2}}\binom{1\ \ \ \ 0}{0\ \ -1}$ that one might
think at first. To see this, observe that the element $\binom{0\ \ 1}{1\ \ 0}$,
of multiplier 1, equals its $\psi$-image with the appropriate choice of square
root: Indeed, it may be obtained from a set of parameters arising from an
anisotropic vector $u$, and the result is independent of $u$. By taking the
square root $-|u|^{2}$ for the reduced norm of $a=u$ (Corollary \ref{NAvn2}
again), we find that $\psi$ just replaces every instance of $u$ by $\tilde{u}$,
and the resulting matrix must therefore be $\binom{0\ \ 1}{1\ \ 0}$ again. Hence
we may use the same argument as for $\varphi$ in order to obtain that $\psi$ may
be evaluated, for example, by applying $(g,t)\mapsto\big(t\overline{g}^{\
-1},t\big)$ to any invertible entry of the matrix in question. However, we shall
stick to our form of Corollary \ref{genform} in what follows.

\smallskip

The relation between the group $GSp_{A}^{\mathbb{F}^{2}}\binom{1\ \ \ \ 0}{0\ \
-1}$ and the space $A^{-} \oplus H$ from Lemma \ref{sp8id1} (embedded in
$M_{2}(A)$ as described above) begins to reveal itself in the following
\begin{lem}
Any anisotropic vector $U \in A^{-} \oplus H$ lies also in
$GSp_{A}^{\mathbb{F}^{2}}\binom{1\ \ \ \ 0}{0\ \ -1}$, with multiplier
$|U|^{2}$. The involution $\hat{\psi}=Id_{H}\oplus\theta$ of $A^{-} \oplus H$
coincides, on anisotropic elements, with the map $\psi$ from Theorem
\ref{GSppsi}, after an appropriate lift of these elements into the double cover
$\widetilde{GSp}_{A}^{\mathbb{F}^{2}}\binom{1\ \ \ \ 0}{0\ \ -1}$. The products
$U\hat{\psi}(U)$ and $\hat{\psi}(U)U$ both equal the scalar $|U|^{2}$, and the
for pairing of two vectors $U$ and $V$ in $A^{-} \oplus H$ we have $2\langle U,V
\rangle=U\hat{\psi}(V)+V\hat{\psi}(U)=\hat{\psi}(U)V+\hat{\psi}(V)U$.
\label{vnorm8}
\end{lem}

\begin{proof}
For the first assertion we evaluate $\binom{u\ \ -p}{q\ \ -\tilde{u}}\binom{1\ \
\ \ 0}{0\ \ -1}\binom{\ \ \tilde{u}\ \ \ \ p}{-q\ \ -u}$ (recall that $u$ and
$\tilde{u}$ are in $A^{-}$), and the result is indeed $|U|^{2}\binom{1\ \ \ \
0}{0\ \ -1}$ for $U=\binom{u\ \ -p}{q\ \ -\tilde{u}}$. In case $p\neq0$, we
multiply $U$ by $\binom{1\ \ \ \ 0}{0\ \ -1}$ (the simplest element with
multiplier $-1$, clearly equals its $\psi$-image) and by $\binom{0\ \ 1}{1\ \
0}$ from the right. The resulting element has multiplier $pq-|u|^{2}$, and it
may be obtained by taking the parameters $a=p$ (with the square root $p^{2}$ of
$N^{A}_{\mathbb{F}}(p)=p^{4}$), $v=0$, $\alpha=\frac{u}{p}$, and
$\beta=\frac{\tilde{u}}{p}$. As $p^{2}\overline{p}^{\ -1}=p$ once again, $\psi$
just applies $\theta$ to the coordinates $u$ and $\tilde{u}$, and multiplying by
$\binom{\ \ 0\ \ 1}{-1\ \ 0}$ (which was seen to equal its $\psi$-image) from
the right again proves the assertion for the case $p\neq0$. Otherwise
$|u|^{2}\neq0$, and we choose the parameters $v=\alpha=0$, $a=u$ (for the
reduced norm of which Corollary \ref{NAvn2} allows us to take $-|u|^{2}$ as a
square root), and $\beta=\frac{q\tilde{u}}{|u|^{2}}$. The resulting $\psi$-image
is once again obtained by just applying $\theta$ to $u$ and to $\tilde{u}$,
completing the verification of this assertion. Now, the products
$U\hat{\psi}(U)=\binom{u\ \ -p}{q\ \ -\tilde{u}}\binom{\tilde{u}\ \ -p}{q\ \
-u}$ and $\hat{\psi}(U)U=\binom{\tilde{u}\ \ -p}{q\ \ -u}\binom{u\ \ -p}{q\ \
-\tilde{u}}$ are easily evaluated, by Lemma \ref{vnorm6}, as $|U|^{2}I$, and the
last assertion now follows from Lemma \ref{Vpol}. This proves the lemma. 
\end{proof}

More importantly, we also have
\begin{prop}
If $U \in A^{-} \oplus H$ and
$\big(g,\psi(g)\big)\in\widetilde{GSp}_{A}^{\mathbb{F}^{2}}\binom{1\ \ \
\ 0}{0\ \ -1}$ then the matrix $gU\psi(g)^{-1}$ also lies in $A^{-} \oplus H$,
and it has the same vector norm as $U$. \label{GSppresHA-}
\end{prop}

\begin{proof}
It suffices to prove the assertion for a generating subset of
$\widetilde{GSp}_{A}^{\mathbb{F}^{2}}\binom{1\ \ \ \ 0}{0\ \ -1}$. By Corollary
\ref{genform}, the set of diagonal and unipotent matrices is a generating set.
$\psi$ operates as $\theta$ on the $A^{-}$-coordinates of the unipotent
generators (by choosing 1 to be the square root of $N^{A}_{\mathbb{F}}(1)$),
while on the diagonal ones it operates as $(a,t)\mapsto\big(t\overline{a}^{\
-1},t\big)$. The composition with inversion as $-\theta$ on the unipotent ones
and $(a,t)\mapsto\big(\frac{\overline{a}}{t},\frac{1}{t}\big)$ on the diagonal
ones. The action of $\binom{1\ \ v}{0\ \ 1}$ thus leaves $q$ invariant, takes
$u$ to $u+qv$ and $\tilde{u}$ to $\tilde{u}+q\tilde{v}$, and $p$ to
$p+u\tilde{v}+v\tilde{u}+qv\tilde{v}$. $\binom{1\ \ 0}{w\ \ 1}$ sends $q$ to
$q+wu+\tilde{u}\tilde{w}+pw\tilde{w}$, $u$ to $u+p\tilde{w}$, $\tilde{u}$ to
$\tilde{u}+pw$, and leaves $p$ invariant. Finally, applying $\binom{a\ \ \ \ 0\
\ \ }{0\ \ m\overline{a}^{-1}}$ multiplies $q$ by $\frac{m}{t}$ and $p$ by
$\frac{t}{m}$, and maps $u$ to $\frac{au\overline{a}}{t}$, $\tilde{u}$ to
$t\overline{a}^{\ -1}\tilde{u}a^{-1}$. The image of $u$ is the the image of
$\tilde{u}$ under $\theta$ in all these cases (this is clear in the first two
operations and uses Lemma \ref{AxA-rel} and the fact that
$t^{2}=N^{A}_{\mathbb{F}}(a)$ for the latter case). In addition, the expressions
we add to $p$ in the first case and $q$ in the second case are $2\langle u,v
\rangle+q|v|^{2}$ and $2\langle u,\tilde{w} \rangle+p|w|^{2}$ by Lemma
\ref{vnorm6} respectively, multiplying which by $q$ (resp. $p$) yields
$|u+qv|^{2}-|u|^{2}$ (resp. $|u+p\tilde{w}|^{2}-|u|^{2}$) by Lemma \ref{Vpol}.
The fact that Lemma \ref{ac6d1} implies that
$\big|\frac{au\overline{a}}{t}\big|^{2}=|u|^{2}$ for the latter generators now
completes the verification of both assertions for all the necessary cases.
This proves the proposition.
\end{proof}

We shall also make use of the following
\begin{lem}
Let $U \in A^{-} \oplus H$ and
$\big(g,\psi(g)\big)\in\widetilde{GSp}_{A}^{\mathbb{F}^{2}}\binom{1\ \ \
\ 0}{0\ \ -1}$ be given. Then the equality
$\hat{\psi}\big(gU\psi(g)^{-1}\big)=\psi(g)\hat{\psi}(U)g^{-1}$ holds.
\label{hpsiGSprel}
\end{lem}

\begin{proof}
First, Proposition \ref{GSppresHA-} shows that $gU\psi(g)^{-1} \in A^{-} \oplus
H$ and it has vector norm $|U|^{2}$. In particular, its $\hat{\psi}$-image
is defined. Now, using Lemma \ref{vnorm8} we write
\[gU\psi(g)^{-1}\hat{\psi}\big(gU\psi(g)^{-1}\big)=|gU\psi(g)^{-1}|^{2}=|U|^{2}
=g|U|^{2}g^{-1}\] (since $|U|^{2}$ is a scalar), and applying Lemma \ref{vnorm8}
again the latter term can be presented as
$gU\hat{\psi}(U)g^{-1}=gU\psi(g)^{-1}\cdot\psi(g)\hat{\psi}(U)g^{-1}$. If $U$
is anisotropic (hence invertible) then so is $gU\hat{\psi}(U)g^{-1}$, and the
assertion follows for such $U$. For the rest we observe that both sides are
linear in $U$ and $A^{-} \oplus H$ is spanned by anisotropic vectors. This
completes the proof of the lemma.
\end{proof}

We can now give more details to the group action in this case:
\begin{lem}
The group $\widetilde{GSp}_{A}^{\mathbb{F}^{2}}\binom{1\ \ \ \ 0}{0\ \ -1}$ maps
to $O(A^{-} \oplus H)$ with kernel $\mathbb{F}^{\times}$, in which the choice of
the square root of the reduced norm of the coordinate $r$ of $rI$ is $r^{2}$.
Let $\tilde{\psi}$ be an element generating a group of order 2. If
$\tilde{\psi}$ operates on $\widetilde{GSp}_{A}^{\mathbb{F}^{2}}\binom{1\ \ \ \
0}{0\ \ -1}$ as the automorphism $\psi$ then sending it to the involution
$\hat{\psi}$ on $A^{-} \oplus H$ yields a map from the semi-direct product of
$\{1,\tilde{\psi}\}$ with $\widetilde{GSp}_{A}^{\mathbb{F}^{2}}\binom{1\ \ \ \
0}{0\ \ -1}$ to $O(A^{-} \oplus H)$. \label{ac8id1}
\end{lem}

\begin{proof}
Proposition \ref{GSppresHA-} defines a map
$\widetilde{GSp}_{A}^{\mathbb{F}^{2}}\binom{1\ \ \ \ 0}{0\ \ -1} \to O(A^{-}
\oplus H)$. Given $r\in\mathbb{F}^{\times}$, the scalar matrix $rI$ lies in
$GSp_{A}^{\mathbb{F}^{2}}\binom{1\ \ \ \ 0}{0\ \ -1}$ (with multiplier $r^{2}$),
and one easily verifies that it equals its $\psi$-image if the square root of
$N^{A}_{\mathbb{F}}(r)=r^{4}$ is taken to be $r^{2}$. This defines an embedding
of $\mathbb{F}^{\times}$ into $\widetilde{GSp}_{A}^{\mathbb{F}^{2}}\binom{1\ \ \
\ 0}{0\ \ -1}$, with image in the kernel of the action on $A^{-} \oplus H$ by
the centrality of such $rI$. In order to show that these are the only elements
operating trivially, let $\binom{a\ \ b}{c\ \ d}$ be an element of
$GSp_{A}^{\mathbb{F}^{2}}\binom{1\ \ \ \ 0}{0\ \ -1}$,
with multiplier $m$, and let $\binom{e\ \ f}{g\ \ h}$ describe the inverse of
its $\psi$-image (with multiplier $\frac{1}{m}$). The action sends the elements
$\binom{0\ \ 0}{1\ \ 0}$ and $\binom{0\ \ 1}{0\ \ 0}$ of $A^{-} \oplus H$ to
$\binom{be\ \ bf}{de\ \ df}$ and $\binom{ag\ \ ah}{cg\ \ ch}$ respectively. If
the action is trivial, we must have $be=bf=0$ and $cg=ch=0$. But then we get,
from the $GSp$ relations for $\binom{e\ \ f}{g\ \ h}$, the equalities
$b=mb(e\overline{h}+f\overline{g})=0$ and
$c=mc(h\overline{e}+g\overline{f})=0$, so that $a \in
A^{(\mathbb{F}^{\times})^{2}}$ with $N^{A}_{\mathbb{F}}(a)=t^{2}$,
$d=m\overline{a}^{\ -1}$, $f=g=0$, $e=\frac{\overline{a}}{t}$, and
$h=\frac{ta^{-1}}{m}$. But the action on $A^{-} \subseteq A^{-} \oplus H$ was
seen in Proposition \ref{GSppresHA-} to be via the map from Lemma \ref{ac6d1},
which shows that the only elements
$(a,t)\in\widetilde{A}^{(\mathbb{F}^{\times})^{2}}$ which act trivially are of
the form $(r,r^{2})$ with $r\in\mathbb{F}^{\times}$. In order for the action of
this element to be trivial also on $H \subseteq A^{-} \oplus H$, the formula
from Proposition \ref{GSppresHA-} implies that $m$ must be $r^{2}$ as well, so
that our element is indeed $rI$ with $\psi$-image also $rI$ (note that using the
other sign for the $\psi$-image yields elements operating as $-Id_{A^{-} \oplus
H}$). The fact that $\hat{\psi}$ clearly lies in $O(A^{-} \oplus H)$ and the
scalar $\frac{1}{m(g)}$ operates trivially now implies, together with Lemma
\ref{hpsiGSprel}, that the map to $O(A^{-} \oplus H)$ is well-defined on the
semi-direct product in question. This proves the lemma.
\end{proof}

Once again, we shall need an assertion about reflections:
\begin{lem}
An anisotropic vector $g \in A^{-} \oplus H$ may be lifted to
$\widetilde{GSp}_{A}^{\mathbb{F}^{2}}\binom{1\ \ \ \ 0}{0\ \ -1}$ such that its
$\psi$-image is $-\hat{\psi}(g)$. The map taking $U \in A^{-} \oplus H$ to the
image of $\hat{\psi}(U)$ under the action of this lift of $g$ is the reflection
in $g$. \label{ref8id1}
\end{lem}

\begin{proof}
By Lemma \ref{vnorm8}, such $g$ lies in $GSp_{A}^{\mathbb{F}^{2}}\binom{1\ \ \ \
0}{0\ \ -1}$, and $\psi(g)$ can be taken from $\{\pm\hat{\psi}(g)\}$. Hence such
a lift to $\widetilde{GSp}_{A}^{\mathbb{F}^{2}}\binom{1\ \ \ \ 0}{0\ \ -1}$
exists, and operates orthogonally on $A^{-} \oplus H$ by Proposition
\ref{GSppresHA-} (or Lemma \ref{ac8id1}). In the evaluation of the result on
$U=g$ the two factors involving $\hat{\psi}(g)$ cancel to give just $-g$. On the
other hand, if $u \in g^{\perp}$ then Lemma \ref{vnorm8} allows us to replace
$g\hat{\psi}(U)$ by $-U\hat{\psi}(g)$, and a similar argument shows that the
final result is just $U$. This proves the lemma.
\end{proof}

We can finally prove
\begin{thm}
The group $\widetilde{GSp}_{A}^{\mathbb{F}^{2}}\binom{1\ \ \ \ 0}{0\ \ -1}$ is
the Gspin group $Gspin(A^{-} \oplus H)$. It is generated by lifts of anisotropic
elements $A^{-} \oplus H$ whose $\psi$-images coincide with their
$-\hat{\psi}$-images, so that $GSp_{A}^{\mathbb{F}^{2}}\binom{1\ \ \ \ 0}{0\ \
-1}$ is generated by $(A^{-} \oplus H) \cap GL_{2}(A)$. The spin group
$spin(A^{-} \oplus H)$ is the double cover
$\widetilde{Sp}_{A}^{\mathbb{F}^{2}}\binom{1\ \ \ \ 0}{0\ \ -1}$ of
$Sp_{A}^{\mathbb{F}^{2}}\binom{1\ \ \ \ 0}{0\ \ -1}$ defined by those pairs in
$\widetilde{GSp}_{A}^{\mathbb{F}^{2}}\binom{1\ \ \ \ 0}{0\ \ -1}$ having
multiplier 1. \label{dim8id1}
\end{thm}

\begin{proof}
Lemma \ref{ref8id1} and Proposition \ref{CDT} show, as in all the previous
cases, that the map from the semi-direct product of Lemma \ref{ac8id1} to
$O(A^{-} \oplus H)$ is surjective. Moreover, $\hat{\psi}$ has determinant $-1$
as an element of $O(A^{-} \oplus H)$ (it inverts a 3-dimensional subspace and
leaves the elements of its orthogonal complement invariant), so that
$\widetilde{GSp}_{A}^{\mathbb{F}^{2}}\binom{1\ \ \ \ 0}{0\ \ -1}$ maps to
$SO(A^{-} \oplus H)$ (Lemma \ref{ref8id1} again), with kernel
$\mathbb{F}^{\times}$, and the map is again surjective (by index
considerations). This shows that $Gspin(A^{-} \oplus
H)=\widetilde{GSp}_{A}^{\mathbb{F}^{2}}\binom{1\ \ \ \ 0}{0\ \ -1}$, and the
structure of the semi-direct product shows (using Lemma \ref{ref8id1} and
Proposition \ref{CDT} again) that this group is generated by those elements of
$\widetilde{GSp}_{A}^{\mathbb{F}^{2}}\binom{1\ \ \ \ 0}{0\ \ -1}$ lying over
$(A^{-} \oplus H) \cap GL_{2}(A)$ whose images under $\psi$ and $-\hat{\psi}$
coincide. The generation of $GSp_{A}^{\mathbb{F}^{2}}\binom{1\ \ \ \ 0}{0\ \
-1}$ follows, as the projection from the double cover is surjective. As the
element $\hat{\psi}$ of $O(A^{-} \oplus H)$ inverts a subspace of determinant 1,
it has spinor norm 1. Hence Lemma \ref{ref8id1} implies that this lift of
invertible $g \in A^{-} \oplus H$ has spinor norm $|g|^{2}$, which coincides
with the multiplier of this element. As these were seen to be a generating set
for $\widetilde{GSp}_{A}^{\mathbb{F}^{2}}\binom{1\ \ \ \ 0}{0\ \ -1}$, the
spinor norm of any element of the latter group is its multiplier (modulo
squares). The fact that this map factors through the projection to
$GSp_{A}^{\mathbb{F}^{2}}\binom{1\ \ \ \ 0}{0\ \ -1}$ is related to the space
$A^{-} \oplus H$ having discriminant 1, so that multiplying by $-Id_{A^{-}
\oplus H}$ does not affect the spinor norms. Therefore $SO^{1}(A^{-} \oplus H)$
consists of the images of those elements whose norm is a square, and
multiplication by suitable elements from the kernel $\mathbb{F}^{\times}$, we
may restrict to elements of multiplier 1. These are the elements of the double
cover $\widetilde{Sp}_{A}^{\mathbb{F}^{2}}\binom{1\ \ \ \ 0}{0\ \ -1}$ of the
symplectic group $Sp_{A}^{\mathbb{F}^{2}}\binom{1\ \ \ \ 0}{0\ \ -1}$ which is
defined by the multiplier 1 condition. As the only scalars with multiplier 1
are $\pm1$, this is the kernel of the (surjective) map from
$\widetilde{Sp}_{A}^{\mathbb{F}^{2}}\binom{1\ \ \ \ 0}{0\ \ -1}$ onto
$SO^{1}(A^{-} \oplus H)$, whence the former group is indeed $spin(A^{-} \oplus
H)$. This proves the theorem.
\end{proof}

\smallskip

Our space $A^{-} \oplus H$ is already assumed to be isotropic. However, we may
consider what happens when it splits more than one hyperbolic plane.
\begin{cor}
In case $A^{-} \oplus H$ splits more than one hyperbolic plane, there is a
quaternion algebra $B$ over $\mathbb{F}$ such that the Gspin and spin groups are
isomorphic to double covers of the subgroups of $GSp_{4}(B)$ and $Sp_{4}(B)$
whose presentation as in Corollary \ref{genform} (but with $M_{2}(B)^{-}$
replaced with $M_{2}^{Her}(B)$ and the lower right entry of the diagonal
generators being $m\iota_{B}(a)^{-t}$ rather than $m\overline{a}^{\ -1}$) uses
parameters from $GL_{2}^{(\mathbb{F}^{\times})^{2}}(B) \subseteq GL_{2}(B)$. If
it splits more than two hyperbolic planes, then our space is the direct sum of 4
hyperbolic planes, and our description of the spin group presents it as a double
cover of the group $SO^{1}\binom{0\ \ I}{I\ \ 0}$ of the direct sum of 4
hyperbolic planes in two ways, which are inequivalent to one another or to the
natural presentation as a double cover of such a group. The Gspin group is a
double cover of a subgroup of the general special orthogonal group of the direct
sum of 4 hyperbolic planes, again in two inequivalent ways. \label{iso8id1}
\end{cor}

\begin{proof}
Conjugation by an arbitrary element $\binom{e\ \ f}{g\ \ h} \in GL_{2}(A)$ takes
the group $GSp_{A}\binom{1\ \ \ \ 0}{0\ \ -1}$ to the $GSp$ group of the matrix
$\binom{e\ \ f}{g\ \ h}\binom{1\ \ \ \ 0}{0\ \ -1}\binom{\ \ \overline{h}\ \
-\overline{f}}{-\overline{g}\ \ \ \ \overline{e}}$, with the same multipliers.
$\varphi$ may be defined through the reverse conjugation, and conjugating any
$\psi$-image by the same matrix yields an order 2 automorphism of a double cover
of the kernel of this $\varphi$ (we may also transfer the action on $A^{-}
\oplus H$ by conjugating its image as well). When we do this with $\binom{e\ \
f}{g\ \ h}=\binom{1\ \ 0}{0\ \ R}$ for $R \in A^{-} \cap A^{\times}$ (and
multiplying by the global scalar $-1$) we get $GSp_{A}\binom{R\ \ 0}{0\ \ R}$.
In this case Corollary \ref{NAvn2} shows that the $GSp_{A}^{\mathbb{F}^{2}}$ (or
$\ker\varphi$) condition keeps its shape, since this conjugation just multiplies
the entries by $R$ or $R^{-1}$. When the $A^{-}$ part of $A^{-} \oplus H$ is
also isotropic, Corollary \ref{iso6d1} shows that we can take $A=M_{2}(B)$, and
we choose the element $R=\binom{0\ \ -1}{1\ \ \ \ 0}$. An application of Lemma
\ref{Sadjt} once on matrices in $M_{2}(A)$ and another time on the entries from
$A=M_{2}(B)$ shows that $GSp_{M_{2}(B)}\binom{R\ \ 0}{0\ \ R}$ is exactly the
group of matrices $g \in M_{4}(B)$ such that $g\binom{0\ \ -I}{I\ \ \ \
0}\iota_{B}(g)^{t}$ is some multiple of $\binom{0\ \ -I}{I\ \ \ \ 0}$. The group
$GSp_{M_{2}(B)}\binom{R\ \ 0}{0\ \ R}$ is thus $GSp_{4}(B)$, and the same for
the symplectic groups: $Sp_{M_{2}(B)}\binom{R\ \ 0}{0\ \ R}=Sp_{4}(B)$. Applying
this to the generators appearing in Corollary \ref{genform} indeed yields the
generators appearing in the parentheses, and as the upper left entry is not
affected, we get the asserted description of the Gspin and spin groups.

We can also conjugate with $\binom{1\ \ 0}{0\ \ S}$ with $S \in A^{+}$, yielding
$GSp_{A}\binom{S\ \ \ \ 0}{0\ \ -S}$. If $A^{-}$ splits more than one hyperbolic
plane then $A^{-} \oplus H$ is the sum of 4 hyperbolic planes (see Corollary
\ref{iso6d1} again), $B=M_{2}(\mathbb{F})$, and $A=M_{4}(\mathbb{F})$. We choose
$S$ to be the tensor product of $\binom{0\ \ -1}{1\ \ \ \ 0}$ with itself, i.e.,
the matrix $\binom{\alpha\ \ \beta}{\gamma\ \ \delta} \in M_{4}(\mathbb{F})$ in
which $\alpha=\delta=0$ and $\gamma=-\beta=\binom{0\ \ -1}{1\ \ \ \ 0} \in
M_{2}(\mathbb{F})$. By the form of the conjugating matrix, the definition of
$\ker\varphi$ remains the same also in this case. After applying Lemma
\ref{Sadjt} twice as in the previous case, plus another time on the entries of
$B=M_{2}(\mathbb{F})$, the group $\widehat{GSp}_{A}\binom{S\ \ \ \ 0}{0\ \ -S}$
is seen to be the group of matrices $g \in M_{8}(\mathbb{F})$ such that
$g\binom{0\ \ I}{I\ \ 0}g^{t}$ is a multiple of $\binom{0\ \ I}{I\ \ 0}$ (i.e.,
the \emph{general orthogonal group} of that matrix). $GSp_{A}\binom{S\ \ \ \
0}{0\ \ -S}$ is then the \emph{general special orthogonal group}, in which the
determinant is the 4th power of the multiplier, and $Sp_{A}\binom{S\ \ \ \ 0}{0\
\ -S}$ is just $SO\binom{0\ \ I}{I\ \ 0}$.

We claim that the map $\varphi$ on $Sp_{M_{4}(\mathbb{F})}\binom{1\ \ \ \ 0}{0\
\ -1}$ corresponds to the spinor norm in this presentation as a special
orthogonal group. It suffices to verify this again on a set of generators, and
we use those from the proof of Proposition \ref{GSppresHA-} once more. Recall
that we must take the conjugates of our elements by $\binom{1\ \ 0}{0\ \ S}$,
and that $S \in M_{4}(\mathbb{F})^{+}$ and satisfies $S^{2}=I$. This conjugation
replaces $\binom{1\ \ v}{0\ \ 1}$ and $\binom{1\ \ 0}{w\ \ 1}$ by $\binom{1\ \
vS}{0\ \ \ 1\ }$ and $\binom{\ 1\ \ \ 0}{Sw\ \ 1}$, with $vS$ being in
$M_{4}^{as}(\mathbb{F})$ by two applications of Lemma \ref{Sadjt} (as in
Corollary \ref{alt6d1}), and $Sw=Sw\overline{S} \cdot S$ lies there as well. For
$\binom{a\ \ \ 0\ \ }{0\ \ \overline{a}^{-1}}$ (we restrict to elements of
multiplier 1, since we do not consider ``spinor norms for general orthogonal
groups'' here), the choice of $S$ and Lemma \ref{Sadjt} imply that after
conjugation, $\overline{a}^{\ -1}$ is replaced by $a^{-t}$. Now, the unipotent
generators (which lie in $\ker\varphi$) are squares in $Sp_{A}\binom{1\ \ \ \
0}{0\ \ -1}$ (a square root is obtained by dividing the entry $vS$ or $Sw$ from
$M_{4}^{as}(\mathbb{F})$ by 2). Hence they have trivial spinor norms since the
range $\mathbb{F}^{\times}/(\mathbb{F}^{\times})^{2}$ of the spinor norm has
exponent 2. For $\binom{a\ \ \ 0\ \ }{0\ \ a^{-t}}$ with $a \in
GL_{4}(\mathbb{F})$, we recall that the latter group is generated by elementary
matrices, which operate only on two of these hyperbolic planes. It thus suffices
to consider the operation of $\binom{g\ \ \ 0\ \ }{0\ \ g^{-t}}$ with $g \in
GL_{2}(\mathbb{F})$ on the direct sum of two hyperbolic planes. But Corollaries
\ref{dim4d1} and \ref{iso4} show that the latter space is isometric to
$M_{2}(\mathbb{F})$ with the determinant as the vector norm, the Gspin group
being the ``equal determinant subgroup'' of the product of two copies of
$GL_{2}(\mathbb{F})$. By considering the first hyperbolic plane as generated by
$\binom{1\ \ 0}{0\ \ 0}$ and $\binom{0\ \ 0}{0\ \ 1}$ and the second one by
$\binom{0\ \ 0}{1\ \ 0}$ and $\binom{0\ \ -1}{0\ \ \ \ 0}$, the action of
$\binom{g\ \ \ 0\ \ }{0\ \ g^{-t}}$ becomes the action of the pair consisting of
$g$ and $\binom{1\ \ \ \ 0\ \ }{0\ \ \det g}$, and the spinor norm is indeed
$\det g$. Thus $Sp_{M_{4}(\mathbb{F})}^{\mathbb{F}^{2}}\binom{1\ \ \ \ 0}{0\ \
-1}$ is isomorphic to $SO^{1}\binom{0\ \ I}{I\ \ 0}$, and our description of the
spin group as $\widetilde{Sp}_{M_{4}(\mathbb{F})}^{\mathbb{F}^{2}}\binom{1\ \ \
\ 0}{0\ \ -1}$ is once again a presented as $spin\binom{0\ \ I}{I\ \ 0}$. We
thus have three representations of this group as a spin group of the direct sum
of 4 hyperbolic planes: The original one, the projection onto
$Sp_{M_{4}(\mathbb{F})}^{\mathbb{F}^{2}}\binom{1\ \ \ \ 0}{0\ \ -1}$, and the
composition of the latter projection with $\psi$. These representations are not
equivalent, since their kernels, all of order 2, are different: The non-trivial
element there is $-I$ with $\psi$-image $-I$ in the original representation, $I$
with $\psi$-image $-I$ in the projection, and $-I$ with $\psi$-image $I$ in the
composition. On the other hand,
$GSp_{M_{4}(\mathbb{F})}^{\mathbb{F}^{2}}\binom{1\ \ \ \ 0}{0\ \ -1}$ is a
subgroup of the general special orthogonal group of $\binom{0\ \ I}{I\ \ 0}$ 
which is defined by some condition which restricts to the triviality of the
spinor norm on $SO\binom{0\ \ I}{I\ \ 0}$, the Gspin group in question is a
double cover of this subgroup, and $\psi$ presents it as a double cover of this
subgroup in an inequivalent way (again, the projections have different kernels).
This completes the proof of the corollary. 
\end{proof}

Note that the groups from Corollary \ref{iso8id1} are not
$\widetilde{GSp}_{A}^{\mathbb{F}^{2}}\binom{1\ \ \ \ 0}{0\ \ -1}$, but
conjugates of the latter group inside $GL_{2}(A)$, with given conjugators. This
yields a definition for $\psi$ on these groups. However, we shall conjugate this
$\psi$ by $\binom{R\ \ 0}{0\ \ R}$ or by $\binom{S\ \ 0}{0\ \ S}$. The formula
for the resulting map looks just like that of $\psi$, but in which
$\overline{a}$ is replaced by $\iota_{B}(a)^{t}$ or just $a^{t}$, while
$\theta:v\mapsto\tilde{v}$ on $A^{-}$ becomes $X\mapsto-adjX$ on
$M_{2}^{Her}(B)$ or $T\mapsto\hat{T}$ on $M_{4}^{as}(\mathbb{F})$ (both having
the property that multiplication of the vector and its image under the
involution yields the vector norm). Hence when the groups from Corollary
\ref{iso8id1} are considered, this is the choice of $\psi$ with which they
come. 

The fact that in the hyperbolic case we get 3 inequivalent 8-dimensional
representations of the spin group, in all of which the image is the $SO^{1}$
group of the direct sum of 4 hyperbolic planes, is an incarnation of
\emph{triality} for this case. Triality exists for more general settings, namely
some non-isotropic spaces of dimension 8 and discriminant 1 (see Section 35 of
\cite{[KMRT]} for more details), but our methods here restrict to the isotropic
case.

We remark that allowing non-trivial spinor norms in the second case in Corollary
\ref{iso8id1} is in some sense dual to allowing multipliers. We have seen that
the Gspin group was mapping to the general special orthogonal group in this
case. On the other hand, we may allow the map $\psi$ from Theorem \ref{GSppsi}
to have a free choice of a scalar (not necessarily squaring to the reduced norm
of an entry), which would extend the definition of $\psi$ to (an
$\mathbb{F}^{\times}$-cover of) all of $GSp_{A}\binom{1\ \ \ \ 0}{0\ \ -1}$, not
only to elements with trivial $\varphi$-image. The group constructed in Theorem
\ref{GSppsi} would then be an $\mathbb{F}^{\times}$-cover of
$GSp_{A}\binom{1\ \ \ \ 0}{0\ \ -1}$, and the action from Proposition
\ref{GSppresHA-} and Lemma \ref{ac8id1} may multiply the bilinear form on $A^{-}
\oplus H$ by a scalar. In the split $A$ case this was seen, when we considered
the generators $\binom{a\ \ \ 0\ \ }{0\ \ \overline{a}^{-1}}$ with multiplier 1,
to produce elements whose image in the projection to $SO(M_{4}(\mathbb{F})^{-}
\oplus H)$, as well as in the composition of this projection with $\psi$, may
have arbitrary spinor norms (but the spinor norm does have to be the same for
these two maps).

In any case, we may have many alternative descriptions of this picture:
\begin{cor}
Let $\Xi$ be an element of $\widetilde{GSp}_{A}^{\mathbb{F}^{2}}\binom{1\ \ \ \
0}{0\ \ -1}$, such that $\Xi\psi(\Xi)$ is a scalar $r\in\mathbb{F}^{\times}$.
Define the map $\psi_{\Xi}:\widetilde{GSp}_{A}^{\mathbb{F}^{2}}\binom{1\ \ \ \
0}{0\ \ -1}\to\widetilde{GSp}_{A}^{\mathbb{F}^{2}}\binom{1\ \ \ \ 0}{0\ \ -1}$
by conjugating $\psi$ by $\Xi$, i.e., $\psi_{\Xi}(g)=\Xi\psi(g)\Xi^{-1}$. Let
$\hat{\psi}_{\Xi}$ be the composition of $\hat{\psi}$ with the operation of
$\Xi$ on $A^{-} \oplus H$, and we embed the latter space into $M_{2}(A)$ by
multiplying the image from above by $\Xi^{-1}$ from the right. Then all the
assertions from Lemma \ref{vnorm8} to Theorem \ref{dim8id1} and Corollary
\ref{iso8id1} hold by replacing every $U$ by $U\binom{0\ \ -1}{1\ \ \ \ 0}$,
$\psi$ by $\psi_{\Xi}$, and $\hat{\psi}$ by $\hat{\psi}_{\Xi}$, up to rescaling
the bilinear forms by the scalar $r$. \label{alt8id1}
\end{cor}

\begin{proof}
The assumption that $\Xi\psi(\Xi)$ is central implies that $\psi_{\Xi}$ again
has order 2 as an automorphism of $\widetilde{GSp}_{A}^{\mathbb{F}^{2}}\binom{1\
\ \ \ 0}{0\ \ -1}$. As an element $V$ of the latter space is $U\Xi^{-1}$ with $U
\in A^{-} \oplus H$ as above, its $\psi_{\Xi}$-image equals
$\Xi\psi(U)\psi(\Xi)^{-1}\Xi^{-1}$. This coincides with our definition of
$\hat{\psi}_{\Xi}$ on $U$ and the modified embedding, so that $\hat{\psi}_{\Xi}$
preserves this embedding of $A^{-} \oplus H$ and is the restriction of a branch
of $\psi_{\Xi}$. In addition, our assumption on $\Xi$ implies that the latter
expression is just $\frac{\Xi\hat{\psi}(U)}{r}$. The original Lemma \ref{vnorm8}
now yields its modified version, with the appropriate rescaling . Furthermore,
multiplying our $V$ by $\psi_{\Xi}(g)^{-1}=\Xi\psi(g)^{-1}\Xi^{-1}$ from the
right gives $U\psi(g)^{-1}\Xi^{-1}$. After left multiplication by $g$, the
original Proposition \ref{GSppresHA-} implies the modified one. All the rest now
follows from these assertions in the same way. This proves the corollary.
\end{proof}
For example, if we take $\Xi=\binom{0\ \ -1}{1\ \ \ \ 0}$ (with $r=1$) in
Corollary \ref{alt8id1} then $A^{-} \oplus H$ becomes the space of matrices of
the form $\binom{p\ \ u}{\tilde{u}\ \ q}$, with the usual $p$, $q$, and $u$ and
with minus the ``bi-quaternionic Moore determinant'' as the vector norm. The map
$\hat{\psi}_{\Xi}$ interchanges $p$ and $q$ with minus one another and leaves
$A^{-}$ pointwise fixed. In general, all the (equivalent) representations we get
in Corollary \ref{alt8id1} are still based on the map $\psi$, which is more
complicated, hence we shall not present any of them explicitly.

\section{Dimension 7, Representing the Discriminant \label{Dim7rd}}

The spaces we consider here are given in
\begin{lem}
The orthogonal complement of a vector in $A^{-} \oplus H$ of some vector norm
$-\delta\neq0$ has discriminant $\delta$ and it contains a vector of norm
$\delta$. Any vector space of dimension 7 containing a vector whose vector norm
equals the discriminant of the space can be obtained in this manner, up to
rescaling . \label{sp7rd}
\end{lem}
We remark that if some vector $Q$ has vector norm which equals the discriminant
of the space then it continues to hold after rescalings.

\begin{proof}
The discriminant and determinant of $A^{-} \oplus H$ is 1. As a vector $Q$ of
vector norm $-\delta$ spans a space of determinant $-\delta$, the complement in
$A^{-} \oplus H$ has the same determinant $-\delta$, and its discriminant is
$\delta$ since $(-1)^{7(7-1)/2}=-1$. As a hyperbolic plane contains vectors of
any given vector norm, the Witt Cancelation Theorem allows us to find some
element of $O(A^{-} \oplus H)$ taking $Q$ to some element of $H \subseteq A^{-}
\oplus H$. The orthogonal complement in $H$ is generated by a vector of vector
norm $\delta$, and its inverse image under the orthogonal map we applied has the
same vector norm. Conversely, adding some vector $Q$ which is perpendicular to
the total space and such that $|Q|^{2}$ is the determinant of the space yields a
space of discriminant 1. The sum of $Q$ with a vector whose vector norm is the
discriminant of the space is then isotropic. Lemma \ref{sp8id1} now completes
the proof of the lemma.
\end{proof}

In view of Lemma \ref{sp7rd}, we write our space as
$A^{-}\oplus\langle\delta\rangle$, using a generator of the orthogonal
complement in $H$ whose norm is $\delta$. Moreover, we embed the space $A^{-}
\oplus H$ into $M_{2}(A)$ as seen after Lemma \ref{sp8id1}, and we choose $Q$ to
be the matrix $\binom{0\ \ -\delta}{1\ \ \ \ 0}$ (of vector norm $-\delta$). We
now have
\begin{lem}
Given such $A$ with $\iota_{B}\otimes\iota_{C}$ and $\delta$, the groups
$Gspin(A^{-}\oplus\langle\delta\rangle)$ and
$spin(A^{-}\oplus\langle\delta\rangle)$ are, up to isomorphism, the stabilizers
of the matrix $\binom{0\ \ \delta}{1\ \ 0}$ in the action given in Proposition
\ref{GSppresHA-} inside the groups
$\widetilde{GSp}_{A}^{\mathbb{F}^{2}}\binom{1\ \ \ 0}{0\ \ -1}$ and
$\widetilde{Sp}_{A}^{\mathbb{F}^{2}}\binom{1\ \ \ 0}{0\ \ -1}$ respectively. The
double cover $\widetilde{GSp}_{A}^{\mathbb{F}^{2}}\binom{1\ \ \ 0}{0\ \ -1}$
splits over $Gspin(A^{-}\oplus\langle\delta\rangle)$. These groups operate by
conjugation on $A^{-}\oplus\langle\delta\rangle$ if we identify the latter space
as the space of matrices of the form $\binom{p\ \ \ \delta u}{\tilde{u}\ \ -p}$
with $u \in A^{-}$ and $p\in\mathbb{F}$, with the vector norm being the
``bi-quaternionic $A^{-}$-Moore determinant'' divided by $-\delta$.
\label{ac7rd}
\end{lem}

\begin{proof}
The first assertion is proved by the same argument used for Lemma \ref{ac5}.
Now, that the action from Proposition \ref{GSppresHA-} is based on the map
$\psi$, hence an element of $GSp_{A}^{\mathbb{F}^{2}}\binom{1\ \ \ 0}{0\ \ -1}$
would stabilize $\binom{0\ \ -\delta}{1\ \ \ \ 0}$ if and only if the
$\psi$-image of one of its lifts equals its conjugate by $\binom{0\ \
-\delta}{1\ \ \ \ 0}$. This yields the splitting of the double cover over
$Gspin(A^{-}\oplus\langle\delta\rangle)$, since every element there comes with a
natural choice of $\psi$-image. Take now $\Xi=\binom{0\ \ -\delta}{1\ \ \ \ 0}$
in Corollary \ref{alt8id1}, which equals its $\psi$-image and squares to
$-\delta$. The space thus obtained is the one written explicitly here, and the
remaining assertions follow since $\psi_{\Xi}(g)=g$ for $g \in
Gspin(A^{-}\oplus\langle\delta\rangle)$ by definition. This proves the lemma.
\end{proof}

In order to give a more detailed description of the groups from Lemma
\ref{ac7rd}, we begin by proving
\begin{lem}
If an element of $\widetilde{GSp}_{A}^{\mathbb{F}^{2}}\binom{1\ \ \ 0}{0\ \ -1}$
lying over $\binom{e\ \ f}{g\ \ h} \in GSp_{A}^{\mathbb{F}^{2}}\binom{1\ \ \
0}{0\ \ -1}$ stabilizes $\binom{0\ \ -\delta}{1\ \ \ \ 0}$ then either $e$ or
$g$ are invertible. \label{stabinv}
\end{lem}

\begin{proof}
First observe that is $(a,t)\in\widetilde{A}^{(\mathbb{F}^{\times})^{2}}$ then
the action of the element $\binom{a\ \ \ \ 0\ \ \ }{0\ \
t\overline{a}^{-1}}$ of $\widetilde{GSp}_{A}^{\mathbb{F}^{2}}\binom{1\ \ \ 0}{0\
\ -1}$, with $\psi$-image $\binom{t\overline{a}^{-1}\ \ 0\ \ }{\ \ 0\ \ ma/t}$,
stabilizes this matrix. The proof of Lemma \ref{genie} thus shows that it
suffices to consider elements in which $e=\binom{1\ \ 0}{0\ \ 0} \in A=M_{2}(B)$
for some quaternion algebra $B$ over $\mathbb{F}$ (for $A$ division the lemma is
immediate). We have seen in the proof of Lemma \ref{genie} that $g$ takes the
form $\binom{\lambda\ \ \mu}{r\ \ 0}$ with $\lambda$ and $\mu$ from $B$ and
$r\in\mathbb{F}$, and a similar argument shows that $f=\binom{\sigma\ \ s}{\tau\
\ 0}$ where $\sigma$ and $\tau$ are in $B$ and $s\in\mathbb{F}$. Moreover,
$\mu\overline{\tau}$ equals the multiplier $m$ (hence both $\mu$ and $\tau$ are
invertible), and $h$ has lower right entry $m+rs$. As parameters for $\binom{e\
\ f}{g\ \ h}$ in Corollary \ref{genform} may be taken to be $v=\binom{\ \ 0\ \
0}{-1\ \ 0}$, $a=\binom{1\ \ 0}{\lambda\ \ \mu}$, $\alpha$ with upper right
entry $s\in\mathbb{F}$, and $\beta=\binom{0\ \ 1}{r\ \ 0}$. The
$GSp_{A}^{\mathbb{F}^{2}}$ condition means that
$N^{B}_{\mathbb{F}}(a)=N^{B}_{\mathbb{F}}(\mu)$ is a square, say $t^{2}$, hence
$t\overline{a}^{\ -1}=\binom{\mu/t\ \ \ \ 0}{\overline{\lambda}\mu/t\ \ t}$, and
the action of $\theta$ leaves $u$, $\beta$, and the $s$-entry of $\alpha$
invariant. The lowest row of $\binom{1\ \ v}{0\ \ 1}\binom{1\ \ 0}{\beta\ \ 1}$
is $(r\ \ 0\ \ 0\ \ 1)$, while in $t\overline{a}^{\ -1}\binom{1\ \
\tilde{\alpha}}{0\ \ 1}$ the most upper right and lower right entries are
$\frac{s\mu}{t}$ and $\frac{m\mu}{t}$ respectively. Hence the most lower right
entry of $\psi\binom{e\ \ f}{g\ \ h}$ is $\frac{(m+rs)\mu}{t}$. As
$\psi\binom{e\ \ f}{g\ \ h}$ has multiplier $m$, its inverse has
$\frac{(m+rs)\mu}{mt}$ as its most upper left entry. Now, the second row of
$\binom{e\ \ f}{g\ \ h}\binom{0\ \ -\delta}{1\ \ \ \ 0}$ is $(\tau\ \ 0\ \ 0\ \
0)$, so that the second row of the action of $\binom{e\ \ f}{g\ \ h}$ on
$\binom{0\ \ -\delta}{1\ \ \ \ 0}$ starts with $\frac{(m+rs)\tau\mu}{mt}$. But
we have assumed that $\binom{e\ \ f}{g\ \ h}$ preserves $\binom{0\ \ -\delta}{1\
\ \ \ 0}$, so that the latter expression must vanish. As $\tau$ and $\mu$ are in
$B^{\times}$ and $m\neq0$, it follows that $r\neq0$, whence $g=\binom{\lambda\ \
\mu}{r\ \ 0}$ is invertible as desired. This completes the proof of the lemma.
\end{proof}

The determination of the group from Lemma \ref{ac7rd} may now be carried out
using the explicit formulae for $\psi$. The result is
\begin{prop}
Any element of $Gspin(A^{-}\oplus\langle\delta\rangle)$ may be presented either
as $\binom{\ a\ \ -t\delta\tilde{\beta}\overline{a}^{-1}}{\beta a\ \ \ \
t\overline{a}^{-1}\ \ }$ with
$(a,t)\in\widetilde{A}^{(\mathbb{F}^{\times})^{2}}$ and $\beta \in A^{-}$, or as
$\binom{wc\ \ -s\delta\overline{c}^{-1}}{\ c\ \ \ \
s\tilde{w}\overline{c}^{-1}}$ where $(c,s)$ is in
$\widetilde{A}^{(\mathbb{F}^{\times})^{2}}$ and $w$ comes from $A^{-}$.
\label{stabdelta}
\end{prop}

\begin{proof}
Lemma \ref{stabinv} implies that for every element $\binom{a\ \ b}{c\ \ d} \in
GSp_{A}^{\mathbb{F}^{2}}\binom{1\ \ \ 0}{0\ \ -1}$ which stabilizes $\binom{0\ \
-\delta}{1\ \ \ \ 0}$, either $\binom{a\ \ b}{c\ \ d}$ or $\binom{0\ \
-\delta}{1\ \ \ \ 0}\binom{a\ \ b}{c\ \ d}$ has an invertible upper right entry.
In the first case $a$ lies under some element
$(a,t)\in\widetilde{A}^{(\mathbb{F}^{\times})^{2}}$, and Lemma \ref{invent}
shows that $b=\beta a$ for some $\beta \in A^{-}$. Moreover, the formula from
Theorem \ref{GSppsi} shows that $\psi\binom{a\ \ b}{c\ \ d}$ has left entry
$\binom{t\overline{a}^{-1}}{t\tilde{\beta}\overline{a}^{-1}}$, and if
$\psi\binom{a\ \ b}{c\ \ d}$ coincides with the image of $\binom{a\ \ b}{c\ \
d}$ under conjugation by $\binom{0\ \ -\delta}{1\ \ \ \ 0}$, then $\binom{a\ \
b}{c\ \ d}$ must have the asserted right colmun. If $a$ is not invertible, then
we may multiply $\binom{a\ \ b}{c\ \ d}$ by $\binom{0\ \ -\delta}{1\ \ \ \ 0}$
from the left, obtain an element of the form just described, and dividing by
$\binom{0\ \ -\delta}{1\ \ \ \ 0}$ back again shows that our element must be of
the second suggested form. This proves the proposition.
\end{proof}
Note that an element of the second form may be uniquely presented as the product
of an anisotropic element of $A^{-}\oplus\langle\delta\rangle=\binom{0\ \
-\delta}{1\ \ \ \ 0}^{\perp} \subseteq A^{-} \oplus H$ in which the lower left
entry is 1 and a diagonal matrix stabilizing $\binom{0\ \ \delta}{1\ \ 0}$ (the
latter multipliers form, as Proposition \ref{stabdelta} shows, a group which is
isomorphic to $\widetilde{A}^{(\mathbb{F}^{\times})^{2}}$): Indeed, such an
element is just the product $\binom{w\ \ \ \ \delta}{1\ \ -\tilde{w}}\binom{c\ \
\ \ \ 0\ \ \ }{0\ \ -r\overline{c}^{-1}}$.

In total, we have
\begin{thm}
The Gspin group of $A^{-}\oplus\langle\delta\rangle$ consists of elements
$\binom{a\ \ b}{c\ \ d}$ of $GSp_{A}^{\mathbb{F}^{2}}\binom{1\ \ \ 0}{0\ \ -1}$
in which $a\overline{d}$ and $b\overline{c}$ are scalars from $\mathbb{F}$,
which square to $N^{A}_{\mathbb{F}}(a)$ (or equivalently
$N^{A}_{\mathbb{F}}(d)$) and $\delta^{2}N^{A}_{\mathbb{F}}(c)$ (which equals
also $\frac{N^{A}_{\mathbb{F}}(b)}{\delta^{2}}$) respectively. It is
characterized by either the two elements $bd^{-1}$ and $-\delta ca^{-1}$ or the
two elements $ac^{-1}$ and $-\delta db^{-1}$ being well-defined elements of
$A^{-}$ which are $\theta$-images of one another. It is generated by anisotropic
vectors of the form $\binom{v\ \ \ \ \delta}{1\ \ -\tilde{v}}$ or $\binom{v\ \ \
\ 0}{0\ \ -\tilde{v}}$. The spin group consists of those elements in which the
two scalars $a\overline{d}$ and $b\overline{c}$ sum to 1. \label{dim7rd}
\end{thm}

\begin{proof}
Any element may be presented in one of the two forms given in Proposition
\ref{stabdelta}. In the first case we have $a\overline{d}=t$ and
$b\overline{c}=t\delta|\beta|^{2}$, while in the second one these numbers are
$-s|w|^{2}$ and $-\delta s$ respectively. The relations with the reduced norms
of $a$, $b$, $c$, and $d$ are easily verified using Corollary \ref{NAvn2} and
the definition of $\widetilde{A}^{(\mathbb{F}^{\times})^{2}}$, and the relation
between $bd^{-1}$ and $-\delta ca^{-1}$ in the first case and $ac^{-1}$ and
$-\delta db^{-1}$ in the third case are also immediate. Conversely, if
$\binom{a\ \ b}{c\ \ d}$ is a matrix in which $a\overline{d}$ and
$b\overline{c}$ are scalars, not both zero, then either $d=t\overline{a}^{\ -1}$
or $b=-s\delta\overline{c}^{\ -1}$ for some scalars $t$ and $s$. If
$t^{2}=N^{A}_{\mathbb{F}}(a)$ then $N^{A}_{\mathbb{F}}(d)$ takes the same value,
while if $s^{2}=N^{A}_{\mathbb{F}}(c)$ then $N^{A}_{\mathbb{F}}(b)$ is obtained
by multiplication by $\delta^{4}$. If we write $c=\beta a$ in the first case and
$a=wc$ in the second case, then the respective values
$b=-t\delta\tilde{\beta}\overline{a}^{\ -1}$ and $d=s\tilde{w}\overline{c}^{\
-1}$ immediately follow. Now, elements with invertible lower left entry were
seen to take the form $\binom{w\ \ \ \ \delta}{1\ \ -\tilde{w}}\binom{c\ \ \ \ \
0\ \ \ }{0\ \ -s\overline{c}^{-1}}$. The generation of these elements by the
asserted set now follows from Theorem \ref{dim6d1}, as the map
$(a,t)\mapsto\binom{a\ \ \ \ 0\ \ }{0\ \ t\overline{a}^{-1}}$ is a group
injection which sends the generators $(v,|v|^{2})$ for $v \in A^{-} \cap
A^{\times}$ to $\binom{v\ \ \ \ 0}{0\ \ -\tilde{v}}$. The other elements are
obtained by multiplying the appropriate elements by the generator $\binom{0\ \
\delta}{1\ \ 0}$, and the assertion about generation follows. For the spin
group, observe that the proof of Lemma \ref{ac5} shows that the spinor norm of
an element of $SO(A^{-}\oplus\langle\delta\rangle)$ is the same when considered
there or in $SO(A^{-} \oplus H)$ (by leaving $\binom{0\ \ -\delta}{1\ \ \ \ 0}$
invariant), and the proof of Theorem \ref{dim8id1} implies that in the latter
group the spinor norms of (the image of) $\binom{a\ \ b}{c\ \ d}$ is just the
multiplier $a\overline{d}+b\overline{c}$. As by the usual scalar multiplication
we may normalize this multiplier to 1 wherever it is a square, the assertion
about the spin group is also established. This proves the theorem.
\end{proof}

After fixing $A$, we have the following assertion about the dependence of the
Gspin and spin groups on $\delta$: 
\begin{prop}
If $\varepsilon\in\delta(\mathbb{F}^{\times})^{2}N^{A}_{\mathbb{F}}(A^{\times})$
then the spin and Gspin groups of $A^{-}\oplus\langle\varepsilon\rangle$ are
isomorphic to those of $A^{-}\oplus\langle\delta\rangle$. \label{deltadep} 
\end{prop}

\begin{proof}
Consider first multiplication from $(\mathbb{F}^{\times})^{2}$. Let
$r\in\mathbb{F}^{\times}$, and we examine the result of conjugation by
$\binom{1\ \ 0}{0\ \ 1/r}$. This operation multiplies the upper right entry by
$r$ and divides the lower left entry by $r$. Hence on elements of
$Gspin(A^{-}\oplus\langle\delta\rangle)$ of the first form of Proposition
\ref{stabdelta} this operation corresponds to leaving $a$ (and $t$) invariant,
dividing $\beta$ by $r$, and multiplying $\delta$ by $r^{2}$, while for elements
of the second form it means dividing $c$ by $r$ (hence dividing $s$ by $r^{2}$),
multiplying $w$ by $r$, and again multiplying $\delta$ by $r^{2}$. Hence this
conjugation takes $Gspin(A^{-}\oplus\langle\delta\rangle)$ into
$Gspin(A^{-}\oplus\langle r^{2}\delta \rangle)$. Conjugation by the inverse
element shows that the map between these two groups is bijective. As for
multiplication by norms from $A^{\times}$, we now consider the conjugation by
$\binom{e\ \ \ 0\ \ }{0\ \ \overline{e}^{-1}}$ for some $e \in A^{\times}$. For
elements having the first form in Proposition \ref{stabdelta}, this operation
sends $a$ to $eae^{-1}$ (hence $t$ remains invariant) and $\beta$ to
$\overline{e}^{\ -1}\beta e^{-1}$, and Lemma \ref{AxA-rel} shows that $\delta$
must be multiplied by $N^{A}_{\mathbb{F}}(e)$. As for the other elements, this
operation takes $c$ to $\overline{e}^{\ -1}ce^{-1}$ (and therefore $s$ is
divided by $N^{A}_{\mathbb{F}}(e)$) and $w$ to $ew\overline{e}$, so that again
$\delta$ has to be multiplied by $N^{A}_{\mathbb{F}}(e)$ (Lemma \ref{AxA-rel} is
just a consistency check in this case). This
$Gspin(A^{-}\oplus\langle\delta\rangle)$ is isomorphic to a subgroup of
$Gspin\big(A^{-}\oplus\big\langle N^{A}_{\mathbb{F}}(e)\delta \big\rangle\big)$,
and an argument using the inverse conjugation show the bijectivity. As
conjugation preserves multipliers and the spin groups are the subgroups of the
Gspin groups which are defined by the multiplier 1 condition, this completes the
proof of the proposition. 
\end{proof}
By Proposition \ref{deltadep}, it suffices to take $\delta$ from a set of
representatives for
$\mathbb{F}^{\times}/(\mathbb{F}^{\times})^{2}N^{A}_{\mathbb{F}}(A^{\times})$.
In addition, Lemma \ref{NM2B} once again shows that if $A=M_{2}(B)$ then this
group involves just classes modulo $N^{B}_{\mathbb{F}}(B^{\times})$. 

The concept of isotropy in this case is the one considered in
\begin{cor}
Assume that $A^{-}\oplus\langle\delta\rangle$ contains an isotropic vector which
is orthogonal to a vector of vector norm which equals the discriminant. Then the
Gspin and spin group consist of matrices $\binom{a\ \ b}{c\ \ d}$ from
$GSp_{4}(B)$ (or $Sp_{4}(B)$), for some quaternion algebra $B$ over
$\mathbb{F}$, in which $a\iota_{B}(d)^{t}$ and $b\iota_{B}(c)^{t}$ lie in
$\mathbb{F}$, and square to
$N^{M_{2}(B)}_{\mathbb{F}}(a)=N^{M_{2}(B)}_{\mathbb{F}}(d)$ and
$\delta^{2}N^{M_{2}(B)}_{\mathbb{F}}(c)=\frac{N^{M_{2}(B)}_{\mathbb{F}}(b)}{
\delta^{2}}$ respectively. In every such matrix, either $bd^{-1}$ and $-\delta
ca^{-1}$ or $ac^{-1}$ and $-\delta db^{-1}$ lie in $M_{2}^{Her}(B)$ and are
minus the adjoints of one another. These groups operate by conjugation on the
space of matrices $\binom{\ pI\ \ \ -\delta X}{adjX\ \ -pI\ } \in M_{4}(B)$,
where $X \in M_{2}^{Her}(B)$ and $p\in\mathbb{F}$, as the Gspin and spin groups
of this space with the with the ``bi-quaternionic $A^{-}$-Moore determinant''
divided by $-\delta$ as the vector norm. The spin group may also be presented as
the spin group of the space of matrices of the sort $\binom{\delta X\ \ \ \ pI\
}{pI\ \ \ adjX}$, with the same ``bi-quaternionic Moore determinant'' as the
vector norm, via $g:N \mapsto gN\iota_{B}(g)^{t}$.

In case $A^{-}\oplus\langle\delta\rangle$ splits three hyperbolic planes, we
get groups of $8\times8$ matrices $\binom{a\ \ b}{c\ \ d}$ which multiply the
bilinear form defined by $\binom{0\ \ I}{I\ \ 0}$ by a scalar, such that
$ad^{t}$ and $bc^{t}$ are scalar $4\times4$ matrices, whose squares are $\det
a=\det d$ and $\delta^{2}\det c=\frac{\det b}{\delta^{2}}$ respectively.
Moreover, either the pair $bd^{-1}$ and $-\delta ca^{-1}$ or the pair $ac^{-1}$
and $-\delta db^{-1}$ are defined, they lie in $M_{4}^{as}(\mathbb{F})$, and
they are sent to one another by the involution $T\mapsto\hat{T}$ which was
described in the paragraph following Corollary \ref{alt6d1}. The space on which
these groups operate by conjugation consists of matrices of the form $\binom{pI\
\ -\delta T}{\hat{T}\ \ \ -pI}$ with $T \in M_{4}^{as}(\mathbb{F})$.
\label{iso7rd}
\end{cor}

\begin{proof}
We have seen in the proof of Corollary \ref{iso8id1} that $GSp_{4}(B)$ is
obtained from $GSp_{M_{2}(B)}\binom{1\ \ \ \ 0}{0\ \ -1}$ through conjugation by
$\binom{1\ \ 0}{0\ \ R}$ with $R=\binom{0\ \ -1}{1\ \ \ \ 0}$. As this operation
takes a matrix $\binom{a\ \ b}{c\ \ d}$ to $\binom{\ a\ \ \ bR^{-1}\ }{Rc\ \
RdR^{-1}}$, Lemma \ref{Sadjt} shows that the relations from Theorem \ref{dim7rd}
become the ones asserted here (the reduced norms are not affected, since
$N^{M_{2}(B)}_{\mathbb{F}}(R)=1$). The assertions involving $bd^{-1}$ and
$-\delta ca^{-1}$ or $ac^{-1}$ and $-\delta db^{-1}$ follow from those appearing
in Theorem \ref{dim7rd} through the fact that $\theta(X)=-adjX$ for $X \in
M_{2}(B)^{-}$, right multiplication by $R$ sends this space to $M_{2}^{Her}(B)$,
and $adjR=R^{-1}$ for our $R$. Recall now that the group
$Gspin(A^{-}\oplus\langle\delta\rangle) \subseteq GSp_{M_{2}(B)}\binom{1\ \ \ \
0}{0\ \ -1}$ operates on the space defined in Lemma \ref{ac7rd} by conjugation.
Conjugating the formula for this action by $\binom{1\ \ 0}{0\ \ R}$ yields the
action of our subgroup of $GSp_{4}(B)$ by conjugation on the asserted space with
the asserted quadratic form. The $Sp_{4}$ condition now shows that multiplying
the latter space from the right by $\binom{0\ \ -I}{I\ \ \ \ 0}$ yields a space
on which $spin\big(M_{2}^{Her}(B)\oplus\langle\delta\rangle\big) \subseteq
Sp_{4}(B)$ operates via $g:N \mapsto gN\iota_{B}(g)^{t}$, and this space is
easily seen to be the one from the last assertion,

In the case of splitting three hyperbolic planes, we apply the same argument
with the matrix $S$ from the proof of Corollary \ref{iso8id1}. Once again Lemma
\ref{Sadjt} yields the desired relations between the squares, and $\det S=1$.
We recall from the proof of Corollary \ref{alt6d1} that multiplication of
$M_{4}(\mathbb{F})^{-}$ by $S$ (from either side) yields anti-symmetric
matrices, and that the vector norm is taken to minus the pfaffian by this
operation. Conjugating the space from Lemma \ref{ac7rd} by $\binom{1\ \ 0}{0\ \
S}$ yields the first space with $T$ and $\hat{T}$. The fact that the spin group
is contained in $O\binom{0\ \ I}{I\ \ 0}$ allows us to multiply our
representation by $\binom{0\ \ I}{I\ \ 0}$ from the right, yielding the second
space with the action $g:L \mapsto gLg^{t}$ of the spin group. This completes
the proof of the corollary.
\end{proof}
Note that the description of the groups from Corollary \ref{iso7rd} is in
correspondence with the choice of $\psi$ on the groups from Corollary
\ref{iso8id1}. For the Gspin groups, the representations which extend those
defined by $g:N \mapsto gN\iota_{B}(g)^{t}$ and $g:L \mapsto gLg^{t}$ from
Corollary \ref{iso7rd} and preserve the bilinear form must include division by
the multiplier $m$. In addition, although in both cases we may obtain natural
8-dimensional representations of these groups by adding $\binom{0\ \ -I}{I\ \ \
\ 0}$ to the first representation and $\binom{0\ \ I}{I\ \ 0}$ to the second
one, this is not dual to preserving $\binom{0\ \ -\delta}{1\ \ \ \ 0}$ since the
full $GSp$ group (of $\binom{Q\ \ 0}{0\ \ Q}$ or $\binom{S\ \ \ \ 0}{0\ \ -S}$)
also preserve this matrix by definition. We also mention the fact that starting
with the representation appearing in Corollary \ref{alt8id1} yields precisely
the representations of $Gspin(A^{-}\oplus\langle\delta\rangle)$ and
$spin(A^{-}\oplus\langle\delta\rangle)$ already given in Lemma \ref{ac7rd} and
Corollary \ref{iso7rd}.

\section{Dimension 8, Isotropic, Any Discriminant \label{Dim8igen}}

Let $d$ be a discriminant, and let $\mathbb{E}=\mathbb{F}(\sqrt{d})$ be the
associated quadratic extension of $\mathbb{F}$, with Galois automorphism $\rho$.
We shall be interested in the spaces given in the following
\begin{lem}
Let $A=B \otimes C$ be a bi-quaternion algebra over $\mathbb{F}$ with the
involution corresponding to this presentation, and let $Q \in A^{-}$ be
anisotropic. The direct sum $(A_{\mathbb{E}}^{-})_{\rho,Q} \oplus H$ of a
hyperbolic plane and the space from Lemma \ref{sp6gen} is 8-dimensional,
isotropic, and has discriminant $d$. Moreover, this yields all the isotropic
8-dimensional quadratic spaces of discriminant $d$ over $\mathbb{F}$.
\label{sp8igen}
\end{lem}

\begin{proof}
The space $(A_{\mathbb{E}}^{-})_{\rho,Q}$ from Lemma \ref{sp6gen} has dimension
6 and discriminant $d$, hence determinant $-d$. Adding the isotropic space $H$,
of determinant $-1$ yields a space with the desired properties. On the other
hand, if an 8-dimensional space is isotropic and has discriminant $d$, then it
splits a hyperbolic plane, and the complement has dimension 6 and discriminant
$d$. The lemma now follows from Lemma \ref{sp6gen} and the fact that hyperbolic
planes are isometric to their rescalings.
\end{proof}

Extending scalars in the space from Lemma \ref{sp8igen} to $\mathbb{E}$, we
obtain an isotropic 8-dimensional space of discriminant 1, which equals
$A_{\mathbb{E}}^{-} \oplus H_{\mathbb{E}}$ by Lemma \ref{sp8id1} and we present
it as the subspace of $M_{2}(A_{\mathbb{E}})$ as before. Therefore our space
$(A_{\mathbb{E}}^{-})_{\rho,Q}$ is isomorphic to the space of matrices
$\binom{u\ \ -p}{q\ \ -\tilde{u}}$ in which $u\in(A_{\mathbb{E}}^{-})_{\rho,Q}$
and $p$ and $q$ are in $\mathbb{F}$, with the restriction of the quadratic form
from $A_{\mathbb{E}}^{-} \oplus H_{\mathbb{E}}$. Proposition \ref{GSppresHA-}
shows that the group $\widetilde{GSp}_{A_{\mathbb{E}}}^{\mathbb{E}^{2}}\binom{1\
\ \ \ 0}{0\ \ -1}$ acts on $A_{\mathbb{E}}^{-} \oplus H_{\mathbb{E}}$, and we
are interested in the subgroup which preserves the subspace
$(A_{\mathbb{E}}^{-})_{\rho,Q} \oplus H$. Observe that the $H$ part of
$(A_{\mathbb{E}}^{-})_{\rho,Q} \oplus H$ is invariant under the automorphism
$\rho$ of $M_{2}(A)$, and the action of $\rho$ on the other part is given in
Lemma \ref{sp6gen}. In particular $\rho$ preserves the subspace
$(A_{\mathbb{E}}^{-})_{\rho,Q} \oplus H$ of $A_{\mathbb{E}}^{-} \oplus
H_{\mathbb{E}}$, and operates as the reflection in $Q$ on this space.

We shall be needing also non-invertible elements of $A$ having the
$A_{\mathbb{E},\rho,\mathbb{F}Q}^{t^{2}}$-property, namely those $a \in A$
which satisfy $aQ\overline{a}^{\rho}=0$ (the reduced norm condition immediately
follows, since $a \not\in A^{\times}$ hence its reduced norm vanishes). We
denote the union of the set of those elements with
$A_{\mathbb{E},\rho,\mathbb{F}Q}^{t^{2}}$ by
$A_{\mathbb{E},\rho,\mathbb{F}Q}^{t^{2},0}$. It is no longer a group, but it is
closed under multiplication, and
$t:A_{\mathbb{E},\rho,\mathbb{F}Q}^{t^{2},0}\to\mathbb{F}$ (including 0) is
multiplicative. Moreover, apart from $(A_{\mathbb{E}}^{-})_{\rho,Q}$, we shall
also be needing the space $(A_{\mathbb{E}}^{-})_{\rho,\tilde{Q}}$. We shall
need a few simple relations between these sets.
\begin{lem}
Let $v\in(A_{\mathbb{E}}^{-})_{\rho,Q}$ and
$w\in(A_{\mathbb{E}}^{-})_{\rho,\tilde{Q}}$ be given. Then we have $(i)$
$\tilde{v}\in(A_{\mathbb{E}}^{-})_{\rho,\tilde{Q}}$. $(ii)$ $Q^{-1}vQ^{-1}$ is
also in $(A_{\mathbb{E}}^{-})_{\rho,\tilde{Q}}$. $(iii)$
$vQ^{-1} \in A_{\mathbb{E},\rho,\mathbb{F}Q}^{t^{2},0}$, with
$t(vQ^{-1})=-\frac{|v|^{2}}{|Q|^{2}}$. $(iv)$
$\tilde{w}\in(A_{\mathbb{E}}^{-})_{\rho,Q}$. $(v)$
$QwQ$ also lies in $(A_{\mathbb{E}}^{-})_{\rho,Q}$. $(vi)$ $Qw$ is an element
of $A_{\mathbb{E},\rho,\mathbb{F}Q}^{t^{2},0}$, whose multiplier $t(Qw)$ is
$-|Q|^{2}|w|^{2}$. \label{QthetaQrel}
\end{lem}

\begin{proof}
Lemma \ref{sp6gen} implies $v^{\rho}=-\frac{Q\tilde{v}Q}{|Q|^{2}}$, hence
$\tilde{v}=-|Q|^{2}Q^{-1}v^{\rho}Q^{-1}=-\frac{\tilde{Q}v^{\rho}\tilde{Q}}{
|\tilde{Q}|^{2}}$ since $Q^{-1}=\frac{\tilde{Q}}{|Q|^{2}}$ and
$|\tilde{Q}|^{2}=|Q|^{2}$. Applying $\rho$ to the latter equation and using the
fact that $\theta^{2}=Id_{A_{\mathbb{E}}^{-}}$ and that $Q$, hence also
$\tilde{Q}$, are $\rho$-invariant, yields the $\tilde{Q}$-based condition from
Lemma \ref{sp6gen} for $\tilde{v}$. This establishes part $(i)$. For part
$(ii)$, recall first that $\rho$ preserves the space
$(A_{\mathbb{E}}^{-})_{\rho,Q}$. It thus preserves also
$(A_{\mathbb{E}}^{-})_{\rho,\tilde{Q}}$. Write $v$ as $(v^{\rho})^{\rho}$, which
equals $-\frac{Q\tilde{v}^{\rho}Q}{|Q|^{2}}$ by Lemma \ref{sp6gen}. Hence
$Q^{-1}vQ^{-1}=-\frac{\tilde{v}^{\rho}}{|Q|^{2}}$, which implies part $(ii)$
since the latter element lies in $(A_{\mathbb{E}}^{-})_{\rho,\tilde{Q}}$ by part
$(i)$. For part $(iii)$ observe that as $Q^{-1}$ and $v$ lie in
$A_{\mathbb{E}}^{-}$ and $Q^{-\rho}=Q^{-1}$, the expression $vQ^{-1} \cdot
Q\cdot\overline{vQ^{-1}}^{\rho}$ equals $vQ^{-1}v^{\rho}$. Substituting the
expression for $v^{\rho}$ from Lemma \ref{sp6gen} again, the latter expression
becomes $-\frac{|v|^{2}}{Q^{2}}Q$, and part $(iii)$ follows since the reduced
norm condition is a consequence of Corollary \ref{NAvn2}. Parts $(iv)$ and
$(v)$ are proved either by applying the necessary changes in the proofs of
parts $(i)$ and $(ii)$ respectively, or since the maps given in parts $(i)$ and
$(ii)$ are injective maps between 6-dimensional vector space over $\mathbb{F}$
and the maps from parts $(iv)$ and $(v)$ are their inverses. For part $(vi)$ we
write $w$ as $Q^{-1}uQ^{-1}$ for some $u\in(A_{\mathbb{E}}^{-})_{\rho,Q}$ using
parts $(ii)$ and $(v)$, and then the assertion for $Qw=uQ^{-1}$ follows from
part $(iii)$ since $u=QwQ=-Qw\overline{Q}$ has vector norm $|Q|^{4}|w|^{2}$ by
Proposition \ref{NAFg} and Corollary \ref{NAvn2}. This completes the proof of
the lemma.
\end{proof}
We have seen that parts $(i)$ and $(iv)$ in Lemma \ref{QthetaQrel}, as well as
parts $(ii)$ and $(v)$ there, are inverses. Moreover, a claim similar to part
$(iii)$ (but without the multiplier) appears in the first assertion of Lemma
\ref{ref6gen}. We remark that the proof of parts $(iii)$ and $(iv)$ in that
lemma shows that if the vector $v$ or $w$ is anisotropic then the converse
implication also holds (cancel $vQ$ from the left in the proof of part $(iii)$,
and applying the same argument to extend it to part $(vi)$). On the other hand,
if $v$ (or $w$) are not isotropic then the converse implications in parts
$(iii)$ and $(iv)$ of Lemma \ref{QthetaQrel} may not hold. Indeed, by taking a
non-zero isotropic vector $v$ in $(A_{\mathbb{E}}^{-})_{\rho,Q}$ and
$z\in\mathbb{E}\setminus\mathbb{F}$, then $vQ^{-1}$ as well as $zvQ^{-1}$ lie in
$A_{\mathbb{E},\rho,\mathbb{F}Q}^{t^{2},0}$, while $zv$ no longer lies in
$(A_{\mathbb{E}}^{-})_{\rho,Q}$ (and the same for
$w\in\in(A_{\mathbb{E}}^{-})_{\rho,\tilde{Q}}$). Note that this argument does
not affect our assertions in the anisotropic case, since in this case $zv$
would have vector norm $z^{2}|v|^{2}$ and the multiplier of $zvQ^{-1}$ (which
does belong to $A_{\mathbb{E},\rho,\mathbb{F}Q}^{\times}$) is
$-N^{\mathbb{E}}_{\mathbb{F}}(z)\frac{|v|^{2}}{|Q|^{2}}$, which do not coincide
if $|v|^{2}\neq0$.

We shall also need the following complement of Lemma \ref{normsq} here:
\begin{lem}
For $\eta\in(A_{\mathbb{E}}^{-})_{\rho,Q}$ and
$\omega\in(A_{\mathbb{E}}^{-})_{\rho,\tilde{Q}}$, the element $1+\eta\omega$ of
$A$ lies in $A_{\mathbb{E},\rho,\mathbb{F}Q}^{t^{2},0}$, with multiplier
$D(\eta,\omega)$. \label{combAErhoFQm20}
\end{lem}

\begin{proof}
First note that the expression $D(\eta,\omega)$ from Lemma \ref{normsq} lies in
$\mathbb{F}$ for such $\eta$ and $\omega$, since $(A_{\mathbb{E}}^{-})_{\rho,Q}$
and $(A_{\mathbb{E}}^{-})_{\rho,\tilde{Q}}$ are quadratic spaces over
$\mathbb{F}$. Now, as in the proof of Lemma \ref{normsq}, we begin by assuming
that $\omega$ is anisotropic, and write $1+\eta\omega$ as
$\big(\tilde{\omega}+|\omega|^{2}\eta\big)\frac{\omega}{|\omega|^{2}}$. Then
$\frac{\omega}{|\omega|^{2}}\in(A_{\mathbb{E}}^{-})_{\rho,Q}$, and
$\tilde{\omega}+|\omega|^{2}\eta\in(A_{\mathbb{E}}^{-})_{\rho,\tilde{Q}}$ by
part $(iv)$ of Lemma \ref{QthetaQrel}. But now parts $(iii)$ and $(iv)$ of that
lemma show that $\big(\tilde{\omega}+|\omega|^{2}\eta\big)Q^{-1}$ and
$Q\frac{\omega}{|\omega|^{2}}$ are both in
$A_{\mathbb{E},\rho,\mathbb{F}Q}^{t^{2},0}$, with multipliers
$-\frac{|\tilde{\omega}+|\omega|^{2}\eta|}{|Q|^{2}}$ and
$-\frac{|Q|^{2}}{|\omega|^{2}}$ respectively, so that the assertion
follows from the multiplicativity of $A_{\mathbb{E},\rho,\mathbb{F}Q}^{t^{2},0}$
and $t:A_{\mathbb{E},\rho,\mathbb{F}Q}^{t^{2},0}\to\mathbb{F}$ and the fact that
the former vector norm was seen to be $|\omega|^{2}D(\eta,\omega)$ in the proof
of Lemma \ref{normsq}. For isotropic $\omega$, we use the polynomial method from
the proof of Lemma \ref{normsq} again, and consider
$\big(1+\eta(\omega+s\xi)\big)Q\overline{\big(1+\eta(\omega+s\xi)\big)}^{\rho}$
and $D(\eta,\omega+s\xi)Q$ as $A$-valued polynomials in $s$, for some fixed,
anisotropic $\xi\in(A_{\mathbb{E}}^{-})_{\rho,\tilde{Q}}$. This means that we
consider $A$ as a vector space over $\mathbb{F}$, we choose a basis for it which
includes $Q$, and consider the two sets of 16 polynomials in $s$ arising as the
coefficients using this basis (in the latter set, 15 polynomials will
identically vanish and one is the coefficient $D(\eta,\omega+s\xi)$). By what we
have proved, both sets of polynomials coincide for every $s$ perhaps maybe $s=0$
and $s=-\frac{2\langle\omega,\xi\rangle}{|\xi|^{2}}$, and the same argument as
in the proof of Lemma \ref{normsq} shows that they coincide for every $s$. The
reduced norm condition required for $A_{\mathbb{E},\rho,\mathbb{F}Q}^{t^{2},0}$
is satisfied by Lemma \ref{normsq} itself. Substituting $s=0$ verifies our
assertion also for isotropic $\omega$, which completes the proof of the lemma.
\end{proof}

\smallskip

Denote $GSp_{A_{\mathbb{E}}}\binom{1\ \ \ \ 0}{0\ \ -1}_{\rho,Q}$ the set of
those elements $\binom{a\ \ b}{c\ \ d}$ of $GSp_{A_{\mathbb{E}}}\binom{1\ \ \ \
0}{0\ \ -1}$, whose multipliers are in $\mathbb{F}^{\times}$, and which arise,
in terms of Corollary \ref{genform}, from parameters $v$ and $\alpha$ from
$(A_{\mathbb{E}}^{-})_{\rho,Q}$, $\beta$ which lies in
$(A_{\mathbb{E}}^{-})_{\rho,\tilde{Q}}$, and where $a$ is assumed to be in
$A_{\mathbb{E},\rho,\mathbb{F}Q}^{t^{2}}$. This definition is independent of the
choice of parameters, as one sees in the following
\begin{lem}
Let $c$, $w$, $\gamma$, and $\delta$ be a set of parameters for some element of
$GSp_{A_{\mathbb{E}}}\binom{1\ \ \ \ 0}{0\ \ -1}_{\rho,Q}$ as in Corollary
\ref{genform}. If $w\in(A_{\mathbb{E}}^{-})_{\rho,Q}$ then
$c \in A_{\mathbb{E},\rho,\mathbb{F}Q}^{t^{2}}$,
$\gamma\in(A_{\mathbb{E}}^{-})_{\rho,Q}$, and
$\delta\in(A_{\mathbb{E}}^{-})_{\rho,\tilde{Q}}$. \label{GSpAEQindep}
\end{lem}

\begin{proof}
By definition, there is a set $a \in A_{\mathbb{E},\rho,\mathbb{F}Q}^{t^{2}}$,
$v\in(A_{\mathbb{E}}^{-})_{\rho,Q}$, $\alpha\in(A_{\mathbb{E}}^{-})_{\rho,Q}$,
and $\beta\in(A_{\mathbb{E}}^{-})_{\rho,\tilde{Q}}$ of parameters from Corollary
\ref{genform} for that element. Given $w$ as an alternative parameter, the other
parameters $c$, $\gamma$, and $\delta$ are determined by the formulae from Lemma
\ref{comp}. The assumptions on $a$, $w$, $v$, and $\beta$ and Lemma
\ref{combAErhoFQm20} shows that the two multipliers in the expression for $c$ in
part $(i)$ of that lemma belong to $A_{\mathbb{E},\rho,\mathbb{F}Q}^{t^{2}}$.
Part $(i)$ of Lemma \ref{QthetaQrel} shows that the expression for $\delta$ in
part $(ii)$ of Lemma \ref{comp} is in $(A_{\mathbb{E}}^{-})_{\rho,\tilde{Q}}$,
as that the denominator $D(\beta,w-v)$ was seen to lie in $\mathbb{F}^{\times}$.
When we examine the expression for $\gamma$ in part $(iii)$ of Lemma \ref{comp},
we get the image of the action of $a^{-1}$ on some vector, which lies in
$(A_{\mathbb{E}}^{-})_{\rho,Q}$ by part $(iv)$ of Lemma \ref{QthetaQrel} and the
$\mathbb{F}$-rationality of $m$ and of the denominator $D(\beta,w-v)$ again.
Since $a$, hence also $a^{-1}$, comes from
$A_{\mathbb{E},\rho,\mathbb{F}Q}^{t^{2}}$, this expression also lies in
$(A_{\mathbb{E}}^{-})_{\rho,Q}$ by Lemma \ref{ac6gen}. This completes the proof
of the lemma.
\end{proof}

We can now prove
\begin{prop}
The set $GSp_{A_{\mathbb{E}}}\binom{1\ \
\ \ 0}{0\ \ -1}_{\rho,Q}$ is a subgroup of $GSp_{A_{\mathbb{E}}}\binom{1\ \ \ \
0}{0\ \ -1}$, which is stable under $\rho$. It is contained in
$GSp_{A_{\mathbb{E}}}^{\mathbb{E}^{2}}\binom{1\ \ \ \ 0}{0\ \ -1}$ and comes
endowed with a splitting map into the double cover
$\widetilde{GSp}_{A_{\mathbb{E}}}^{\mathbb{E}^{2}}\binom{1\ \ \
\ 0}{0\ \ -1}$. \label{GSpAErhoFQ}
\end{prop}

\begin{proof}
As in the proofs of Proposition \ref{GSphom} and Theorem \ref{GSppsi}, given
two elements of $GSp_{A_{\mathbb{E}}}\binom{1\ \ \ \ 0}{0\ \ -1}_{\rho,Q}$, we
may assume that in the form of Corollary \ref{genform}, the left multiplier $g$
has parameters $a$, $v$, $\alpha$ and $\beta$, the right multiplier $h$ arises
from $e$, $z$, $\kappa$, and $\nu$, and $x$, the same $v$, $\xi$, and $\zeta$
are parameters for the product $gh$. We may further assume that $v$ and $z$ lie
in $(A_{\mathbb{E}}^{-})_{\rho,Q}$. Lemma \ref{GSpAEQindep} implies that
$\alpha$ and $\kappa$ also come from the same space, $\beta$ and $\nu$ are in
$(A_{\mathbb{E}}^{-})_{\rho,\tilde{Q}}$, and $a$ and $e$ are elements of
$A_{\mathbb{E},\rho,\mathbb{F}Q}^{t^{2}}$. We have to show that the remaining
parameters for $gh$ lie in the appropriate sets. For $x$ this follows from the
properties of $a$, $e$, $\alpha$, $z$, and $\nu$ by Lemma \ref{combAErhoFQm20}.
Invoking part $(iv)$ of Lemma \ref{QthetaQrel}, as well as Lemma \ref{ac6gen}
for the action of $e^{-1} \in A_{\mathbb{E},\rho,\mathbb{F}Q}^{t^{2}}$, verifies
the assertion for $\xi$, since the multiplier $n$ and the denominator
$D(\alpha+z,\nu)$ lie in $\mathbb{F}^{\times}$ by the proof of Lemma
\ref{combAErhoFQm20}. Now, $\zeta$ is the sum of $\beta$ and the image of the
action of $\overline{a}^{\ -1}$ on a vector, which is contained in
$(A_{\mathbb{E}}^{-})_{\rho,\tilde{Q}}$ by part $(i)$ of Lemma
\ref{QthetaQrel}. By parts $(i)$ and $(iv)$ of the latter lemma, it is
sufficient to show that
$\tilde{\zeta}-\tilde{\beta}\in(A_{\mathbb{E}}^{-})_{\rho,Q}$. But using Lemma
\ref{AxA-rel} and part $(i)$ of Lemma \ref{QthetaQrel} again, the latter vector
is obtained by an element of $(A_{\mathbb{E}}^{-})_{\rho,Q}$ by the action of
$a \in A_{\mathbb{E},\rho,\mathbb{F}Q}^{t^{2}}$ (up to scalars from
$\mathbb{F}^{\times}$), which verifies the assertion for $\zeta$ as well. Hence
$GSp_{A_{\mathbb{E}}}\binom{1\ \ \ \ 0}{0\ \ -1}_{\rho,Q}$ is a subgroup of
$GSp_{A_{\mathbb{E}}}\binom{1\ \ \ \ 0}{0\ \ -1}$. For the stability under
$\rho$, we may write any element of $GSp_{A_{\mathbb{E}}}\binom{1\ \ \ \ 0}{0\ \
-1}_{\rho,Q}$ as in Corollary \ref{genform}, with parameters as in Lemma
\ref{GSpAEQindep}. The fact that $\rho$ preserves
$(A_{\mathbb{E}}^{-})_{\rho,Q}$, $(A_{\mathbb{E}}^{-})_{\rho,\tilde{Q}}$, and
$A_{\mathbb{E},\rho,\mathbb{F}Q}^{t^{2}}$ implies the preservation of
$GSp_{A_{\mathbb{E}}}\binom{1\ \ \ \ 0}{0\ \ -1}_{\rho,Q}$ as well. Now, the
value of the map $\varphi$ from Proposition \ref{GSphom} is based only on the
parameter from $A_{\mathbb{E}}^{\times}$ in Corollary \ref{genform}. As for
matrices in $GSp_{A_{\mathbb{E}}}\binom{1\ \ \ \ 0}{0\ \ -1}_{\rho,Q}$ these
parameters come from $A_{\mathbb{E},\rho,\mathbb{F}Q}^{t^{2}}$ and the latter
group is contained in $A_{\mathbb{E}}^{(\mathbb{E}^{\times})^{2}}$ by the
definition of the former group in Lemma \ref{ind2int}, our subgroup is contained
in $GSp_{A_{\mathbb{E}}}^{\mathbb{E}^{2}}\binom{1\ \ \ \ 0}{0\ \ -1}$. In
addition, the double cover
$\widetilde{GSp}_{A_{\mathbb{E}}}^{\mathbb{E}^{2}}\binom{1\ \ \ \ 0}{0\ \ -1}$
is defined by adding a choice of a square root for the reduced norm of the
parameter from $A_{\mathbb{E}}^{(\mathbb{E}^{\times})^{2}}$, so that we get a
parameter from $\widetilde{A}_{\mathbb{E}}^{(\mathbb{E}^{\times})^{2}}$ for
elements of this double cover. But
$\widetilde{A}_{\mathbb{E}}^{(\mathbb{E}^{\times})^{2}}$ was seen to split over
$A_{\mathbb{E},\rho,\mathbb{F}Q}^{t^{2}}$ via $g\mapsto\big(g,t(g)\big)$.
Moreover, this splitting map is compatible with parameter changes inside
$A_{\mathbb{E},\rho,\mathbb{F}Q}^{t^{2}}$,
$(A_{\mathbb{E}}^{-})_{\rho,Q}$, and
$(A_{\mathbb{E}}^{-})_{\rho,\tilde{Q}}$ by part $(i)$ of Lemma \ref{comp} and
Lemma \ref{combAErhoFQm20}. This establishes the splitting of
$\widetilde{GSp}_{A_{\mathbb{E}}}^{\mathbb{E}^{2}}\binom{1\ \ \ \ 0}{0\ \ -1}$
over $GSp_{A_{\mathbb{E}}}\binom{1\ \ \ \ 0}{0\ \ -1}_{\rho,Q}$ as well, and
completes the proof of the proposition.
\end{proof}

In case we wish to evaluate $\psi$-images of elements of
$GSp_{A_{\mathbb{E}}}\binom{1\ \ \ \ 0}{0\ \ -1}_{\rho,Q}$ using other entries,
we cannot use the matrix $\binom{0\ \ 1}{1\ \ 0}$, as it does not belong to this
group. However, the element $\binom{0\ \ Q}{\tilde{Q}\ \ 0}$, of multiplier
$-|Q|^{2}$, does belong there: By choosing some non-zero $h\in\mathbb{E}_{0}$,
we recall that $hQ$ and $\frac{Q}{h|Q|^{2}}$ are in
$(A_{\mathbb{E}}^{-})_{\rho,Q}$,
$\frac{\tilde{Q}}{h|Q|^{2}}\in(A_{\mathbb{E}}^{-})_{\rho,Q}$, and $h|Q|^{2}$ is
in $A_{\mathbb{E},\rho,\mathbb{F}Q}^{t^{2}}$, and our element may be obtained
from the parameters $v=-hQ$, $a=h|Q|^{2}$, $\alpha=\frac{Q}{h|Q|^{2}}$, and
$\beta=\frac{\tilde{Q}}{h|Q|^{2}}$. Recall that the multiplier $t(h|Q|^{2})$ of
$h|Q|^{2} \in A_{\mathbb{E},\rho,\mathbb{F}Q}^{t^{2}}$ is $-h^{2}|Q|^{4}$, so
that as an element of $GSp_{A_{\mathbb{E}}}\binom{1\ \ \ \ 0}{0\ \ -1}_{\rho,Q}$
we must have $\psi\binom{0\ \ Q}{\tilde{Q}\ \ 0}=\binom{\ \ 0\ \ -\tilde{Q}}{-Q\
\ \ \ 0}$. In any case, we may use this element in order to transfer the formula
for $\psi$ on $GSp_{A_{\mathbb{E}}}\binom{1\ \ \ \ 0}{0\ \ -1}_{\rho,Q}$ to be
based on the other matrix entries.

As in the situation we encountered in dimension 6 and general discriminant, we
remark that unless $Q^{2}$ is a scalar and $Q$ and $\tilde{Q}$ span the same
vector space, the space $(A_{\mathbb{E}}^{-})_{\rho,Q} \oplus H$ is not
invariant under the linear automorphism $\hat{\psi}$ from Lemma \ref{vnorm8}. In
addition, the automorphism $\psi$ does now preserve the subgroup
$GSp_{A_{\mathbb{E}}}\binom{1\ \ \ \ 0}{0\ \ -1}_{\rho,Q}$ (embedded via the
splitting map from Proposition \ref{GSpAErhoFQ}) in this case. However, we can
still pursue the usual route by using the following
\begin{lem}
Let $\hat{Q}$ denote the element of
$\widetilde{GSp}_{A_{\mathbb{E}}}^{\mathbb{E}^{2}}\binom{1\ \ \
\ 0}{0\ \ -1}$ lying over $\binom{Q\ \ \ \ 0}{0\ \ -\tilde{Q}}$ in which the
square root $t$ of $N^{A_{\mathbb{E}}}_{\mathbb{E}}(Q)$ is chosen to be
$|Q|^{2}$. Then the element $\hat{Q}\tilde{\psi}$ of the semi-direct product of
$\{1,\tilde{\psi}\}$ and
$\widetilde{GSp}_{A_{\mathbb{E}}}^{\mathbb{E}^{2}}\binom{1\ \ \ \ 0}{0\ \ -1}$
from Lemma \ref{ac8id1} squares to $-|Q|^{2}$ (with $|Q|^{4}$ as the square root
of its reduced norm) in that semi-direct product. Its action on
$(A_{\mathbb{E}}^{-})_{\rho,Q} \oplus H$ coincides with that of $\rho$ (which
preserves it), and conjugation by this element operates on
$\widetilde{GSp}_{A_{\mathbb{E}}}^{\mathbb{E}^{2}}\binom{1\ \ \ \ 0}{0\ \ -1}$
via $g\mapsto\hat{Q}\psi(g)\hat{Q}^{-1}$. This automorphism of
$\widetilde{GSp}_{A_{\mathbb{E}}}^{\mathbb{E}^{2}}\binom{1\ \ \ \ 0}{0\ \ -1}$
has order 2, and it preserves the subgroup $GSp_{A_{\mathbb{E}}}\binom{1\ \ \ \
0}{0\ \ -1}_{\rho,Q}$ embedded through the splitting map from Proposition
\ref{GSpAErhoFQ} as its operation on the latter group is the same as that of
$\rho$. \label{Qhatpsi}
\end{lem}

\begin{proof}
The multiplier of $\hat{Q}$ is $|Q|^{2}$, and as $Q \in A_{\mathbb{E}}^{-}$ and
$t=|Q|^{2}$ indeed satisfies $t^{2}=N^{A}_{\mathbb{F}}(Q)$ by Corollary
\ref{NAvn2}, the definition of the map $\psi$ in Theorem \ref{GSppsi} shows that
$\psi(\hat{Q})$ is $\binom{-\tilde{Q}\ \ 0}{\ \ 0\ \ Q}$ (with the same square
root $-|Q|^{2}$). The first assertion follows immediately from the fact that
$\hat{Q}\psi(\hat{Q})=-|Q|^{2}I$ and the product of the square roots is
$|Q|^{4}$ (the number $D(\alpha+z,w)$ appearing in the proof of Theorem
\ref{GSppsi} in the square root of the reduced norm of $x$ from Lemma
\ref{GSpprod} is just 1). Now, $\tilde{\psi}$ operates as $\hat{\psi}$ on
$A_{\mathbb{E}}^{-} \oplus H_{\mathbb{E}}$ (which is just $\theta$ on the
$A_{\mathbb{E}}^{-}$ part), and the operation of the diagonal element $\hat{Q}$
was evaluated in Proposition \ref{GSppresHA-}: The $H_{\mathbb{E}}$ part is
pointwise fixed since $t=m$ for our element, and the combined operation on
$A_{\mathbb{E}}^{-}$ is via $u\mapsto-\frac{Q\tilde{u}Q}{|Q|^{2}}$. But Lemma
\ref{sp6gen} shows that on $(A_{\mathbb{E}}^{-})_{\rho,Q}$ the latter map
coincides with $\rho$, and as $\rho$ leaves $H$ also pointwise fixed, this
establishes the second assertion. The formula for the conjugation on
$\widetilde{GSp}_{A_{\mathbb{E}}}^{\mathbb{E}^{2}}\binom{1\ \ \ \ 0}{0\ \ -1}$
follows directly from the structure of the semi-direct product in Lemma
\ref{ac8id1}, and it is of order 2 either since $-|Q|^{2}I$ (with the square
root $|Q|^{4}$ of $N^{A_{\mathbb{E}}}_{\mathbb{E}}(-|Q|^{2})$) operates
trivially or by a direct evaluation. In order to examine the action of this
automorphism on $GSp_{A_{\mathbb{E}}}\binom{1\ \ \ \ 0}{0\ \ -1}_{\rho,Q}$ in
view of that of $\rho$, we recall the relation between $\hat{Q}$ and
$\psi(\hat{Q})$, so that we evaluate $-\frac{1}{|Q|^{2}}\hat{Q}g\psi(\hat{Q})$
for $g \in GSp_{A_{\mathbb{E}}}\binom{1\ \ \ \ 0}{0\ \ -1}_{\rho,Q}$ and
compare it with $g^{\rho}$. It suffices to take $g$ from a set of generators of
the latter group, namely $\binom{1\ \ v}{0\ \ 1}$ with
$v\in(A_{\mathbb{E}}^{-})_{\rho,Q}$, $\binom{1\ \ 0}{w\ \ 1}$ with
$w\in(A_{\mathbb{E}}^{-})_{\rho,\tilde{Q}}$, and $\binom{a\ \ \ \ 0\ \ \ }{0\ \
m\overline{a}^{-1}}$ where $a$ lies in
$A_{\mathbb{E},\rho,\mathbb{F}Q}^{t^{2}}$. Now, $\psi$ replaces $v$ and $w$ by
their $\theta$-images and $a$ by $t(a)\overline{a}^{\ -1}$, and after
conjugating by $\hat{Q}$ and using Lemma \ref{vnorm6} we find that $v$ is sent
to $-\frac{Q\tilde{v}Q}{|Q|^{2}}$, $w$ is taken to
$-\frac{\tilde{Q}\tilde{w}\tilde{Q}}{|Q|^{2}}$, the diagonal entries of the
unipotent generators remain invariant, and $a$ is mapped to
$t(a)Q\overline{a}^{\ -1}Q^{-1}$ (the other entry $\overline{a}^{\ -1}$ becomes
$\frac{Q^{-1}aQ}{t(a)}$). But as we assume that
$v\in(A_{\mathbb{E}}^{-})_{\rho,Q}$,
$w\in(A_{\mathbb{E}}^{-})_{\rho,\tilde{Q}}$, and $a \in
A_{\mathbb{E},\rho,\mathbb{F}Q}^{t^{2}}$, the first two expressions are
$v^{\rho}$ and $w^{\rho}$ by Lemma \ref{sp6gen}, while the proof of Lemma
\ref{Qtheta} shows that the latter expression is just $a^{\rho}$ (and the one
in parentheses is $\overline{a}^{\ -\rho}$). Since the multiplier $m$ lies in
$\mathbb{F}$ and is thus $\rho$-invariant, this establishes the coincidence of
$\rho$ and conjugation by $\hat{Q}\tilde{\psi}$ on
$GSp_{A_{\mathbb{E}}}\binom{1\ \ \ \ 0}{0\ \ -1}_{\rho,Q}$, and as $\rho$ was
seen to preserve this group in Proposition \ref{GSpAErhoFQ}, conjugation by
$\hat{Q}\tilde{\psi}$ does the same. This completes the proof of the lemma.
\end{proof}

Note that $\Xi=\hat{Q}$ satisfies the conditions of Corollary \ref{alt8id1},
and the operations of $\hat{Q}\tilde{\psi}$ appearing in Lemma \ref{Qhatpsi} are
just $\hat{\psi}_{\hat{Q}}$ and $\psi_{\hat{Q}}$ respectively. Now, Lemma
\ref{Qhatpsi} gives an intrinsic description of $GSp_{A_{\mathbb{E}}}\binom{1\ \
\ \ 0}{0\ \ -1}_{\rho,Q}$, and we can also establish some properties of its
entries and the elements from the $GSp$ relations:
\begin{cor}
The subgroup $GSp_{A_{\mathbb{E}}}\binom{1\ \ \ \ 0}{0\ \ -1}_{\rho,Q}$ is
characterized as the set of elements of
$\widetilde{GSp}_{A_{\mathbb{E}}}^{\mathbb{E}^{2}}\binom{1\ \ \ \ 0}{0\ \ -1}$
on which the automorphism $\psi_{\hat{Q}}$ operates as $\rho$. \label{GSprhoQst}
\end{cor}

\begin{proof}
The fact that $\psi_{\hat{Q}}(g)=g^{\rho}$ for $g \in
GSp_{A_{\mathbb{E}}}\binom{1\ \ \ \ 0}{0\ \ -1}_{\rho,Q}$ was seen in the proof
of Lemma \ref{Qhatpsi}. Conversely, let
$g\in\widetilde{GSp}_{A_{\mathbb{E}}}^{\mathbb{E}^{2}}\binom{1\ \ \ \ 0}{0\ \
-1}$ be an element satisfying $\psi_{\hat{Q}}(g)=g^{\rho}$. Then the multipliers
of both sides coincide, and as $\psi$ commutes with $m$ and
$m(g^{\rho})=m(g)^{\rho}$ we find that $m(g)\in\mathbb{F}^{\times}$. Morever,
Lemma \ref{genie} allows us to find a set of parameters for Corollary
\ref{genform} for $g$ in which $v\in(A_{\mathbb{E}}^{-})_{\rho,Q}$. As
$\binom{1\ \ v}{0\ \ 1}$ lies in $GSp_{A_{\mathbb{E}}}\binom{1\ \ \ \ 0}{0\ \
-1}_{\rho,Q}$, it suffices to consider elements
$g\in\widetilde{GSp}_{A_{\mathbb{E}}}^{\mathbb{E}^{2}}\binom{1\ \ \ \ 0}{0\ \
-1}$ having invertible upper left entry. But these elements have a unique
decomposition as in Lemma \ref{invent}. Comparing these decompositions for
$g^{\rho}$ and $\psi_{\hat{Q}}(g)$ and recalling that the unipotent matrices are
assumed to have unipotent $\psi$-images (and not with $-1$ on the diagonal)
reduces the verification to the multipliers appearing in Lemma \ref{invent},
lifted into $\widetilde{GSp}_{A_{\mathbb{E}}}^{\mathbb{E}^{2}}\binom{1\ \ \ \
0}{0\ \ -1}$. But these verifications are carried out in the proof of Lemma
\ref{Qhatpsi}. This proves the corollary.
\end{proof}
One can also show that if $\binom{a\ \ b}{c\ \ d}$ is an element of
$GSp_{A_{\mathbb{E}}}\binom{1\ \ \ \ 0}{0\ \ -1}_{\rho,Q}$ then $a$, $bQ^{-1}$,
$Qc$, and $QdQ^{-1}$ all lie in $A_{\mathbb{E},\rho,\mathbb{F}Q}^{t^{2},0}$, the
elements $a\overline{b}=-b\overline{a}$ and $\overline{b}d=-\overline{d}b$ of
$A_{\mathbb{E}}^{-}$ belong to $(A_{\mathbb{E}}^{-})_{\rho,Q}$, and 
$\overline{a}c=-\overline{c}a$ and $\overline{c}d=-\overline{d}c$ come from
$(A_{\mathbb{E}}^{-})_{\rho,\tilde{Q}}$. In fact, conjugating by $\binom{1\ \
0}{0\ \ Q}$ as in Corollary \ref{iso8id1} yields a subgroup
$GSp_{A_{\mathbb{E}}}\binom{Q\ \ 0}{0\ \ Q}_{\rho,Q}$ of
$GSp_{A_{\mathbb{E}}}\binom{Q\ \ 0}{0\ \ Q}$ with a simpler description: All the
entries of elements $\binom{e\ \ f}{g\ \ h}$ of that group come from
$A_{\mathbb{E},\rho,\mathbb{F}Q}^{t^{2},0}$, and for which the elements
$eQ\overline{f}=fQ\overline{e}$ and $gQ\overline{h}=hQ\overline{g}$ of
$A_{\mathbb{E}}^{-}$ lie in $(A_{\mathbb{E}}^{-})_{\rho,Q}$, while
$\overline{e}Q^{-1}g=\overline{g}Q^{-1}e$ and
$\overline{f}Q^{-1}h=\overline{h}Q^{-1}f$ are in
$(A_{\mathbb{E}}^{-})_{\rho,\tilde{Q}}$ (use parts $(ii)$ and $(v)$ of Lemma
\ref{QthetaQrel} again). Moreover, when many of these terms are
invertible, some of the latter assertions follow from one another---see Lemma
\ref{QthetaQrel} and the remarks following it. However, once may non-zero
entries which are not invertible are involved, there may be some matrices
satisfying these conditions which do not belong to
$GSp_{A_{\mathbb{E}}}\binom{1\ \ \ \ 0}{0\ \ -1}_{\rho,Q}$. Hence we content
ourselves with the description appearing in Corollary \ref{GSprhoQst}. 

The reason for considering $GSp_{A_{\mathbb{E}}}\binom{1\ \ \ \ 0}{0\ \
-1}_{\rho,Q}$ is given in the following
\begin{lem}
Consider the operation of $GSp_{A_{\mathbb{E}}}\binom{1\ \ \ \ 0}{0\ \
-1}_{\rho,Q}$, viewed as a subgroup of
$\widetilde{GSp}_{A_{\mathbb{E}}}^{\mathbb{E}^{2}}\binom{1\ \ \ \ 0}{0\ \ -1}$
via the lift from Proposition \ref{GSpAErhoFQ}, on the $\mathbb{E}$-vector space
$A_{\mathbb{E}}^{-} \oplus H_{\mathbb{E}}$. This action preserves the
$\mathbb{F}$-subspace $(A_{\mathbb{E}}^{-})_{\rho,Q} \oplus H$. The group
generated by $GSp_{A_{\mathbb{E}}}\binom{1\ \ \ \ 0}{0\ \ -1}_{\rho,Q}$ and the
element $\hat{Q}\tilde{\psi}$ from Lemma \ref{Qhatpsi} containg the latter
group with index 2, and maps into $O\big((A_{\mathbb{E}}^{-})_{\rho,Q} \oplus
H\big)$ as well. \label{ac8igen}
\end{lem}

\begin{proof}
As in Proposition \ref{GSppresHA-}, it suffices to prove the assertion for a
generating subset of the subgroup. Corollary \ref{genform} and the definition of
$GSp_{A_{\mathbb{E}}}\binom{1\ \ \ \ 0}{0\ \ -1}_{\rho,Q}$ show that the set
consisting of unipotent matrices $\binom{1\ \ v}{0\ \ 1}$ with
$v\in(A_{\mathbb{E}}^{-})_{\rho,Q}$ and $\binom{1\ \ 0}{w\ \ 1}$ where
$w\in(A_{\mathbb{E}}^{-})_{\rho,\tilde{Q}}$ together with the subgroup of
diagonal matrices $\binom{a\ \ \ \ 0\ \ \ }{0\ \ m\overline{a}^{-1}}$ with $a
\in A_{\mathbb{E},\rho,\mathbb{F}Q}^{t^{2}}$ and $m\in\mathbb{F}^{\times}$ is a
such a generating set. The action of these generators on an arbitrary element
$\binom{u\ \ -p}{q\ \ -\tilde{u}} \in A_{\mathbb{E}}^{-} \oplus H_{\mathbb{E}}$
was seen in the proof of Proposition \ref{GSppresHA-} to be as follows: The $u$
coordinate becomes $u+qv$, $u+p\tilde{w}$, and $\frac{au\overline{a}}{t(a)}$
(recall the choice of the $\psi$-image), $p$ is sent to $p+2\langle u,v
\rangle+q|v|^{2}$, $p$, and $\frac{t(a)p}{m}$, and $q$ is mapped to $q$,
$q+2\langle u,\tilde{w} \rangle+q|w|^{2}$, and $\frac{mq}{t(a)}$, respectively.
Our assumptions on $v$, $w$, $m$, and $a$ show, using part $(iv)$ of Lemma
\ref{QthetaQrel} for $w$ and Lemma \ref{ac6gen} for $a$, that if $p$ and $q$ are
from $\mathbb{F}$ and $u\in(A_{\mathbb{E}}^{-})_{\rho,Q}$ then the same
assertion holds for their images. The remaining assertions follow, as in the
proof of Lemma \ref{ac6gen}, from Lemma \ref{Qhatpsi}, Lemma \ref{ac8id1} over
$\mathbb{E}$, and the fact that $(A_{\mathbb{E}}^{-})_{\rho,Q} \oplus H$
inherits its quadratic structure from $A_{\mathbb{E}}^{-} \oplus
H_{\mathbb{E}}$. This proves the lemma.
\end{proof}

The assertion about reflections in this case appears in the following
\begin{lem}
Let $g$ be an anisotropic element of $(A_{\mathbb{E}}^{-})_{\rho,Q} \oplus H$.
Then the product $g\hat{Q}^{-1}$ belongs to $GSp_{A_{\mathbb{E}}}\binom{1\ \ \ \
0}{0\ \ -1}_{\rho,Q}$, and if we compose the action of the element
$\hat{Q}\tilde{\psi}$ from Lemma \ref{Qhatpsi} with that of the latter
composition we obtain the reflection in $g$ on $(A_{\mathbb{E}}^{-})_{\rho,Q}
\oplus H$. \label{ref8igen}
\end{lem}

\begin{proof}
First, the fact that $(A_{\mathbb{E}}^{-})_{\rho,Q} \oplus H$ is a quadratic
space over $\mathbb{F}$ means that all the multipliers are from $\mathbb{F}$.
Consider the element $g\hat{Q}^{-1}$, for anisotropic $g=\binom{u\ \ -p}{q\ \
-\tilde{u}}\in(A_{\mathbb{E}}^{-})_{\rho,Q} \oplus H$. If $p\neq0$, we multiply
it by $\binom{0\ \ Q}{\tilde{Q}\ \ 0}$ from the right. As the product of
$\hat{Q}^{-1}=\frac{\psi(\hat{Q})}{|Q|^{2}}$ and the latter element yields the
matrix $\binom{\ \ 0\ \ 1}{-1\ \ 0}$ from the proof of Lemma \ref{vnorm8}, we
may use the parameters given in that lemma for this case. As the scalar
$p\in\mathbb{F}^{\times}$ lies in $A_{\mathbb{E},\rho,\mathbb{F}Q}^{t^{2}}$, the
vector $\frac{u}{p}$ is in $(A_{\mathbb{E}}^{-})_{\rho,Q}$ by our assumption on
$U$, and $\frac{\tilde{u}}{p}\in(A_{\mathbb{E}}^{-})_{\rho,\tilde{Q}}$ by part
$(i)$ of Lemma \ref{QthetaQrel}, we indeed get $g\hat{Q}^{-1} \in
GSp_{A_{\mathbb{E}}}\binom{1\ \ \ \ 0}{0\ \ -1}_{\rho,Q}$ if $p\neq0$. with
$p=0$, so that $u \in A_{\mathbb{E}}^{\times}$, we get a matrix of multiplier
$\frac{|u|^{2}}{|Q|^{2}}$ for which the parameters may be taken to be
$v=\alpha=0$, $a=uQ^{-1}$, and $\beta=\frac{q\tilde{u}}{|u|^{2}}$. Parts $(i)$
and $(iii)$ of Lemma \ref{QthetaQrel} show that $g\hat{Q}^{-1} \in
GSp_{A_{\mathbb{E}}}\binom{1\ \ \ \ 0}{0\ \ -1}_{\rho,Q}$ in case $p=0$ as well.
We must, however, consider these elements in the double cover
$\widetilde{GSp}_{A_{\mathbb{E}}}\binom{1\ \ \ \ 0}{0\ \ -1}$. Recall from the
proof of Lemma \ref{vnorm8} that the branch of $\psi$ which coincides with
$\hat{\psi}$ is determined for $p\neq0$ by the condition that after right
multiplication by $\binom{\ \ 0\ \ 1}{-1\ \ 0}$ (considered as its own
$\psi$-image) the square root of $N^{A}_{\mathbb{F}}(p)=p^{4}$ is $p^{2}$, while
for $p=0$ we take $-|u|^{2}$ for the square root of $N^{A}_{\mathbb{F}}(u)$. One
verifies that the chosen $\psi$-image in the definition of $\hat{Q}$ and
$\hat{Q}^{-1}$ in Lemma \ref{Qhatpsi} and the $\psi$-image of $\binom{0\ \
Q}{\tilde{Q}\ \ 0}$ as an element of $GSp_{A_{\mathbb{E}}}\binom{1\ \ \ \ 0}{0\
\ -1}_{\rho,Q}$ combine to give $\binom{\ \ 0\ \ 1}{-1\ \ 0}$ with itself as its
$\psi$-image, and indeed $p^{2}$ is the multiplier of $p \in
A_{\mathbb{E},\rho,\mathbb{F}Q}^{t^{2}}$ as a parameter of $g\binom{\ \ 0\ \
1}{-1\ \ 0}$. On the other hand, part $(iii)$ of Lemma \ref{QthetaQrel} shows
that the multiplier of the parameter $uQ^{-1}$ is $-\frac{|u|^{2}}{|Q|^{2}}$,
and multiplying it by the chosen square root $|Q|^{2}$ of the diagonal element
$\hat{Q}$ yields the desired value $-|u|^{2}$. It follows that the composition
of $g\hat{Q}^{-1}$ and $\hat{Q}\tilde{\psi}$ is the combination $g\tilde{\psi}$,
with $g \in \widetilde{GSp}_{A_{\mathbb{E}}}\binom{1\ \ \ \ 0}{0\ \ -1}$ whose
$\psi$-image coincides with $\hat{\psi}(g)$. The action of this element on
$A_{\mathbb{E}}^{-} \oplus H_{\mathbb{E}}$ was seen in Lemma \ref{ref8id1} (over
$\mathbb{E}$) to be the reflection in $g$, and this assertion descends to
$(A_{\mathbb{E}}^{-})_{\rho,Q} \oplus H$ for anisotropic $g$ which is taken from
the latter space. This completes the proof of lemma.
\end{proof}

Observe that Lemma \ref{ref8igen} in fact uses the alternative representation
appearing in Corollary \ref{alt8id1} with $\Xi=\hat{Q}$, with right
multiplication by $\hat{Q}^{-1}$ on the space and $\psi_{\hat{Q}}$ as the
automorphism of the group.

Now we are in a position to prove
\begin{thm}
We have $Gspin\big((A_{\mathbb{E}}^{-})_{\rho,Q} \oplus
H\big)=GSp_{A_{\mathbb{E}}}\binom{1\ \ \ \ 0}{0\ \ -1}_{\rho,Q}$, and the spin
group $spin\big((A_{\mathbb{E}}^{-})_{\rho,Q} \oplus H\big)$ is the subgroup
$Sp_{A_{\mathbb{E}}}\binom{1\ \ \ \ 0}{0\ \ -1}_{\rho,Q}$ consisting of those
elements of $GSp_{A_{\mathbb{E}}}\binom{1\ \ \ \ 0}{0\ \ -1}_{\rho,Q}$
having multiplier 1 (namely the intersection of $GSp_{A_{\mathbb{E}}}\binom{1\ \
\ \ 0}{0\ \ -1}_{\rho,Q}$ with $\widetilde{Sp}_{A_{\mathbb{E}}}\binom{1\ \ \ \
0}{0\ \ -1}$). \label{dim8igen}
\end{thm}

\begin{proof}
We have a surjective map from the semi-direct product of $\{1,\tilde{\psi}\}$
with $\widetilde{GSp}_{A_{\mathbb{E}}}\binom{1\ \ \ \ 0}{0\ \ -1}$ (which we may
consider as generated by $\hat{Q}\tilde{\psi}$ and the latter group) onto
$O(A_{\mathbb{E}}^{-} \oplus H_{\mathbb{E}})$, whose kernel is
$\mathbb{E}^{\times}$, and such that the inverse image of
$SO(A_{\mathbb{E}}^{-} \oplus H_{\mathbb{E}})$ is
$\widetilde{GSp}_{A_{\mathbb{E}}}\binom{1\ \ \ \ 0}{0\ \ -1}$. By Lemma
\ref{ac8igen}, the subgroup which $\hat{Q}\tilde{\psi}$ generates with
$GSp_{A_{\mathbb{E}}}\binom{1\ \ \ \ 0}{0\ \ -1}_{\rho,Q}$ is sent to the
subgroup $O\big((A_{\mathbb{E}}^{-})_{\rho,Q} \oplus H\big)$ of
$O(A_{\mathbb{E}}^{-} \oplus H_{\mathbb{E}})$. Invoking Lemma \ref{ref8igen} and
Proposition \ref{CDT}, we obtain that this subgroup surjects onto
$O\big((A_{\mathbb{E}}^{-})_{\rho,Q} \oplus H\big)$. By Theorem \ref{dim8id1}
and the fact that the determinant commutes with the injection of
$O\big((A_{\mathbb{E}}^{-})_{\rho,Q} \oplus H\big)$ into $O(A_{\mathbb{E}}^{-}
\oplus H_{\mathbb{E}})$, we find that an element of the group generated by
$\hat{Q}\tilde{\psi}$ and $GSp_{A_{\mathbb{E}}}\binom{1\ \ \ \ 0}{0\ \
-1}_{\rho,Q}$ is sent to $SO\big((A_{\mathbb{E}}^{-})_{\rho,Q} \oplus H\big)$
if and only if it comes from $GSp_{A_{\mathbb{E}}}\binom{1\ \ \ \ 0}{0\ \
-1}_{\rho,Q}$ alone. It follows that the map $GSp_{A_{\mathbb{E}}}\binom{1\ \ \
\ 0}{0\ \ -1}_{\rho,Q} \to O\big((A_{\mathbb{E}}^{-})_{\rho,Q} \oplus H\big)$
is surjective. An element $GSp_{A_{\mathbb{E}}}\binom{1\ \ \ \ 0}{0\ \
-1}_{\rho,Q}$ acts trivially if and only if it is a scalar matrix $rI$ (with
$r\in\mathbb{E}^{\times}$) such that its $\rho$-image is the same as
$\hat{Q}\psi(rI)\hat{Q}^{-1}$, and moreover we require $\psi(rI)$ to be $+rI$ in
the double cover $\widetilde{GSp}_{A_{\mathbb{E}}}\binom{1\ \ \ \ 0}{0\ \ -1}$.
This happens if and only if $r\in\mathbb{F}^{\times}$ (non-zero elements
$r\in\mathbb{E}_{0}$ satisfy the first condition but have the wrong sign of
$\psi(rI)$, hence they do not operate trivially but rather as
$-Id_{(A_{\mathbb{E}}^{-})_{\rho,Q} \oplus H}$). This proves that
$Gspin\big((A_{\mathbb{E}}^{-})_{\rho,Q} \oplus
H\big)=GSp_{A_{\mathbb{E}}}\binom{1\ \ \ \ 0}{0\ \ -1}_{\rho,Q}$. The spinor
norm of an element of $\widetilde{GSp}_{A_{\mathbb{E}}}\binom{1\ \ \ \ 0}{0\ \
-1}$ was seen in Theorem \ref{dim8id1} to be the image of its multiplier $m$, so
that the same assertion holds for the images of $GSp_{A_{\mathbb{E}}}\binom{1\ \
\ \ 0}{0\ \ -1}_{\rho,Q}$ in $SO\big((A_{\mathbb{E}}^{-})_{\rho,Q} \oplus
H\big)$. It follows that $SO^{1}\big((A_{\mathbb{E}}^{-})_{\rho,Q} \oplus
H\big)$ consists of images of those $g \in GSp_{A_{\mathbb{E}}}\binom{1\ \ \ \
0}{0\ \ -1}_{\rho,Q}$ such that $m(g)\in(\mathbb{F}^{\times})^{2}$. Dividing by
scalars, we restrict attention to those $g$ with $m(g)=1$, i.e., to $g \in
Sp_{A_{\mathbb{E}}}\binom{1\ \ \ \ 0}{0\ \ -1}_{\rho,Q}$, and as the kernel of
the restriction to $Sp_{A_{\mathbb{E}}}\binom{1\ \ \ \ 0}{0\ \ -1}_{\rho,Q}$ is
$\{\pm1\}$, the latter group is indeed the
$spin\big((A_{\mathbb{E}}^{-})_{\rho,Q} \oplus H\big)$. This completes the proof
of the theorem.
\end{proof}

\smallskip

When we wish to consider the independence of our groups of the choices which we
made, we first observe that by the Witt Cancelation Theorem, the complement of
any hyperbolic plane inside an isotropic space is independent (up to
isomorphism) of the specific hyperbolic plane we took. Hence it remains to
see what happens when we change the choices for the complement
$(A_{\mathbb{E}}^{-})_{\rho,Q}$. For this  we again denote by $\sigma$ the map
on $A_{\mathbb{E}}$ which was previously denoted $\rho$, and let $Q$ and
$\sigma$ vary. The resulting independence assertion appears in
\begin{prop}
Let $\sigma$, $\tau$, and $\eta$ be ring automorphisms of $A_{\mathbb{E}}$, all
of order 2, which restrict to $\rho$ on $\mathbb{E}$, and let $Q$, $R$, $S$ and
$T$ be elements of $A_{\mathbb{E}}^{-}$ which satisfy the conditions of parts
$(i)$, $(ii)$, and $(iii)$ of Proposition \ref{NQF2NAE}. Assume further that
the element $T=Qb^{-1}$ of $A_{\mathbb{E}}^{-}$ has vector norm
$\frac{|Q|^{2}}{b\overline{b}}$. Then the four groups
$GSp_{A_{\mathbb{E}}}\binom{1\ \ \ \ 0}{0\ \ -1}_{\sigma,Q}$,
$GSp_{A_{\mathbb{E}}}\binom{1\ \ \ \ 0}{0\ \ -1}_{\sigma,R}$,
$GSp_{A_{\mathbb{E}}}\binom{1\ \ \ \ 0}{0\ \ -1}_{\tau,S}$, and
$GSp_{A_{\mathbb{E}}}\binom{1\ \ \ \ 0}{0\ \ -1}_{\eta,T}$ are isomorphic in
such a way that the isomorphisms take the $Sp$ subgroups to one another.
\label{indepQ}
\end{prop}

\begin{proof}
The proof of Proposition \ref{NQF2NAE} shows that
$A_{\mathbb{E},\sigma,\mathbb{F}Q}^{t^{2}}$ coincides with
$A_{\mathbb{E},\eta,\mathbb{F}S}^{t^{2}}$, and conjugation by $e$ takes
this group to $A_{\mathbb{E},\tau,\mathbb{F}eQ\overline{e}}^{t^{2}}$. Moreover,
this proof yields the existence of two elements $c$ and $d$ of
$A_{\mathbb{E}}^{\times}$, with $c^{\sigma}=c$ and $d^{\tau}=d$, such that $c$
conjugates $A_{\mathbb{E},\sigma,\mathbb{F}Q}^{t^{2}}$ to
$A_{\mathbb{E},\sigma,\mathbb{F}R}^{t^{2}}$ and $d$ conjugates
$A_{\mathbb{E},\tau,\mathbb{F}eQ\overline{e}}^{t^{2}}$ to
$A_{\mathbb{E},\tau,\mathbb{F}S}^{t^{2}}$. We first claim that that the action
of $c$ (resp. $e$) on $A_{\mathbb{E}}^{-}$ sends the space
$(A_{\mathbb{E}}^{-})_{\sigma,Q}$ to $(A_{\mathbb{E}}^{-})_{\sigma,R}$ (resp.
$(A_{\mathbb{E}}^{-})_{\tau,eQ\overline{e}}$). Indeed, for any $v \in
A_{\mathbb{E}}^{-}$ we have $(cv\overline{c})^{\sigma}=cv^{\sigma}\overline{c}$
and $(ev\overline{e})^{\tau}=ev^{\sigma}\overline{e}$ (this is clear for the
$\sigma$-invariant element $c$, and for $e$ we use the fact that
$x^{\tau}=ee^{-\sigma}x^{\sigma}e^{\sigma}e^{-1}$ and the invariance of
$ee^{-\sigma}$ under $x\mapsto\overline{x}^{\sigma}$). Hence if
$v\in(A_{\mathbb{E}}^{-})_{\sigma,Q}$ then we substitute the value of
$v^{\sigma}$ from Lemma \ref{sp6gen}, and then using Lemma \ref{AxA-rel} and
Proposition \ref{NAFg} we get the asserted result (recall that
$R=rcQ\overline{c}$). It follows that the operation of $d$ maps
$(A_{\mathbb{E}}^{-})_{\tau,eQ\overline{e}}$ to
$(A_{\mathbb{E}}^{-})_{\tau,S}$. Next, we show that
$(A_{\mathbb{E}}^{-})_{\eta,T}=(A_{\mathbb{E}}^{-})_{\sigma,Q}$ as well. To see
this, first observe that if $\overline{b}^{\sigma}=b$ and
$b\overline{b}\in\mathbb{E}^{\times}$ then applying $\sigma$ to the latter
scalar yields $\overline{b}b$. As $b$ is invertible (with inverse
$\frac{\overline{b}}{b\overline{b}}$), this shows that the latter scalar lies in
fact in $\mathbb{F}^{\times}$. The fact that both $Q$ and $T=Qb^{-1}$ are in
$A_{\mathbb{E}}^{-}$ allows us now to write $T$ also as
$\frac{bQ}{b\overline{b}}$. Hence if $v\in(A_{\mathbb{E}}^{-})_{\sigma,Q}$ then
using the definition of $\eta$ and Lemma \ref{sp6gen} we find that
$v^{\eta}=-b\overline{b}\frac{T\tilde{v}T}{|Q|^{2}}$, which equals
$-\frac{T\tilde{v}T}{|T|^{2}}$ by our assumption on $|T|^{2}$. This proves that
$v$ indeed lies in $(A_{\mathbb{E}}^{-})_{\eta,T}$. The equality of the two
spaces, as well as the fact that the entire spaces with $R$, $eQ\overline{e}$,
and $S$ are in the image of $(A_{\mathbb{E}}^{-})_{\sigma,Q}$, follows either
by inverting the above argument or by comparing dimensions (and using the
injectivity of all the operations considered here).

Let now $g$ be an element of $GSp_{A_{\mathbb{E}}}\binom{1\ \ \ \ 0}{0\ \
-1}_{\sigma,Q}$, and we present it by using the appropriate parameters in
Corollary \ref{genform}. The previous paragraph shows that by using the same
parameters we also get that $g \in GSp_{A_{\mathbb{E}}}\binom{1\ \ \ \ 0}{0\ \
-1}_{\eta,T}$. For the other groups, observe that when we conjugate the
generators $\binom{1\ \ v}{0\ \ 1}$, $\binom{1\ \ 0}{w\ \ 1}$, and $\binom{a\ \
\ \ 0\ \ \ }{0\ \ m\overline{a}^{-1}}$ of the group
$GSp_{A_{\mathbb{E}}}\binom{1\ \ \ \ 0}{0\ \ -1}$ by some diagonal element of
the form $\binom{x\ \ \ 0\ \ }{0\ \ \overline{x}^{-1}}$ (of multiplier 1), then
$v$, $w$, and $a$ are taken to $xv\overline{x}$, $\overline{x}^{\ -1}wx^{-1}$,
and $xax^{-1}$ respectively. Moreover, if $w=\tilde{u}$ for some $u \in
A_{\mathbb{E}}^{-}$ then its image is $\widetilde{xu\overline{x}}$ (up to the
scalar $N^{A}_{\mathbb{F}}(x)$) by Lemma \ref{AxA-rel}. For the generators of
$GSp_{A_{\mathbb{E}}}\binom{1\ \ \ \ 0}{0\ \ -1}_{\sigma,Q}$, in which
$v\in(A_{\mathbb{E}}^{-})_{\sigma,Q}$,
$w\in(A_{\mathbb{E}}^{-})_{\sigma,\tilde{Q}}$, and $a \in
A_{\mathbb{E},\sigma,\mathbb{F}Q}^{t^{2}}$, we have seen that for $x=c$ (resp.
$x=e$) the images of $v$ and $a$ lie in $(A_{\mathbb{E}}^{-})_{\sigma,R}$ and
$A_{\mathbb{E},\sigma,\mathbb{F}Q}^{t^{2}}$ (resp.
$(A_{\mathbb{E}}^{-})_{\tau,eQ\overline{e}}$ and
$A_{\mathbb{E},\tau,\mathbb{F}eQ\overline{e}}^{t^{2}}$) respectively. Moreover,
$w=\tilde{u}$ for $u\in(A_{\mathbb{E}}^{-})_{\sigma,Q}$ by parts $(i)$ and
$(iv)$ of Lemma \ref{QthetaQrel}, and the fact that both $c$ and $e$ have
reduced norms in $\mathbb{F}$ shows that its image is the $\theta$-image of an
element which lies in $(A_{\mathbb{E}}^{-})_{\sigma,R}$ (resp.
$(A_{\mathbb{E}}^{-})_{\tau,eQ\overline{e}}$). Hence part $(i)$ of Lemma
\ref{QthetaQrel} proves the assertion for $w$ as well. This shows that
conjugation by $\binom{c\ \ \ 0\ \ }{0\ \ \overline{c}^{-1}}$ (resp. $\binom{e\
\ \ 0\ \ }{0\ \ \overline{e}^{-1}}$) sends $GSp_{A_{\mathbb{E}}}\binom{1\ \ \ \
0}{0\ \ -1}_{\sigma,Q}$ to $GSp_{A_{\mathbb{E}}}\binom{1\ \ \ \ 0}{0\ \
-1}_{\sigma,R}$ (resp. $GSp_{A_{\mathbb{E}}}\binom{1\ \ \ \ 0}{0\ \
-1}_{\tau,eQ\overline{e}}$), and further conjugation of the latter group by
$\binom{d\ \ \ 0\ \ }{0\ \ \overline{d}^{-1}}$ takes the latter group to
$GSp_{A_{\mathbb{E}}}\binom{1\ \ \ \ 0}{0\ \ -1}_{\tau,S}$. The fact that these
identifications and conjugations yields the full groups is established either
by inverting the above argument, or by observing that since the maps and
identifications from the previous paragraph are all surjective, we get all the
generators of the required groups in this process. Since our identifications
and isomorphisms commute with the multiplier maps to $\mathbb{F}^{\times}$, the
$Sp$ groups are sent to one another in this process. This completes the proof of
the proposition.
\end{proof}
Note that as both $\eta$ and $v\mapsto-\frac{T\tilde{v}T}{|T|}^{2}$ are
involutions which separate $A_{\mathbb{E}}^{-}$ to $\pm1$-eigenspaces, the proof
of Proposition \ref{indepQ} already shows that the vector norm of $T=Qb^{-1}$
must be $\pm\frac{|Q|^{2}}{b\overline{b}}$. It seems likely that only the $+$
sign is possible (making the additional assumption in Proposition \ref{indepQ}
redundant), but we have not checked this out in detail. The remarks about the
transitivity of the relations from Proposition \ref{NQF2NAE} on the possible
choices of ring automorphisms of $A_{\mathbb{E}}^{-}$ (commuting with
$\iota_{B}\otimes\iota_{C}$ and reducing to $\rho$ on $\mathbb{E}$ as usual)
extend to similar assertions for Proposition \ref{indepQ}.

Going back to our previous notation, with $\rho$ on $A_{\mathbb{E}}$ as well as
on $M_{2}(A_{\mathbb{E}})$, the spaces with more isotropic vectors which appear
in this case are considered in the following
\begin{cor}
Assume that after extending scalars to $\mathbb{E}$, the quadratic space
$\big((A_{\mathbb{E}}^{-})_{\rho,Q} \oplus
H\big)_{\mathbb{E}}=A_{\mathbb{E}}^{-} \oplus H_{\mathbb{E}}$ splits more than
one hyperbolic plane. Then there exists some quaternion algebra $B$ over
$\mathbb{F}$ and some number $\delta\in\mathbb{F}^{\times}$ representing a class
modulo $N^{B}_{\mathbb{F}}(B^{\times})$ such that the following assertions hold:
The quadratic space is $\mathbb{E} \oplus B \oplus H$ with the norms from
$\mathbb{E}$ multiplied by $-\delta$, and its Gspin and spin groups, denoted
$GSp_{4}(B_{\mathbb{E}})_{\rho,\delta}$ and
$Sp_{4}(B_{\mathbb{E}})_{\rho,\delta}$ respectively, consist of those elements
of $GSp_{4}(B_{\mathbb{E}})$ and $Sp_{4}(B_{\mathbb{E}})$ on which $\psi$ is
defined as in the remark following Corollary \ref{iso8id1}, and conjugating it
by the diagonal matrix with diagonal entries $-\delta$, 1, $-1$, and $\delta$
operates in the same way as $\rho$. When the latter $\mathbb{E}$-space splits
more than two hyperbolic planes, meaning that it is the direct sum of 4
hyperbolic planes, then there exists representatives $\varepsilon$ of a class
modulo $N^{\mathbb{E}}_{\mathbb{F}}(\mathbb{E}^{\times})$ and $\delta$ of a
class with respect to $N^{B}_{\mathbb{F}}(B^{\times})$ for
$B=(\mathbb{E},\rho,\varepsilon)$, for which the following holds: The space is
the direct sum of a hyperbolic plane and three copies of $\mathbb{E}$, with the
norms in two copies are multiplied by $-\varepsilon$ and $-\delta$, and the spin
group is a the subgroup of a double cover of $SO^{1}\binom{0\ \ I}{I\ \ 0}$ over
$\mathbb{E}$ on which $\rho$ coincides with the conjugation of $\psi$ by the
diagonal $8\times8$ matrix whose diagonal entries are $\delta\varepsilon$,
$-\delta$, $-\varepsilon$, 1, $\varepsilon$, $-1$, $-\delta\varepsilon$, and
$\delta$. The Gspin group is a the subgroup of the spinor norm related subgroup
of the general special orthogonal group of $\binom{0\ \ I}{I\ \ 0}$ which is
defined by the same relation between $\rho$ and $\psi$. \label{iso8igen}
\end{cor}

\begin{proof}
The existence of $B$ and $\delta$, as well as $\varepsilon$, is a consequence of
Corollary \ref{iso6gen}, from which we also adopt the choice of $Q$ to be
$\binom{0\ \ \delta}{1\ \ 0}$. The form of the spaces is then given in Corollary
\ref{alt6gen}. The description of the Gspin and spin groups now follows from
Theorem \ref{dim8igen}, Corollary \ref{GSprhoQst}, and Corollary \ref{iso8id1},
taking into consideration the conjugation by $\binom{1\ \ 0}{0\ \ R}$ or
$\binom{1\ \ 0}{0\ \ S}$ (which are $\rho$-invariant and operate on our
$\hat{Q}$ in the same way), the additional conjugation by $\binom{R\ \ 0}{0\ \
R}$ or $\binom{S\ \ 0}{0\ \ S}$ in the definition of $\psi$ on these groups, and
the action of $\rho$ on $(\mathbb{E},\rho,\varepsilon)$ for the latter case.
This proves the corollary. 
\end{proof}
We remark that in the second case considered in Corollary \ref{iso8igen} the
ambient spin group of the direct sum of 4 hyperbolic planes over $\mathbb{E}$
comes, as seen in Corollary \ref{iso8id1}, with three inequivalent maps onto
the associated $SO^{1}$ group. The subgroup considered in Corollary
\ref{iso8igen} in this case is defined by the images in two of these
representations (those which is not the defining representation as a spin group)
being $\rho$-images of one another, after replacing one of them by an
equivalent one. The fact that the double cover splits over this group, but it
still maps to the defining representation with an order 2 kernel, is related to
the fact that $-I$ equals its $\rho$-image. Indeed, of the order 2 elements in
the kernels of the three representations, only $-I$ with $\psi$-image $-I$
satisfies the $\psi_{\hat{Q}}=\rho$ condition. The associated Gspin group maps
onto the full special orthogonal group (with kernel $\mathbb{F}^{\times}$) in
the first representation, but remains injective in the other two
representations, allowing the bilinear form to be multiplied by scalars from
$\mathbb{F}^{\times}$ (the same scalar in both representations). 

\smallskip

Every subspace of dimension 7 is contained in a space as we considered in Lemma
\ref{sp8igen}. Indeed, just choose a non-zero vector norm from the space, and
add a vector whose vector norm is the additive inverse of the chosen one.
However, there is then a relation between the choice of the vector norm and the
discriminant of the resulting 8-dimensional space, meaning that the description
of the groups thus obtained depends on many choices. Hence we content ourselves
with 7-dimensional spaces in which we can make the resulting space have
discriminant 1, yielding a ``canonical'' complement, and do not pursue the rest
(until a good description of non-isotropic 8-dimensional spaces becomes
available). 

\section{Fields with Many Squares \label{ManySq}}

Both the discriminant and the spinor norm take values in the group
$\mathbb{F}^{\times}/(\mathbb{F}^{\times})^{2}$. It is thus worthwhile to
consider explicitly the cases where this group is very small. A field of
characteristic different from 2 for which this group is trivial is called
\emph{quadratically closed}. This is the case, for example, where $\mathbb{F}$
is algebraically closed, e.g., $\mathbb{F}=\mathbb{C}$. Over such a field
$\mathbb{F}$ (with characteristic different from 2), every two non-degenerate
quadratic spaces of the same dimension $n$ are isomorphic, hence the special
orthogonal group can be denoted simply $SO(n,\mathbb{F})$. The group
$SO^{1}(n,\mathbb{F})$ always coincides with $SO(n,\mathbb{F})$, as the spinor
norm is trivial. Note that $\mathbb{F}$ admits neither non-trivial quadratic
field extensions nor non-split quaternion algebras, so that we may always take
$\mathbb{E}=\mathbb{F}\times\mathbb{F}$ and $B$ (or $C$) to be
$M_{2}(\mathbb{F})$. Gathering the results of Theorems \ref{dim12}, \ref{dim3},
\ref{dim4}, \ref{dim6d1}, \ref{dim5}, \ref{dim6gen}, \ref{dim8id1},
\ref{dim7rd}, and \ref{dim8igen} and their corollaries, we establish the
following assertions:

$SO(1,\mathbb{F})=\{1\}$, with $spin(1,\mathbb{F})=\{\pm1\}$ and
$Gspin(1,\mathbb{F})=\mathbb{F}^{\times}$.

$SO(2,\mathbb{F})\cong\mathbb{F}^{\times}$, with $spin(2,\mathbb{F})$ being also
$\mathbb{F}^{\times}$ and
$Gspin(2,\mathbb{F})=\mathbb{F}^{\times}\times\mathbb{F}^{\times}$.

$spin(3,\mathbb{F})=SL_{2}(\mathbb{F})$, $SO(3,\mathbb{F})=PSL_{2}(\mathbb{F})$,
and $Gspin(3,\mathbb{F})=GL_{2}(\mathbb{F})$.

$spin(4,\mathbb{F})=SL_{2}(\mathbb{F}) \times SL_{2}(\mathbb{F})$,
$SO(4,\mathbb{F})$ is the quotient by $\{\pm(I,I)\}$, and $Gspin(4,\mathbb{F})$
is a subgroup $GL_{2}(\mathbb{F}) \times GL_{2}(\mathbb{F})$ determined by the
equal determinant condition.

$spin(5,\mathbb{F})=Sp_{4}(\mathbb{F})$, $SO(5,\mathbb{F})=PSp_{4}(\mathbb{F})$,
and $Gspin(5,\mathbb{F})=GSp_{4}(\mathbb{F})$.

$spin(6,\mathbb{F})=SL_{4}(\mathbb{F})$, $SO(6,\mathbb{F})$ is the quotient by
$\pm I$ (this is not $PSL_{4}(\mathbb{F})$, since in order to obtain the latter
group one must also divide by the two square roots of $-1$), and
$Gspin(6,\mathbb{F})=GL_{4}(\mathbb{F})$ since the $(\mathbb{F}^{\times})^{2}$
condition on the determinant is vacuous over a quadratically closed field.

$spin(7,\mathbb{F})$ is the subgroup of those elements $\binom{a\ \ b}{c\ \ d}$
of $SO\binom{0\ \ I}{I\ \ 0} \subseteq GL_{8}(\mathbb{F})$ in which $ad^{t}$ and
$bc^{t}$ are in $\mathbb{F}$ and square to the determinants of the corresponding
$4\times4$ blocks, and either $bd^{-1}$ and $ca^{-1}$ or $ac^{-1}$ and $db^{-1}$
are anti-symmetric and multiply to give minus their Pfaffian.
$Gspin(7,\mathbb{F})$ is the group defined by the same conditions, but in which
the bilinear form arising from $\binom{0\ \ I}{I\ \ 0}$ may be multiplied by a
scalar.

Finally, $spin(8,\mathbb{F})$ has 3 inequivalent representations in which it
maps onto the $SO\binom{0\ \ I}{I\ \ 0}$ with different order 2 kernels.
$Gspin(8,\mathbb{F})$ maps to $SO\binom{0\ \ I}{I\ \ 0}$ with kernel
$\mathbb{F}^{\times}$, but extending the two other representations of
$spin(8,\mathbb{F})$ to it results in transformations which may multiply the
bilinear form defined by $\binom{0\ \ I}{I\ \ 0}$ by non-trivial scalars.

\smallskip

We now present the case where $\mathbb{F}^{\times}/(\mathbb{F}^{\times})^{2}$
has order 2 (and $ch\mathbb{F}\neq2$). In this case there is a unique quadratic
field extension of $\mathbb{F}$, which we denote $\mathbb{E}$. Hence all the
(non-trivial) unitary groups defined over $\mathbb{F}$ are based on $\mathbb{E}$
with its Galois automorphism $\rho$ over $\mathbb{F}$. One family of fields
having this property is the family of finite fields of odd cardinality. As
another example, recall that a field $\mathbb{F}$ is \emph{Euclidean} if it is
ordered and every positive element is a square. For these fields we have
\begin{prop}
If $\mathbb{F}$ is Euclidean then $\mathbb{E}$ equals $\mathbb{F}(\sqrt{-1})$,
it  is quadratically closed, and we have
$N^{\mathbb{E}}_{\mathbb{F}}(\mathbb{E}^{\times})=(\mathbb{F}^{\times})^{2}$.
Moreover, the quaternion algebra $\big(\frac{-1,-1}{\mathbb{F}}\big)$, which we
denote $\mathbb{H}$ from now on, is not split and satisfies
$N^{\mathbb{H}}_{\mathbb{F}}(\mathbb{H}^{\times})=(\mathbb{F}^{\times})^{2}$ as
well. \label{Euc}
\end{prop}

\begin{proof}
We know that $\mathbb{E}=\mathbb{F}(\sqrt{-1})$ since $-1$ cannot be a square
in an ordered field. Thus, for any $z\in\mathbb{E}$,
$N^{\mathbb{E}}_{\mathbb{F}}(z)$ may be presented as the sum of two squares in
$\mathbb{E}$, not both zero if $z\neq0$. The assertion
$N^{\mathbb{E}}_{\mathbb{F}}(\mathbb{E}^{\times})=(\mathbb{F}^{\times})^{2}$ now
follows from the properties of orderings. A similar argument shows that
$N^{\mathbb{H}}_{\mathbb{F}}(\alpha)$ is the sum of four squares, hence
non-zero for $\alpha\neq0$. Hence $\mathbb{H}$ is a division algebra and
$N^{\mathbb{H}}_{\mathbb{F}}(\mathbb{H}^{\times})=(\mathbb{F}^{\times})^{2}$.

It remains to show that $\mathbb{E}$ is quadratically closed. Any
$z\in\mathbb{E}^{\times}$ can be written as $r^{2}u$ with
$r\in\mathbb{F}^{\times}$ and $u\in\mathbb{E}^{1}$: As
$N^{\mathbb{E}}_{\mathbb{F}}(z)\in(\mathbb{F}^{\times})^{2}$, it is the square
of some $t\in\mathbb{F}^{\times}$, and by replacing $t$ by $-t$ is necessary,
we may assume $t>0$. But then $t=r^{2}$, and
$u=\frac{z}{r^{2}}\in\mathbb{E}^{1}$ as desired. But then Hilbert's Theorem 90
implies that $u=\frac{w^{\rho}}{w}$ for $w\in\mathbb{E}^{\times}$, and as this
element is the quotient of
$N^{\mathbb{E}}_{\mathbb{F}}(w)\in(\mathbb{F}^{\times})^{2}$ and
$w^{2}\in(\mathbb{E}^{\times})^{2}$, we find that $u$ (hence also $z$) is a
square in $\mathbb{E}$. Thus $\mathbb{E}$ is indeed quadratically closed, which
completes the proof of the proposition.
\end{proof}

Note that $\mathbb{H}$ and $M_{2}(\mathbb{F})$ are the only quaternion algebras
over $\mathbb{F}$: For any symbol $\big(\frac{\alpha,\beta}{\mathbb{F}}\big)$
we may get an isomorphic quaternion algebra with $\alpha$ and $\beta$ taken
from $\{\pm1\}$, yielding $\mathbb{H}$ if $\alpha=\beta=-1$ and a split algebra
otherwise. Thus, there are also two bi-quaternion algebras over $\mathbb{F}$,
namely $M_{2}(\mathbb{H})$ and $M_{4}(\mathbb{F})$. All these split over the
quadratically closed extension $\mathbb{E}$. We also remark that the
Artin--Schreier Theorem states that every field $\mathbb{F}$ such that the
algebraic closure of $\mathbb{F}$ has non-trivial finite degree over
$\mathbb{F}$ (like $\mathbb{R}$) must be Euclidean.

For a non-degenerate quadratic form over a Euclidean field one defines the
\emph{signature}: An orthogonal basis can be normalized to have norms in $\pm1$,
and then the signature is the number of $+1$s and the number of $-1$s. These are
two numbers $(p,q)$ which sum to the dimension of the space, and we have
\begin{prop}
The signature classifies quadratic spaces over $\mathbb{F}$ up to isometry.
\label{sign}
\end{prop}

\begin{proof}
We can distinguish the spaces of signature $(n,0)$ and $(0,n)$ (these spaces
are called \emph{definite}, the others \emph{indefinite}) from the others by the
fact that the first one has only positive vector norms, the second has only
negative vector norms, and all the rest are isotropic. This completes the
verification for dimensions 1 and 2. For a larger dimension, it suffices to
prove that isotropic spaces with different signatures are not isometric. But
each such space splits a hyperbolic plane, and the complements must have
different signatures since a hyperbolic plane has signature $(1,1)$. As these
complements are not isometric by the induction hypothesis, the original spaces
cannot be isometric by the Witt Cancelation Theorem. This proves the
proposition.
\end{proof}

\smallskip

For a finite field $\mathbb{F}$ of odd cardinality,
$\mathbb{F}^{\times}/(\mathbb{F}^{\times})^{2}$ has order 2, but it is not
Euclidean as it cannot be ordered (only fields of characteristic 0 may admit
orderings). In fact, we have the following complement to Proposition \ref{Euc}:
\begin{prop}
Let $\mathbb{F}$ be a field such that
$\mathbb{F}^{\times}/(\mathbb{F}^{\times})^{2}$ has order 2. Then if
$\mathbb{F}$ is not Euclidean then
$N^{\mathbb{E}}_{\mathbb{F}}(\mathbb{E}^{\times})=\mathbb{F}^{\times}$,
$\mathbb{E}$ is not quadratically closed, there are no non-split quaternion
algebras over $\mathbb{F}$, and every quadratic form over $\mathbb{F}$ is
determined by its dimension and discriminant.
\label{F/F2=2}
\end{prop}

\begin{proof}
The norm group $N^{\mathbb{E}}_{\mathbb{F}}(\mathbb{E}^{\times})$ is a subgroup
of $\mathbb{F}^{\times}$ which contains $(\mathbb{F}^{\times})^{2}$. Hence if
the latter has index 2 in the former,
$N^{\mathbb{E}}_{\mathbb{F}}(\mathbb{E}^{\times})$ must coincide with one of
these two groups. Now, the only quaternion algebra which may not split is
$\mathbb{H}=\big(\frac{d,d}{\mathbb{F}}\big)$ where $d\in\mathbb{F}^{\times}$ is
not a square. This algebra does split if and only if $d \in
N^{\mathbb{E}}_{\mathbb{F}}(\mathbb{E}^{\times})$, which is equivalent to
$N^{\mathbb{E}}_{\mathbb{F}}(\mathbb{E}^{\times})$ being the full group
$\mathbb{F}^{\times}$ by what we have seen above. $\mathbb{E}$ cannot be
quadratically closed in the this case, since the norm of a square in
$\mathbb{E}$ is a square in $\mathbb{F}$.

We claim that if $\mathbb{H}$ does not split then $\mathbb{F}$ must be
Euclidean. First observe that the product of the two generators of $\mathbb{H}$
squares to $-d^{2}$, so that if $\mathbb{H}$ does not split then
$-1\not\in(\mathbb{F}^{\times})^{2}$, and we may take $d=-1$. But then norms
from $\mathbb{E}$ are sums of two squares. It thus follows from the assumption
$N^{\mathbb{E}}_{\mathbb{F}}(\mathbb{E}^{\times})=(\mathbb{F}^{\times})^{2}$
that $(\mathbb{F}^{\times})^{2}$ is closed under addition, and as
$\mathbb{F}^{\times}=(\mathbb{F}^{\times})^{2}\times\{\pm1\}$ we find that
$(\mathbb{F}^{\times})^{2}$ defines an ordering on $\mathbb{F}$. This proves
that $\mathbb{F}$ is indeed Euclidean.

It remains to prove that in the non-Euclidean case, a quadratic form is
determined by its dimension and discriminant. In dimension 1 the discriminant
characterizes the quadratic space over any field. In dimension 2 it follows from
Lemma \ref{sp2} that any space contains vectors of any given non-zero vector
norm: For trivial discriminant this is clear (a hyperbolic plane), and for the
non-trivial discriminant this follows from our assumption on
$N^{\mathbb{E}}_{\mathbb{F}}(\mathbb{E}^{\times})$. The assertion for dimension
2 is now clear, by taking a vector of vector norm 1 and knowing what the
orthogonal complement must be. Every quadratic space of dimension 3 (hence also
larger) must therefore be isotropic: Fix an anisotropic vector $v$, the space
$v^{\perp}$ must contain some $u$ with $|u|^{2}=-|v|^{2}$, and then $u+v$ is
isotropic. The assertion now follows by induction: If two quadratic spaces have
the same dimension at least 3 and the same discriminant, then both are
isotropic, both split hyperbolic planes, and in both the complement has the same
dimension and the same discriminant. As the complements are isometric by the
induction hypothesis, the same assertion holds for the original ones. This
completes the proof of the proposition.
\end{proof}

In addition to finite fields, every \emph{quasi-finite} field (i.e., a perfect
field which admits a unique extension of every finite order) of characteristic
different from 2 satisfies the conditions of Proposition \ref{F/F2=2}, and is
not Euclidean. The same holds for (perfect) fields whose absolute Galois group
misses some factors $\mathbb{Z}_{p}$ for odd $p$ in the pro-finite completion
of $\mathbb{Z}$, such as direct limits of finite fields where the power of 2 in
the exponent is bounded (otherwise the result is quadratically closed). Hence we
call a non-Euclidean field of characteristic different from 2 such that
$\mathbb{F}^{\times}/(\mathbb{F}^{\times})^{2}$ has order 2 \emph{quadratically
finite}.

In Proposition \ref{Euc} we have seen that the quadratic extension of a
Euclidean field is quadratically closed. As a complementary claim, Proposition
\ref{F/F2=2} has the following 
\begin{cor}
If $\mathbb{F}$ is quadratically finite then so is its (unique) quadratic
extension $\mathbb{E}$. \label{quadfin}
\end{cor}

\begin{proof}
First observe that if $z\in\mathbb{E}$ satisfies
$N^{\mathbb{E}}_{\mathbb{F}}(z)\in(\mathbb{F}^{\times})^{2}$ then an argument
similar to the last part of the proof of Proposition \ref{Euc} shows that $z$
can be presented as $\frac{tN^{\mathbb{E}}_{\mathbb{F}}(w)}{w^{2}}$ for some
$t\in\mathbb{F}^{\times}$ and $w\in\mathbb{F}^{\times}$. It follows that
$z\in(\mathbb{E}^{\times})^{2}$ since we clearly have
$\mathbb{F}^{\times}\subseteq(\mathbb{E}^{\times})^{2}$. Thus the index of
$(\mathbb{E}^{\times})^{2}$ in $\mathbb{E}^{\times}$ can be at most 2, but it
has to be precisely 2 since Proposition \ref{F/F2=2} shows that
$\mathbb{E}^{\times}/(\mathbb{E}^{\times})^{2}$ cannot be of order 1. Since
$-1\in\mathbb{F}^{\times}$ is a square in $\mathbb{E}$, the latter field cannot
be Euclidean, hence it is quadratically finite by Proposition \ref{F/F2=2}. This
proves the corollary. 
\end{proof}

For notational purposes, it will be convenient to identify the group
$\mathbb{F}^{\times}/(\mathbb{F}^{\times})^{2}$, in both the Euclidean and
quadratically finite cases,  with $\{\pm1\}$ (this is the Legendre symbol over
$q$ in the finite case). As the spinor norm takes values in this group, we
write $SO^{+}$ for the $SO^{1}$ groups.

We know that $SO^{+}(V)=SO(V)$ wherever $V$ has dimension 1 (and this
group is trivial). But this is almost the only case where this may happen:
\begin{prop}
Let $V$ be a quadratic space of dimension $>1$ over a field $\mathbb{F}$ which
is either Euclidean or quadratically finite. Then $SO^{+}(V)$ has index 2 in
$SO(V)$ unless $\mathbb{F}$ is Euclidean and the space is definite, a case in
which the index is 1. \label{SO+ind2} \end{prop}

\begin{proof}
As $SO^{+}(V)$ is the kernel of a map $SO(V)\mapsto\{\pm1\}$, it either
coincides with $SO(V)$ (if the map is trivial) or has index 2 in it. Therefore
it suffices to construct an element having non-trivial spinor norm, and prove
that there is no such element in the exceptional case. Now, if $\mathbb{F}$ is
quadratically finite then the proof of Proposition \ref{F/F2=2} shows that
there $V$ contains vectors with vector norms in $(\mathbb{F}^{\times})^{2}$ as
well as anisotropic vectors whose vector norms do not lie in
$(\mathbb{F}^{\times})^{2}$. The same assertion clearly holds for indefinite
spaces in the Euclidean case. Hence the composition of reflections in one
vector of each vector norm yields the desired element of $SO(V)$. On the other
hand, if $\mathbb{F}$ is Euclidean and the space is definite then every
reflection has the same vector norm. As Proposition \ref{CDT} implies the
elements of $SO(V)$ are products of an even number of reflections, the
triviality of the spinor norm in this case follows. This proves the proposition.
\end{proof}

Using Proposition \ref{sign}, all the special orthogonal groups over a
Euclidean field take the form $SO(p,q,\mathbb{F})$ for some natural numbers $p$
and $q$, where we have $SO(p,q,\mathbb{F})=SO(q,p,\mathbb{F})$ by a global sign
inversion on the space. By Proposition \ref{SO+ind2}, the subgroup
$SO^{+}(p,q,\mathbb{F})$ has index 2 in $SO(p,q,\mathbb{F})$ unless $pq=0$. On
the other hand, it follows from Proposition \ref{F/F2=2} that a special
orthogonal group over a quadratically finite field takes the form
$SO(n,\varepsilon,\mathbb{F})$, where $n$ is the dimension and
$\varepsilon\in\{\pm\}$ represents the discriminant. Moreover, as rescaling may
change the discriminant in odd dimensions, we write just $SO(n,\mathbb{F})$ for
odd $n$. Proposition \ref{SO+ind2} shows that $SO^{+}(n,\varepsilon,\mathbb{F})$
or $SO^{+}(n,\mathbb{F})$ always has index 2 there if $n>1$ (and otherwise the
groups are trivial). The spin and Gspin groups are denoted with $SO$ or $SO^{+}$
replaced by $spin$ and $Gspin$ respectively. For the finite fields
$\mathbb{F}_{q}$, where we may replace $\mathbb{F}_{q}$ by simply $q$ in the
notation, this means that $spin(n,\varepsilon,q)$ and $SO(n,\varepsilon,q)$ have
the same cardinality (for $n>1$), while the cardinality of
$SO^{+}(n,\varepsilon,q)$ is half that number and to get the cardinality of
$Gspin(n,\varepsilon,q)$ we multiply by $q-1$.

Recall that in some settings we write results in terms of unitary groups.
Arguments which are similar to the orthogonal groups over quadratically closed
and Euclidean fields show that unitary spaces over quadratically finite fields
are determined by their dimensions, while over Euclidean field they have
signatures just like the quadratic ones. Hence we shall use notations like
$U_{\mathbb{E},\rho}(n)$ in the former case and $U_{\mathbb{E},\rho}(p,q)$ in
the latter case, as well as $U$ replaced by $SU$, $GU$, or $GSU$ when required.
The groups $Sp_{\mathbb{H}}(p,q)$ are defined in a similar manner in the
Euclidean case.

\smallskip

When we apply Theorems \ref{dim12}, \ref{dim3}, \ref{dim4}, \ref{dim6d1},
\ref{dim5}, \ref{dim6gen}, \ref{dim8id1}, \ref{dim7rd}, and \ref{dim8igen}, and
their corollaries, to  the quadratically finite field case, the results we get
are as follows:

$SO(1,\mathbb{F})$ as well as $SO^{+}(1,\mathbb{F})$ are $\{1\}$, while
$spin(1,\mathbb{F})$ equals $\{\pm1\}$ and $Gspin(1,\mathbb{F})$ is
$\mathbb{F}^{\times}$.

$SO(2,+,\mathbb{F})$ and $spin(2,+,\mathbb{F})$ are both isomorphic to
$\mathbb{F}^{\times}$, $SO^{+}(2,+,\mathbb{F})$ is $(\mathbb{F}^{\times})^{2}$,
and $Gspin(2,+,\mathbb{F})=\mathbb{F}^{\times}\times\mathbb{F}^{\times}$. On
the other hand, both $SO(2,-,\mathbb{F})$ and $Spin(2,-,\mathbb{F})$ are
isomorphic to $\mathbb{E}^{1}$ (or equivalently $U_{\mathbb{E},\rho}(1)$),
$SO^{+}(2,-,\mathbb{F})$ is $(\mathbb{E}^{1})^{2}$, and $Gspin(2,-,\mathbb{F})$
equals $\mathbb{E}^{\times}$ (or $GU_{\mathbb{E},\rho}(1)$).

$spin(3,\mathbb{F})=SL_{2}(\mathbb{F})=SU_{\mathbb{E},\rho}(2)$,
$Gspin(3,\mathbb{F})=GL_{2}(\mathbb{F})=GSU_{\mathbb{E},\rho}(2)$,
$SO^{+}(3,\mathbb{F})=PSL_{2}(\mathbb{F})$,
and $SO(3,\mathbb{F})=PGL_{2}(\mathbb{F})$.

$spin(4,+,\mathbb{F})$ is $SL_{2}(\mathbb{F}) \times SL_{2}(\mathbb{F})$,
$Gspin(4,+,\mathbb{F})$ is a subgroup of the product $GL_{2}(\mathbb{F}) \times
GL_{2}(\mathbb{F})$ determined by the equal determinant condition, and
$SO^{+}(4,+,\mathbb{F})$ and $SO(4,+,\mathbb{F})$ are the appropriate quotients.
For the non-trivial discriminant, we get $SL_{2}(\mathbb{E})$ for
$spin(4,-,\mathbb{F})$ and $PSL_{2}(\mathbb{E})$ for $SO^{+}(4,-,\mathbb{F})$,
while $SO(4,-,\mathbb{F})$ is obtained from the latter group as a direct
product with $\{\pm1\}$ and
$Gspin(4,-,\mathbb{F})=GL_{2}^{\mathbb{F}^{\times}}(\mathbb{E})$.

The group $spin(5,\mathbb{F})$ is $Sp_{4}(\mathbb{F})$,
$SO^{+}(5,\mathbb{F})$ is $PSp_{4}(\mathbb{F})$, $Gspin(5,\mathbb{F})$ equals
$GSp_{4}(\mathbb{F})$, and $SO(5,\mathbb{F})$ is $PGSp_{4}(\mathbb{F})$.

$spin(6,+,\mathbb{F})=SL_{4}(\mathbb{F})$, $SO^{+}(6,+,\mathbb{F})$ equals
$PSL_{4}(\mathbb{F})$ if $-1\not\in(\mathbb{F}^{\times})^{2}$ (but not
otherwise!), $SO(6,+,\mathbb{F})$ is the direct product with $\{\pm1\}$, and
$Gspin(6,+,\mathbb{F})$ equals $GL_{4}^{(\mathbb{F}^{\times})^{2}}(\mathbb{F})$.
With the other discriminant we get the groups $SU_{\mathbb{E},\rho}(4)$ for
$spin(6,-,\mathbb{F})$ and $GSU_{\mathbb{E},\rho}(4)$ for
$Gspin(6,-,\mathbb{F})$, with the groups $SO^{+}(6,-,\mathbb{F})$ and
$SO(6,-,\mathbb{F})$ being the appropriate quotients.

$spin(7,\mathbb{F})$ may again be described as the subgroup of $SO^{+}\binom{0\
\ I}{I\ \ 0} \subseteq GL_{8}(\mathbb{F})$ in which the elements, presented as
block matrices $\binom{a\ \ b}{c\ \ d}$, satisfy the conditions that $ad^{t}$
and $bc^{t}$ are scalars (with squares some block determinants) and $bd^{-1}$
and $ca^{-1}$, or $ac^{-1}$ and $db^{-1}$, are anti-symmetric and related to
one another via the Pfaffian being their product. For $Gspin(7,\mathbb{F})$ we
relax the $SO\binom{0\ \ I}{I\ \ 0}$ condition to allow scalar multiplication
of the underlying bilinear form.

The group $spin(8,+,\mathbb{F})$ admits 3 inequivalent representations as
double covers of $SO^{+}\binom{0\ \ I}{I\ \ 0}$ groups. These representations
are restrictions of representations of $Gspin(8,+,\mathbb{F})$, in one of which
the kernel becomes $\mathbb{F}^{\times}$ and in the other two the bilinear form
may be multiplied by scalars. $spin(8,-,\mathbb{F})$ is the subgroup of
$spin(8,+,\mathbb{E})$ (which has the structure we just considered by Corollary
\ref{quadfin}, whence also the notation) in which two of the representations to
$SO^{+}\binom{0\ \ I}{I\ \ 0}$ (over $\mathbb{E}$) become isomorphic, with
$\rho$ and conjugating by $\binom{I\ \ \ \ 0}{0\ \ -I}$ being the isomorphism.
$Gspin(8,-,\mathbb{F})$ is defined by the same condition on
$Gspin(8,+,\mathbb{E})$, with the two representations which become isomorphic
being those in which the bilinear form is multiplied by scalars (which are only
from $\mathbb{F}^{\times}$ in this subgroup). 

\smallskip

When we consider the case of Euclidean fields, recall that
$(\mathbb{F}^{\times})^{2}$ is the set of positive elements. Hence we denote it
$\mathbb{F}^{\times}_{+}$, and furthermore replace any such superscript by
simply $+$. Recall that the determinant of a space of signature $(p,q)$ is
$(-1)^{q}$, but for the discriminant, which matters to us only for even $n$, we
must multiply by $(-1)^{n/2}$. Note that double covers which are based on
choosing a square root split here, since we have the canonical choice of the
positive square root. The results of Theorems \ref{dim12}, \ref{dim3},
\ref{dim4}, \ref{dim6d1}, \ref{dim5}, \ref{dim6gen}, \ref{dim8id1},
\ref{dim7rd}, and \ref{dim8igen} then take the form appearing below:

$SO(1,0,\mathbb{F})=SO^{+}(1,0,\mathbb{F})=\{1\}$,
$Spin(1,0,\mathbb{F})=\{\pm1\}$, and
$GSpin(1,0,\mathbb{F})$ equals $\mathbb{F}^{\times}$.

$SO(2,0,\mathbb{F})=SO^{+}(2,0,\mathbb{F})$ as well as $spin(2,0,\mathbb{F})$
are $\mathbb{E}^{1}=U_{\mathbb{E},\rho}(1,0)$ and
$GSpin(2,0,\mathbb{F})=\mathbb{E}^{\times}=GU_{\mathbb{E},\rho}(1,0)$. On the
other hand, $SO(1,1,\mathbb{F})=\mathbb{F}^{\times}$,
$SO^{+}(1,1,\mathbb{F})=\mathbb{F}^{\times}_{+}$,
$Spin(1,1,\mathbb{F})=\mathbb{F}^{\times}$ as well, and
$GSpin(1,1,\mathbb{F})=\mathbb{F}^{\times}\times\mathbb{F}^{\times}$.

$spin(3,0,\mathbb{F})$ is $\mathbb{H}^{1}$ or equivalently
$SU_{\mathbb{E},\rho}(2,0)$, $SO(3,0,\mathbb{F})=SO^{+}(3,0,\mathbb{F})$ is
obtained as the quotient by $\{\pm1\}$, and
$Gspin(3,0,\mathbb{F})=\mathbb{H}^{\times}=GSU_{\mathbb{E},\rho}(2,0)$. On the
other hand, $spin(1,2,\mathbb{F})$ is $SL_{2}(\mathbb{F})$ (or equivalently
$SU_{\mathbb{E},\rho}(1,1)$) hence $SO^{+}(1,2,\mathbb{F})=PSL_{2}(\mathbb{F})$,
$GSpin(1,2,\mathbb{F})$ equals $GL_{2}(\mathbb{F})$ (or
$GSU_{\mathbb{E},\rho}(1,1)$), and $SO^{+}(1,2,\mathbb{F})$ is
$PGL_{2}(\mathbb{F})$.

The group $spin(4,0,\mathbb{F})$ is $\mathbb{H}^{1}\times\mathbb{H}^{1}$ or
its isomorph $SU_{\mathbb{E},\rho}(2,0) \times SU_{\mathbb{E},\rho}(2,0)$,
$Gspin(4,0,\mathbb{F})$ is the subgroup of
$\mathbb{H}^{\times}\times\mathbb{H}^{\times}$ consisting of pairs of
quaternions with the same norm, and $SO^{+}(4,0,\mathbb{F})=SO(4,0,\mathbb{F})$
is the corresponding quotient. Similarly, but with the split algebra,
$spin(2,2,\mathbb{F})$ is $SL_{2}(\mathbb{F}) \times SL_{2}(\mathbb{F})$ (or
equivalently $SU_{\mathbb{E},\rho}(1,1) \times SU_{\mathbb{E},\rho}(1,1)$),
$GSpin(2,2,\mathbb{F})$ is the ``same determinant subgroup'' of
$GL_{2}(\mathbb{F}) \times GL_{2}(\mathbb{F})$, and $SO^{+}(2,2,\mathbb{F})$
and $SO(2,2,\mathbb{F})$ are the appropriate quotients. On the other hand,
$spin(1,3,\mathbb{F})=SL_{2}(\mathbb{E})$, $SO^{+}(1,3,\mathbb{F})$ is
$PSL_{2}(\mathbb{E})$, $SO^{+}(1,3,\mathbb{F})$ is the direct product of the
latter group with $\{\pm1\}$, and $GSpin(1,3,\mathbb{F})$ equals
$GL_{2}^{\mathbb{F}^{\times}}(\mathbb{E})$.

We also have $spin(5,0,\mathbb{F})=Sp_{\mathbb{H}}(2,0)$ and
$Gspin(5,0,\mathbb{F})=GSp_{\mathbb{H}}(2,0)$, with
$SO^{+}(5,0,\mathbb{F})=SO(5,0,\mathbb{F})$ being the quotient
$GSp_{\mathbb{H}}(2,0)$ of the former modulo $\{\pm1\}$ or the latter modulo
$\mathbb{F}^{\times}$. In a similar manner, $spin(4,1,\mathbb{F})$ is
$Sp_{\mathbb{H}}(1,1)$, $Gspin(4,1,\mathbb{F})$ is $GSp_{\mathbb{H}}(1,1)$, and
$SO^{+}(4,1,\mathbb{F})$ and $SO(4,1,\mathbb{F})$ are the usual quotients. In
addition, $Sp_{4}(\mathbb{F})$ is $spin(2,3,\mathbb{F})$ so that
$SO^{+}(2,3,\mathbb{F})$ is $PSp_{4}(\mathbb{F})$, and $GSpin(2,3,\mathbb{F})$
equals $GSp_{4}(\mathbb{F})$.

The group $spin(5,1,\mathbb{F})$ is $GL_{2}^{1}(\mathbb{H})$,
$GSpin(5,1,\mathbb{F})$ equals $GL_{2}(\mathbb{H})\times\{\pm1\}$ (the double
cover splits, and the superscript $+$ is unnecessary by Lemma \ref{NM2B} and the
fact that
$N^{\mathbb{H}}_{\mathbb{F}}(\mathbb{H}^{\times})=(\mathbb{F}^{\times})^{2}$),
and $SO^{+}(5,1,\mathbb{F})$ and $SO(5,1,\mathbb{F})$ are the quotients (the
latter being the direct product of the former with $\{\pm1\}$). Using the split
algebra, $Spin(3,3,\mathbb{F})$ is just $SL_{4}(\mathbb{F})$,
$SO^{+}(3,3,\mathbb{F})$ is $PSL_{4}(\mathbb{F})$ as
$-1\not\in(\mathbb{F}^{\times})^{2}$, $GSpin(3,3,\mathbb{F})$ is isomorphic to
$GL_{4}^{+}(\mathbb{F})\times\{\pm1\}$ (a split double cover), and
$SO(3,3,\mathbb{F})$ equals $PSL_{4}(\mathbb{F})\times\{\pm1\}$.
$spin(6,0,\mathbb{F})$ equals $SU_{\mathbb{E},\rho}(4,0)$,
$GSpin(6,0,\mathbb{F})$ is $GSU_{\mathbb{E},\rho}(4,0)$, and
$SO^{+}(6,0,\mathbb{F})=SO(6,0,\mathbb{F})$ is obtained as both the appropriate
quotients. Finally, $spin(4,2,\mathbb{F})$ is isomorphic to
$SU_{\mathbb{E}}(2,2)$, $GSpin(4,2,\mathbb{F})$ is $GSU_{\mathbb{E}}(2,2)$, and
the usual quotients give $SO^{+}(4,2,\mathbb{F})$ and $SO(4,2,\mathbb{F})$.

The group $spin(4,3,\mathbb{F})$ is the subgroup of $SO^{+}\binom{0\ \ I}{I\ \
0}$ ($8\times8$ matrices) consisting of those block matrices $\binom{a\ \ b}{c\
\ d}$ in which $ad^{t}$ and $bc^{t}$ are in $\mathbb{F}$ and square to $\det
a=\det d$ and $\det b=\det c$ respectively, and where either $bd^{-1}$ and
$ca^{-1}$ or $ac^{-1}$ and $db^{-1}$ belong to $M_{4}^{as}(\mathbb{F})$ and
multiply to minus their common Pfaffian. $Gspin(4,3,\mathbb{F})$ is described by
the same condition on the group of matrices in $GL_{8}(\mathbb{F})$ whose action
multiplies $\binom{0\ \ I}{I\ \ 0}$ by a scalar (with some extra condition
extending the $SO^{1}$ condition). For $Gspin(5,2,\mathbb{F})$ we get the
subgroup of $GSp_{4}(\mathbb{H})$, elements $\binom{a\ \ b}{c\ \ d}$ of which
satisfy the conditions that $a\overline{d}^{t}$ and $b\overline{c}^{t}$ are
scalars squaring to $N^{M_{2}(B)}_{\mathbb{F}}(a)=N^{M_{2}(B)}_{\mathbb{F}}(d)$
and $N^{M_{2}(B)}_{\mathbb{F}}(b)=N^{M_{2}(B)}_{\mathbb{F}}(c)$ respectively,
and either the pair $bd^{-1}$ and $ca^{-1}$ or the pair $ac^{-1}$ and $db^{-1}$
is a pair of matrices in $M_{2}^{Her}(B)$ which are minus the adjoints of one
another. $spin(5,2,\mathbb{F})$ is the group of the matrices in 
$Sp_{4}(\mathbb{H})$ having these properties. $Gspin(6,1,\mathbb{F})$ and
$spin(6,1,\mathbb{F})$ are similar subgroups of $GSp_{4}(\mathbb{H})$ and
$Sp_{4}(\mathbb{H})$, but in which the pairs of minus adjoint matrices are
$bd^{-1}$ and $-ca^{-1}$ or $ac^{-1}$ and $-db^{-1}$.

$spin(6,2,\mathbb{F})$ and $Gspin(6,2,\mathbb{F})$ are double covers of
$Sp_{4}(\mathbb{H})$ and $GSp_{4}(\mathbb{H})$ respectively in two, inequivalent
ways. We have omitted the superscript $\mathbb{F}^{2}$ since the reduced norms
from $\mathbb{H}$, hence also from $M_{2}(\mathbb{H})$ by Lemma \ref{NM2B}, are
non-negative, whence the map $\varphi$ from Proposition \ref{GSphom} is trivial
in this case. $spin(4,4,\mathbb{F})$ maps in 3 inequivalent ways to
$SO^{+}\binom{0\ \ I}{I\ \ 0}$ with kernels of order 2, and
$Gspin(4,4,\mathbb{F})$ maps in one representation to $SO\binom{0\ \ I}{I\ \ 0}$
with kernel $\mathbb{F}^{\times}$ but multiplies the bilinear form by arbitrary
scalars in the other two representations (where its kernel remains $\{\pm1\}$).
The group $spin(7,1,\mathbb{F})$ is defined by the condition on
$spin(8,\mathbb{E})$ (recall that $\mathbb{E}$ is quadratically closed by
Proposition \ref{Euc}) which states that conjugating one representations to the
group $SO^{+}\binom{0\ \ I}{I\ \ 0}$ over $\mathbb{E}$ by $\binom{I\ \ \ \ 0}{0\
\ -I}$ yields the $\rho$-image of the other representation. For
$Gspin(7,1,\mathbb{F})$ we apply the same condition on $Gspin(8,\mathbb{E})$
using two representations in which the operation multiplies the bilinear form by
scalars (from $\mathbb{F}^{\times}$ here). The groups $spin(5,3,\mathbb{F})$ and
$spin(5,3,\mathbb{F})$ are obtained in the same manner but with each $4\times4$
identity matrix replaced by $\binom{I\ \ \ \ 0}{0\ \ -I}$ (involving $2\times2$
identity matrices). 

We remark that the case of the spin group in signature $(6,2)$ over
$\mathbb{F}=\mathbb{R}$ was considered in \cite{[SH]}, using Clifford algebras,
Eichler transformations, and some real, complex, and quaternionic analytic
tools. The homomorphism denoted $\phi$ in Lemma 6.10 of that reference is just
$a \mapsto t\overline{a}^{\ -1}$ on
$GL_{2}(\mathbb{H})=GL_{2}^{(\mathbb{F}^{\times})^{2}}(\mathbb{H})$, with the
square roots being positive. Proposition 6.11 there is the projectivization of
our Proposition \ref{GSppsi}. Now, the notion of positive definiteness extends
from $\mathbb{F}=\mathbb{R}$ to any Euclidean field. It thus seems reasonable
that the action of $Sp_{4}(\mathbb{H})$ on the subset of
$M_{2}^{Her}(\mathbb{H})\otimes\mathbb{E}$ in which the ``imaginary part''
(which is also well-defined) is positive definite, as well as the $\psi$-images
of the elements of $\widetilde{Sp}_{4}(\mathbb{H})$ lying over $g \in
Sp_{4}(\mathbb{H})$ being those which send $Z^{t}$ (for $Z$ in the latter space)
to $g(Z)^{t}$, also extend to this more general setting. However, as we have
seen, these aspects of the theory are not required for obtaining the general
result.

We do not get presentations of the definite spin groups $spin(7,0,\mathbb{F})$
and $spin(8,0,\mathbb{F})$ here, since a definite space of dimension 7 does not
represent its discriminant, and a definite space of dimension 8 is not
isotropic, and our methods for spaces of dimensions 7 and 8 require these
properties. 

Unitary groups preserving sesqui-linear forms of dimension 3 do not arise in
the context of orthogonal groups since the dimension 8 of such special unitary
groups is not the dimension of any orthogonal group. Sesqui-linear forms of
signature $(3,1)$ also do not appear here because of the discriminant 1
condition. For $\mathbb{F}=\mathbb{R}$ we can also derive this fact from Lie
theory: The dimension of $SU_{\mathbb{C}}(3,1)$ is indeed 15, but its maximal
compact subgroup $S\big(U(3) \times U(1)\big)$ has dimension 9, which does not
equal the dimension of $SO(p) \times SO(q)$ for any pair $(p,q)$ with sum 6
(the required dimension is $\frac{p(p-1)}{2}+\frac{q(q-1)}{2}$, which attains
only the values 15, 10, 7, and 6).

\smallskip

The splitting of the double covers here comes from the splitting of the sequence
\[1\to\{\pm1\}\to\mathbb{F}^{\times}\to(\mathbb{F}^{\times})^{2}\to1.\] This
happens wherever the Abelian group $\mathbb{F}^{\times}$ contains $\{\pm1\}$ as
a direct summand (e.g., when $\mathbb{F}$ may be ordered, when
$(\mathbb{F}^{\times})^{2}$ is free like in number fields of class number 1
with no complex roots of unity, etc.). Note that when the double cover
$\widetilde{A}^{(\mathbb{F}^{\times})^{2}}$ splits, the elements $(g,|g|^{2})$
with $g \in A^{-} \cap A^{\times}$ will never lie all in the splitting group, as
they generate the full double cover by Theorem \ref{dim6d1}.

On the other hand, for quadratically closed and Euclidean fields the sequence
\[1\to(\mathbb{F}^{\times})^{2}\to\mathbb{F}^{\times}\to\mathbb{F}^{\times}
/(\mathbb{F}^{\times})^{2}\to1\] splits. This fact is related to the full
special orthogonal group admitting a double cover inside the Gspin group: For
the quadratically closed case as well as the definite case this is just the
usual spin group. In the indefinite case the double cover is obtained by
imposing the condition that a certain reduced norm, determinant, or multiplier
takes only the values $\pm1$. In fact, one can show that the only additional
case where this sequence splits is for a quadratically finite field in which
$-1\not\in(\mathbb{F}^{\times})^{2}$ (for the finite field $\mathbb{F}_{q}$ this
happens if and only if $q\equiv3(\mathrm{mod\ }4)$). The double covers in this
case are obtained in a way similar to the indefinite case.

\smallskip

We remark that the cardinality of the group
$\mathbb{F}^{\times}/(\mathbb{F}^{\times})^{2}$ can be either infinite or any
finite power of 2. To see this, observe that if $\mathbb{K}=\mathbb{F}((X))$
then $\mathbb{K}^{\times}/(\mathbb{K}^{\times})^{2}$ is generated by the
injective image of $\mathbb{F}^{\times}/(\mathbb{F}^{\times})^{2}$ and the class
of $X$ (another element of order 2 in
$\mathbb{K}^{\times}/(\mathbb{K}^{\times})^{2}$ which is independent of
$\mathbb{F}^{\times}/(\mathbb{F}^{\times})^{2}$). For the case the group has
order 4 we again have different types of fields. Indeed, if $\mathbb{K}$ is our
field and it takes the form $\mathbb{F}((X))$ then if $\mathbb{F}$ is
quadratically finite then $\mathbb{K}$ admits only one non-split quaternion
algebra (this is also the case for the $p$-adic numbers for an odd prime $p$ or
their finite extensions), while for Euclidean $\mathbb{F}$ there are 3
non-isomorphic quaternion algebras over $\mathbb{K}$ which are not split. The
description of the quadratic spaces over these fields will thus be different,
with the first case probably resembling the results appearing in Section 2 of
Chapter IV of \cite{[S]}. We leave these questions, as well as the question
whether every field $\mathbb{F}$ in which
$\mathbb{F}^{\times}/(\mathbb{F}^{\times})^{2}$ has order 4 resembles one of
these families, for future research.

\noindent\textsc{Fachbereich Mathematik, AG 5, Technische Universit\"{a}t
Darmstadt, Schlossgartenstrasse 7, D-64289, Darmstadt, Germany}

\noindent E-mail address: zemel@mathematik.tu-darmstadt.de

\end{document}